\documentclass[11pt]{article}
\usepackage{amsthm,amsmath,latexsym,amssymb}

\newcommand{\bP}{{\rm |\kern-.15em P}}
\newcommand{\Q}{\kern.3em\rule{.07em}{.65em}\kern-.3em{\rm Q}}
\newcommand{\R}{{\rm I\kern-.15em R}}
\newcommand{\D}{{\rm |\kern-.15em D}}
\newcommand{\h}{{\rm |\kern-.15em H}}
\newcommand{\C}{\kern.3em\rule{.07em}{.65em}\kern-.3em{\rm C}}
\newcommand{\T}{{\rm T\kern-.35em T}}

\theoremstyle{plain}
\newtheorem{theorem}{Theorem}[section]

\newtheorem{proposition}[theorem]{Proposition}
\newtheorem{corollary}[theorem]{Corollary}

\theoremstyle{definition}
\newtheorem{definition}[theorem]{Definition}

\theoremstyle{remark}
\newtheorem{remark}[theorem]{Remark}

\begin{document}
\title{\bf Fractals and the two dimensional Jacobian Conjecture}
\author{Ronen Peretz}
 
\maketitle

\section{Introduction}

In this paper we prove the two dimensional Jacobian Conjecture. Namely, let $F(X,Y)\in\mathbb{C}[X,Y]^2$ satisfy the Jacobian Condition: 
$\det J_F\equiv 1$, where $J_F$ is the Jacobian matrix of $F$. Then $F(X,Y)$ considered as a polynomial mapping 
$\mathbb{C}^2\rightarrow\mathbb{C}^2$ is injective, surjective and $F^{-1}\in\mathbb{C}[X,Y]^2$. Such a polynomial mapping
is called an invertible mapping.

We define the geometrical degree of $F$ by $d_F:=\max\{|F^{-1}(a,b)|\,|\,(a,b)\in\mathbb{C}^2\}$. Here if $A$ is a set then $|A|$
is the cardinality of the set $A$. Then $F$ is invertible if and
only if $d_F=1$. Thus it suffices to prove that if $F$ satisfies the Jacobian condition, then $d_F=1$. Our method accomplishes
that by proving an integral identity of the following form:
\begin{equation}
\int_{\{F\in {\rm et}_P(\mathbb{C}^2)\}}d_F^{-s}d\mu(F)\equiv 1,
\label{eq1}
\end{equation}
for all real positive values of the parameter $s$ in a non-degenerate interval $0\le\alpha<s<\beta$. The measure space 
${\rm et}_P(\mathbb{C}^2)$ is composed of all the counterexamples to the two dimensional Jacobian Conjecture in a
distinguished set that generates all the possible such counterexamples. This generating set of counterexamples is the parallel
of the set of prime integers (in $\mathbb{Z}^+$) and it generates just like the prime numbers generate all the natural numbers 
as described in the Fundamental Theorem of the Arithmetic. This is the reason that we call the \'etale mappings in ${\rm et}_P(\mathbb{C}^2)$, the 
prime \'etale mappings. Clearly $\forall\,F\in {\rm et}_P(\mathbb{C}^2)$ we have $d_F\ge 2$. The measure $\mu$ is a positive measure and so 
the integral above is a strictly decreasing function of the real parameter $s$ unless there are no prime \'etale mappings. It follows that 
there are no polynomial mappings $F$ that satisfy the Jacobian Condition and for which $d_F>1$.

Thus our method works on the full set of possible counterexamples to the two dimensional Jacobian Conjecture. Unlike some
other ideas it does not work on just a single such a mapping and in this sense it is novel. It makes an extensive use of a fractal
structure that can be put on the set of all the counterexamples to the two dimensional Jacobian Conjecture. So the geometry
in our method is not at the level of a single mapping, but at more abstract level of the total set of counterexamples. This set
carries the structure of a very special metric space and it is the geometry of this metric space that is being used in our proof.

To describe this metric space we consider the set of all the polynomial mappings $G(X,Y)=(G_1(X,Y),G_2(X,Y))\in\mathbb{C}[X,Y]^2$ which 
satisfy the following two conditions: \\
1) $\det J_G\equiv 1$, and \\
2) $\deg G_1=\deg_Y G_1$ and $\deg G_2=\deg_Y G_2$. \\
This set is denoted by ${\rm et}(\mathbb{C}^2)$. With composition of mappings, $\circ$, as a binary operation, the pair $({\rm et}(\mathbb{C}^2),\circ)$
is a semigroup. It contains the group $({\rm Aut}(\mathbb{C}^2),\circ)$ of all the invertible polynomial mappings (subject to the
conditions, 1 and 2 above). The geometrical degree function: $d:\,\,{\rm et}(\mathbb{C}^2)\rightarrow\mathbb{Z}^+$, $d(F)=d_F$ is a multiplicative function on
the semigroup $({\rm et}(\mathbb{C}^2),\circ)$, i.e., $d_{F\circ G}=d_F\cdot d_G$. We turn the set ${\rm et}(\mathbb{C}^2)$ into a metric space $({\rm et}(\mathbb{C}^2),\rho_D)$, where the metric $\rho_D$ is determined by a set $D\subseteq\mathbb{C}^2$ with the aid of the following,
\begin{definition}
Let $D$ be an open subset of $\mathbb{C}^2$ with respect to the Euclidean (the strong) topology, that satisfies the following conditions: \\
1) ${\rm int}(\overline{D})=D$ ($D$ has no "slits"). \\
2) $\overline{D}$ is a compact subset of $\mathbb{C}^2$ (in the strong topology). \\
3) $\forall\,F,G\in {\rm et}(\mathbb{C}^2)$, $F(D)=G(D)\Leftrightarrow F=G$. \\
We define the following real valued function:
$$
\rho_D\,:\,{\rm et}(\mathbb{C}^2)\times {\rm et}(\mathbb{C}^2)\rightarrow\mathbb{R}^+\cup\{0\},
$$
$$
\rho_D(F,G)={\rm the\,volume\,\,of}\,\,F(D)\Delta G(D).
$$
Here we use the standard set-theoretic notation of the symmetric difference between two sets $A$ and $B$, i.e.
$A\Delta B=(A-B)\cup (B-A)$.
\end{definition}

\begin{remark}
It is not clear how to construct an open subset $D$ of $\mathbb{C}^2$ that will satisfy the three properties
that are required in definition 1.1. We call such a set $D$ a characteristic set of the family of maps ${\rm et}(\mathbb{C}^2)$.
The reason for this name is the property 3 in definition 1.1. Section 7 of this paper is devoted to a construction
of characteristic sets of rigid families of local diffeomorphic mappings (such as ${\rm et}(\mathbb{C}^2)$ or such as the family of all the entire mappings
$F\,:\,\mathbb{C}^2\rightarrow\mathbb{C}^2$ such that $\det J_F(a,b)\ne 0$, $\forall\,(a,b)\in\mathbb{C}^2$).
\end{remark}

\begin{proposition}
$\rho_D$ is a metric on ${\rm et}(\mathbb{C}^2)$.
\end{proposition}
\noindent
{\bf Proof.} \\
1) $\rho_D(G_1,G_2)=0 \Leftrightarrow {\rm the\,volume\,of}\,\,G_1(D)\Delta G_2(D)=0 \Leftrightarrow G_1(D)=G_2(D)$
(where the last equivalence follows by the fact that $G_1$ and $G_2$ are local homeomorphisms in the strong topology and because
of condition 1 in definition 1.1) $\Leftrightarrow G_1=G_2$ (by condition 3 in definition 1.1). \\
2) By $G_1(D)\Delta G_2(D)=G_2(D)\Delta G_1(D)$ it follows that $\rho_D(G_1,G_2)=\rho_D(G_2,G_1)$. \\
3) Here we use a little technical set-theoretic containment. Namely, for any three sets $A,B$ and $C$ we have, 
$$
A\Delta C\subseteq (A\Delta B)\cup (B\Delta C).
$$
This implies that $G_1(D)\Delta G_3(D)\subseteq (G_1(D)\Delta G_2(D))\cup (G_2(D)\Delta G_3(D))$ from which it follows that
$$
({\rm the\,volume\,of}\,\,G_1(D)\Delta G_3(D))\le ({\rm the\,volume\,of}\,\,G_1(D)\Delta G_2(D))+
$$
$$
+({\rm the\,volume\,of}\,\,G_2(D)\Delta G_3(D)).
$$
Hence the triangle inequality $\rho_D(G_1,G_3)\le\rho_D(G_1,G_2)+\rho_D(G_2,G_3)$ holds. $\qed $ \\
\\
We recall the following standard notions of the composition mappings.

\begin{definition}
Let $F\in {\rm et}(\mathbb{C}^2)$. The right (composition) mapping induced by $F$ is defined by:
$$
R_F\,:\,{\rm et}(\mathbb{C}^2)\rightarrow {\rm et}(\mathbb{C}^2),\,\,\,R_F(G)=G\circ F.
$$
The left (composition) mapping induced by $F$ is defined by:
$$
L_F\,:\,{\rm et}(\mathbb{C}^2)\rightarrow {\rm et}(\mathbb{C}^2),\,\,\,L_F(G)=F\circ G.
$$
\end{definition}

It is natural to inquire when these two mappings are surjective and if they are injective. It turns out that deciding the surjectivity
of these mappings is not hard. So is the fact that $R_F$ is always an injective mapping. However, the question of the injectivity of
$L_F$ is much harder. It turns out that it is always injective. This is a corollary of the fact that $L_F$ is a $\rho_D$ bi-Lipschitz
mapping. The higher bi-Lipschitz constant is $1$. The lower bi-Lipschitz constant can be made as close as we please to $1/d_F$ by taking
the characteristic set $D$ large enough. To be precise - we will see that we can construct $D$ to be almost an open ball $B((a,b),R)$ of radius $R$.
The only difference between $D$ and the ball $B((a,b),R)$ will be in a small portion of the boundary which will have a dense set of spikes. 
Any dilation of $D$ by a factor $t>0$ will also be a characteristic set, $D_t$ of ${\rm et}(\mathbb{C}^2)$ and it will also be very close to the
dilated ball $B((a,b),tR)$. It is clear that $\lim_{t\rightarrow\infty} D_t=\mathbb{C}^2$. What was said above about the lower bi-Lipschitz constant
can now be accurately stated as follows: $\forall\,\epsilon>0$, $\exists\,T=T(\epsilon)$ such that for $t>T$ we have 
$\forall\,G_1,G_2\in {\rm et}(\mathbb{C}^2)$
$$
\left(\frac{1}{d_F}-\epsilon\right)\cdot\rho_{D_t}(G_1,G_2)\le\rho_{D_t}(F\circ G_1,F\circ G_2).
$$
Moreover, we have the following identity:
$$
\lim_{t\rightarrow\infty}\frac{\rho_{D_t}(F\circ G_1,F\circ G_2)}{\rho_{D_t}(G_1,G_2)}=\frac{1}{d_F}.
$$
This means, geometrically that $L_F$ tends to a similarity transformation on ${\rm et}(\mathbb{C}^2)$ with the similarity constant equals to $1/d_F$.
This is a corner stone in putting on ${\rm et}(\mathbb{C}^2)$ an approximate fractal structure with respect to the metric $\rho_{D_t}$ and letting $t\rightarrow\infty$ to make it as close as we want to a classical fractal structure.
The above theorem about the fact that $L_F$ tends to a similarity transformation on ${\rm et}(\mathbb{C}^2)$ with the similarity constant equals to $1/d_F$
is in agreement with the following theorem: $\forall\,F\in {\rm Aut}(\mathbb{C}^2)$ the left (composition) mapping $L_F$ is an isometry of the
metric space $({\rm et}(\mathbb{C}^2),\rho_D)$ (for any characteristic set $D$ of the family of the \'etale mappings ${\rm et}(\mathbb{C}^2)$).

We are now ready to give the approximate fractal structure of the metric space $({\rm et}(\mathbb{C}^2),\rho_D)$. It is composed of the following two
theorems: \\
a) $\forall\,F\in {\rm et}(\mathbb{C}^2)$ the left (composition) mapping $L_F\,:\,{\rm et}(\mathbb{C}^2)\rightarrow L_F({\rm et}(\mathbb{C}^2))$ is an
homeomorphism of metric spaces (the metrics are $\rho_D$) which is approximately the $1/d_F$ similarity mapping (the closer $D$ is to $\mathbb{C}^2$, the better
is this approximation). \\
b) Assuming that ${\rm et}(\mathbb{C}^2)\ne {\rm Aut}(\mathbb{C}^2)$, there exists an infinite index set $I$ and a family of \'etale mappings indexed by
$I$, $\{F_i\,|\,i\in I\}\subset {\rm et}(\mathbb{C}^2)-{\rm Aut}(\mathbb{C}^2)$, such that
$$
{\rm et}(\mathbb{C}^2)={\rm Aut}(\mathbb{C}^2)\cup\bigcup_{i\in I}L_{F_i}({\rm et}(\mathbb{C}^2)),
$$
so that $\forall\,i\in I$, ${\rm Aut}(\mathbb{C}^2)\cap L_{F_i}({\rm et}(\mathbb{C}^2))=\emptyset$, and $\forall\,i,j\in I$, if $i\ne j$ then
$F_i\not\in L_{F_j}({\rm et}(\mathbb{C}^2))$ and $F_j\not\in L_{F_i}({\rm et}(\mathbb{C}^2))$. \\
\\
We can refine the description of the approximate fractal structure of ${\rm et}(\mathbb{C}^2)$ by pointing at a very concrete index set $I$. For that
we borrow the arithmetic multiplicative notion of a prime integer (i.e. a prime number $p\in \mathbb{Z}^+$) to the context of the \'etale semigroup
$({\rm et}(\mathbb{C}^2),\circ)$. Namely, an \'etale mapping $F\in {\rm et}(\mathbb{C}^2)$ is called a composite mapping, 
if $\exists\,G,H\in  {\rm et}(\mathbb{C}^2)-{\rm Aut}(\mathbb{C}^2)$ such that $F=G\circ H$. If $F\in {\rm et}(\mathbb{C}^2)$ is not a composite mapping
then we call it a prime mapping. It follows that if $d_F$ is a prime integer then $F$ is a prime mapping. We denote by ${\rm et}_P(\mathbb{C}^2)$ the set
of all the prime mappings in ${\rm et}(\mathbb{C}^2)$. We have the following theorem (which is parallel to the existence part of the so called The Fundamental Theorem of Arithmetic): If ${\rm et}(\mathbb{C}^2)\ne {\rm Aut}(\mathbb{C}^2)$, then, \\
a) ${\rm et}_P(\mathbb{C}^2)\ne\emptyset$, \\
and \\
b) $\forall\,F\in {\rm et}(\mathbb{C}^2)$, $\exists\,k\in\mathbb{Z}^+\cup\{0\}$, $\exists\,A_0\in {\rm Aut}(\mathbb{C}^2)$, $\exists\,
P_1,\ldots,P_k\in {\rm et}_P(\mathbb{C}^2)$ such that $F=A_0\circ P_1\circ\ldots\circ P_k$. \\
\\
Returning to the approximate fractal representation of ${\rm et}(\mathbb{C}^2)$, we have the following theorem: We can choose the index set
$I={\rm et}_P(\mathbb{C}^2)$. This means that we have the identity,
$$
{\rm et}(\mathbb{C}^2)={\rm Aut}(\mathbb{C}^2)\cup\bigcup_{F\in {\rm et}_P(\mathbb{C}^2)}L_F({\rm et}(\mathbb{C}^2)).
$$
We will see that the following identity holds:
$$
\sum_{F\in {\rm et}_P(\mathbb{C}^2)}\frac{1}{d_F}=1.
$$
The conclusion is the following theorem: \\
a) $|{\rm et}_P(\mathbb{C}^2)|=\aleph_0$, and \\
b) The primes ${\rm et}_P(\mathbb{C}^2)$ form a discrete subset of the metric space $({\rm et}(\mathbb{C}^2),\rho_D)$ (for any characteristic set
of the family of the \'etale mappings, ${\rm et}(\mathbb{C}^2)$). \\
c) The semigroup $({\rm et}(\mathbb{C}^2),\circ)$ is generated by ${\rm Aut}(\mathbb{C}^2)$ and the countable set of primes ${\rm et}_P(\mathbb{C}^2)$. \\

\begin{remark}
One can prove the following theorem: \\
a) $\forall\,F\in {\rm et}(\mathbb{C}^2)$, the right (composition) mapping $R_F\,:\,{\rm et}(\mathbb{C}^2)\rightarrow R_F({\rm et}(\mathbb{C}^2))$ is
an homeomorphism of metric spaces (the metrics are the $\rho_D$). \\
 b) Assuming that ${\rm et}(\mathbb{C}^2)\ne {\rm Aut}(\mathbb{C}^2)$, there exists an infinite index set $I$ and a family of \'etale mappings indexed by
$I$, $\{F_i\,|\,i\in I\}\subset {\rm et}(\mathbb{C}^2)-{\rm Aut}(\mathbb{C}^2)$, such that
$$
{\rm et}(\mathbb{C}^2)={\rm Aut}(\mathbb{C}^2)\cup\bigcup_{i\in I}R_{F_i}({\rm et}(\mathbb{C}^2)),
$$
so that $\forall\,i\in I$, ${\rm Aut}(\mathbb{C}^2)\cap R_{F_i}({\rm et}(\mathbb{C}^2))=\emptyset$, and $\forall\,i,j\in I$, if $i\ne j$ then
$F_i\not\in R_{F_j}({\rm et}(\mathbb{C}^2))$ and $F_j\not\in R_{F_i}({\rm et}(\mathbb{C}^2))$. \\
\\
This resembles very much the theorem on the approximate fractal structure induced on ${\rm et}(\mathbb{C}^2)$ by the left (composition)
mappings $L_F$. Moreover, proving that the $R_F$'s are injective is easy, while proving the injectivity of the $L_F$'s is hard and uses
the fact that $L_F$ is a $\rho_D$ bi-Lipschitz mapping. This raises the question: why do we prefer the left (composition) mappings? The reason is 
exactly that same bi-Lipschitz property, which tends to a $1/d_F$ similarity when $D\rightarrow\mathbb{C}^2$. For these are exactly the properties needed 
to establish the approximate fractal structure on ${\rm et}(\mathbb{C}^2)$. But why don't we try to find another metric on ${\rm et}(\mathbb{C}^2)$
with respect to which the $R_F$'s will tend to similarity mappings? The reason is two folded. First the metric $\rho_D$ is a special 
metric that is tailored to capture the fact that the mappings in ${\rm et}(\mathbb{C}^2)$ are locally volume preserving mappings. That is encoded
in the Jacobian Condition $\det J_F\equiv 1$. Secondly, since the left (composition) mappings $L_F$ tend to $1/d_F$ similarities on
the metric space $({\rm et}(\mathbb{C}^2),\rho_D)$, when $D\rightarrow\mathbb{C}^2$, our computations of Hausdorff dimensions of
${\rm et}(\mathbb{C}^2)$ and of ${\rm Aut}(\mathbb{C}^2)$ will lead us to the fundamental integral identity (\ref{eq1}), which implies the
Two Dimensional Jacobian Conjecture. We do not know of a "right (composition)" replacement to the metric $\rho_D$ and (unfortunately) the right (composition)
mappings, $R_F$, are not controlled by our metric $\rho_D$ (they are not $\rho_D$ bi-Lipschitz for example).
\end{remark}

We will denote the $s$-dimensional Hausdorff measure of the set $A$, by $H^{s}(A)$. This is the standard notion and terminology
(see for example page 29 of \cite{falconer}). If we fix an \'etale mapping $F\in {\rm et}(\mathbb{C}^2)$ we can define a density
function $f(F,s)$ of the non-negative parameter $s$. We think of $s$ as the dimension with respect to which we calculate our Hausdorff
measures, and require that $H^s(L_F({\rm et}(\mathbb{C}^2))=f(F,s)\cdot H^s({\rm et}(\mathbb{C}^2))$. So $f(F,s)$ is a scaling factor
of the left (composition) mapping $L_F$. In particular if $s_0=\dim_H {\rm et}(\mathbb{C}^2)$ is the Hausdorff dimension of ${\rm et}(\mathbb{C}^2)$
($0<s_0\le\infty$), then one can prove by using our approximate fractal structure on ${\rm et}(\mathbb{C}^2)$, that:
$$
1\le\int_{F\in {\rm et}_P(\mathbb{C}^2)}f(F,s_0)d\mu(F).
$$
We note that we view the integrand $f(F,s_0)$ as a function of the variable $F\in {\rm et}_P(\mathbb{C}^2)$ because $s=s_0$ is fixed. We recall
that when $D\rightarrow\mathbb{C}^2$, then $f(F,s_0)\rightarrow d_F^{-s_0}$ which easily implies the geometric inequality,
\begin{equation}
1\le\int_{F\in {\rm et}_P(\mathbb{C}^2)}\frac{d\mu(F)}{d_F^{s_0}}.
\label{eq2}
\end{equation}
A wealth of arithmetical inequalities follow from this fundamental inequality. At this point we will note that these inequalities become
identities provided that we knew that the so called, disjointness property, was valid. This property is the following:
$$
\forall\,F,G\in {\rm et}_P(\mathbb{C}^2),\,\,\,F\ne G\Rightarrow L_F({\rm et}(\mathbb{C}^2))\cap L_G({\rm et}(\mathbb{C}^2))=\emptyset.
$$
The disjointness property is known in the literature of the theory of fractals as the strong separation property. See \cite{hutchinson}
for that theory. The strong separation property is stronger than the more commonly used, open set condition, which in
most cases suffices to derive the desired results. In our particular case the disjointness property in conjunction with the inequality
(\ref{eq2}) (which becomes an equality) will imply that for every value of the parameter $s$ in the non-degenerate interval
$0\le H^{s_0}({\rm Aut}(\mathbb{C}^2))<s\le\dim_H {\rm et}(\mathbb{C}^2)\le\infty$ we have the identity,
$$
\int_{F\in {\rm et}_P(\mathbb{C}^2)}\frac{d\mu(F)}{d_F^s}\equiv 1.
$$
This is our fundamental identity in equation (\ref{eq1}).

To summarize we mention two central ideas that emerge from the theory of fractals and which come handy in our proof. The first
is the relatively new theory of invariant sets with respect to infinite systems of contractions. Here is a very partial list
of related articles: \cite{anderson,fernau,gm,hutchinson,mw,wicks}. In our case the set of generators might, apriori, be uncountable.
The second idea is that our target metric space is a limiting value of the metric spaces $({\rm et}(\mathbb{C}^2),\rho_D)$ where the
characteristic domain $D$ tends to $\mathbb{C}^2$. In that limiting process the lower bi-Lipschitz constants of the generators $F\in
{\rm et}_p(\mathbb{C}^2)$ tend to the reciprocals of the geometric degrees $d_F^{-1}$ which are reciprocals of natural numbers
that are greater than or equal to $2$.

In order to prove some of the results above we will make use of algebro-geometric tools. These tools are mostly emerging from the theory of the
asymptotic values of \'etale mappings, i.e. the mappings that constitute one of our main objects, the semigroup
$({\rm et}(\mathbb{C}^2),\circ)$. These were developed by several mathematicians in order, among other things, to tackle the two dimensional Jacobian 
Conjecture by more conservative (or direct) methods.

In the next section we will outline the main preparatory such results that will be the most useful for our purposes. We give references to the 
literature in which the reader can find those results.

\section{The asymptotic values of (polynomial) \'etale mappings, \cite{picard}}

We recall the following standard,
\begin{definition}
Let $M$ and $N$ be two manifolds. A differentiable mapping $f\,:\,M\rightarrow N$ is called a diffeomorphism if it is a bijection and if
the inverse mapping $f^{-1}\,:\,N\rightarrow M$ is also differentiable. If both mappings are $r$ times continuously differentiable, then
$f$ is called a $C^r$-diffeomorphism (here $1\le r\le\infty$). Similarly for real or complex analytic diffeomorphisms.
\end{definition}

\begin{definition}
Let $X$ and $Y$ be differentiable manifolds. A mapping $f\,:\,X\rightarrow Y$ is a local diffeomorphism, if every $x\in X$ has an open 
neighborhood $U\subseteq X$, such that $f(U)$ is open in $Y$, and the restriction mapping $f|_U\,:\,U\rightarrow f(U)$ is a diffeomorphism.
\end{definition}

A beautiful theorem of Jacques Hadamard, \cite{h}, gives a necessary and sufficient condition on a local $C^1$-diffeomorphism, 
$f\,:\,\mathbb{R}^n\rightarrow\mathbb{R}^n$ to be, in fact, a (global) diffeomorphism. The condition is written in terms of the
divergence to $\infty$ of a certain improper integral of the first kind, whose (positive) integrand involves $|J_f^{-1}|$, where $J_f$ is the
$n\times n$ Jacobian matrix of $f$, and $|J_f^{-1}|$ is the norm of the operator which is determined by the inverse matrix. It is possible to 
rephrase Hadamard's condition in terms of asymptotic values of the mapping $f\,:\,\mathbb{R}^n\rightarrow\mathbb{R}^n$. We first recall the following,
\begin{definition}
An asymptotic value of the local diffeomorphism $f\,:\,\mathbb{R}^n\rightarrow\mathbb{R}^n$ is a finite limit 
$(a_1,\ldots,a_n)=\lim_{t\rightarrow\infty} f(\sigma(t))$, where $\sigma\,:\,(0,\infty)\rightarrow\mathbb{R}^n$ is a piecewise smooth curve that
goes to infinity, i.e. $\lim_{t\rightarrow\infty} ||\sigma(t)||_2=\infty$. We denote by $||\cdot||_2$ the standard $L^2$ norm on $\mathbb{R}^n$.
The curve $\sigma$ is called an asymptotic tract of $f$ that corresponds to its asymptotic value $(a_1,\ldots,a_n)$. The set of all the asymptotic
values of $f$, is called the asymptotic variety of $f$ and we denote it by $A(f)$.
\end{definition}

The conclusion referred to above that follows from Hadamard's theorem is: A local diffeomorphism $f\,:\,\mathbb{R}^n\rightarrow\mathbb{R}^n$ is, in fact, a
(global) diffeomorphism, if and only if $A(f)=\emptyset$.

\begin{remark}
We deliberately avoid using the more modern and general (and topological) notion of a proper mapping in relation to the above theorem of Hadamard. The
reason is that the notion of asymptotic values and their asymptotic tracts, form computational objects. We will make use of the
formulas that are the result of this computation for the particular case of polynomial \'etale mappings $F\,:\,\mathbb{C}^2\rightarrow\mathbb{C}^2$.
\end{remark}

Let $F\in {\rm et}(\mathbb{C}^2)$. We will now give a complete description of the asymptotic variety $A(F)$ and point at a canonical set
of asymptotic tracts of $F$ that generate the full asymptotic variety. This canonical set of curves will be called the canonical geometric 
basis of $F$ and will be denoted by $R_0(F)$. This basis consists of finitely many rational mappings of the following form:
$$
R(X,Y)=(X^{-\alpha},X^{\beta}Y+X^{-\alpha}\Phi(X)),
$$
where $\alpha\in\mathbb{Z}^+$, $\beta\in\mathbb{Z}^+\cup\{ 0\}$, $\Phi(X)\in\mathbb{C}[X]$ and $\deg\Phi<\alpha+\beta$.
Also the effective $X$ powers (the powers with non-zero coefficients) in $X^{\alpha+\beta}Y+\Phi(X)$ have a gcd which equals $1$. Finally, 
$2\le\gamma\le\beta-\alpha$ where the role of $\gamma\in\mathbb{Z}^+$ will be explained below. The cardinality
of the geometric basis, $|R_0(F)|$, equals the number of all the irreducible components of the affine algebraic curve $A(F)$.
$\forall\,R\in R_0(F)$ we have the so-called, double asymptotic identity $F\circ R=G_R\in\mathbb{C}[X,Y]^2$ where the
polynomial mapping $G_R$ is called the $R$-dual of $F$. Each $R\in R_0(F)$ generates exactly one component of $A(F)$.
This component is normally parametrized by $\{G_R(0,Y)\,|\,Y\in\mathbb{C}\}$. This means that this parameterization is surjective. We will 
denote by $H_R(X,Y)=0$ an implicit representation of this component in terms of the irreducible polynomial $H_R\in\mathbb{C}[X,Y]$. Then there
exists a natural number $\gamma(R)\ge 2$ and a polynomial $S_R(X,Y)\in\mathbb{C}[X,Y]$ such that $S_R(X,Y)=e_R+X\cdot T_R(X,Y)$
for some non-zero polynomial $e_R\in\mathbb{C}[Y]$ and $T_R(X,Y)\in\mathbb{C}[X,Y]-\mathbb{C}[X]$.
The affine curve $S_R(X,Y)=0$ is called the $R$-phantom curve of $F$. The $R$-component of $A(F)$, $H_R(X,Y)=0$, 
is a polynomial curve which is not isomorphic to $\mathbb{A}^1$, and hence in particular must be a singular
irreducible curve. We have the relation:
$$
H_R(G_R(X,Y))=X^{\gamma(R)}S_R(X,Y)=X^{\gamma(R)}(e_R+X\cdot T_R(X,Y)).
$$
The exponent $\gamma(R)$ is the number $\gamma $ that appears above in the double inequality 
$2\le\gamma\le\beta-\alpha$. In our case of the canonical rational mappings $R\in R_0(F)$, we have
${\rm sing}(R)=\{ X=0\}$ (the singularity set of the mapping $R$). The following is true:
$$
G^{-1}_R(H_R(X,Y)=0)=G^{-1}_R(G_R({\rm sing}(R)))={\rm sing}(R)\cup\{S_R(X,Y)=0\}.
$$
Thus the $R$-dual preimage (i.e., the $G_R$-preimage) of the $R$-component of $A(F)$ (which is the $G_R$-image of ${\rm sing}(R)$) is
the union of two curves: The first is ${\rm sing}(R)$ and the second is the $R$-phantom curve of $F$. It can be shown that even if for 
a single $R(X,Y)$ the $R$-phantom curve is empty then the two dimensional Jacobian Conjecture follows.

These are few facts about the asymptotic variety $A(F)$ and the canonical geometric basis $R_0(F)$ of an \'etale mapping
$F\in {\rm et}(\mathbb{C}^2)$. We now seek for more global connections, ${\rm et}(\mathbb{C}^2)$ wide.
In particular we would like to understand how the binary operation in the semigroup ${\rm et}(\mathbb{C}^2)$ affects
the structures of the geometric objects $A(F)$ and of the algebraic objects $R_0(F)$.

\begin{proposition}
If $F,G\in {\rm et}(\mathbb{C}^2)$ then $R_0(G)\subseteq R_0(F\circ G)$, $F(A(G))\subseteq A(F\circ G)$.
\end{proposition}
\noindent
{\bf Proof.} \\
$R\in R_0(G)\Rightarrow G\circ R\in\mathbb{C}[X,Y]^2\Rightarrow F\circ (G\circ R)\in\mathbb{C}[X,Y]^2
\Rightarrow (F\circ G)\circ R\in\mathbb{C}[X,Y]^2\Rightarrow R\in R_0(F\circ G)$. Next we have \\
$(a,b)\in F(A(G))\Rightarrow \exists\,R\in R_0(G)\exists\,Y\in \mathbb{C}\,\,{\rm such}\,{\rm that}\,\,
(a,b)=F((G\circ R))(0,Y))\Rightarrow\exists\,R\in R_0(F\circ G)\exists\,Y\in\mathbb{C}\,\,{\rm such}\,{\rm that}\,\,
(a,b)=((F\circ G)\circ R)(0,Y)\Rightarrow (a,b)\in A(F\circ G)$. $\qed $ \\
\\
The proposition tells us that compositions of \'etale mappings do not decrease the geometric basis of the 
right factor and consequently do not decrease the left image of its asymptotic variety. We naturally ask,
under what conditions the geometric basis of $F\circ G$ is actually larger than that of $G$? In other
words we would like to know when is it true that $R_0(G)\subset R_0(F\circ G)$? This happens exactly
when $\exists\,R\in R_0(F\circ G)-R_0(G)$. This means that $(F\circ G)\circ R\in\mathbb{C}[X,Y]^2$,
$G\circ R\not\in\mathbb{C}[X,Y]^2$. Let $R(X,Y)=(X^{-\alpha},X^{\beta}Y+X^{-\alpha}\Phi(X))$,
$G(X,Y)=(P(X,Y),Q(X,Y))$. Then 
$$
(G\circ R)(X,Y)=(P(X^{-\alpha},X^{\beta}Y+X^{-\alpha}\Phi(X)),Q(X^{-\alpha},X^{\beta}Y+X^{-\alpha}\Phi(X))\in
$$
$$
\in\mathbb{C}(X,Y)^2-\mathbb{C}[X,Y]^2\,\,{\rm in}\,{\rm fact}\,\,\in\mathbb{C}[X^{-1},X,Y]^2-\mathbb{C}[X,Y]^2.
$$
We clearly have ${\rm sing}(G\circ R)\subseteq {\rm sing}(R)$ and so ${\rm sing}(G\circ R)=\{ X=0\}$. By 
$F\circ (G\circ R)=(F\circ G)\circ R\in\mathbb{C}[X,Y]^2$ we have $G\circ R\in R(F)$. The set $R(F)$ is the set of all the asymptotic tracts of $F$. 
The asymptotic tract $G\circ R$ is not necessarily a member of the canonical geometric basis of $F$. We recall that the canonical geometric basis 
of $F$, $R_0(F)$ contains finitely many rational mappings of the form:
$$
S(X,Y)=(X^{-a},X^bY+X^{-a}\Psi(X)).
$$
Since $G\in {\rm et}(\mathbb{C}^2)$ it follows that $|\mathbb{C}^2-G(\mathbb{C}^2)|<\infty$ (a similar phenomenon
to the Picard Theorem for entire functions of a single complex variable). If $L$ is an asymptotic tract of $F$ then $G^{-1}(L)$ 
can not be a bounded subset of $\mathbb{C}^2$. The reason is that if $\overline{G^{-1}(L)}$ is compact, then $G(\overline{G^{-1}(L)})$ is
compact and since $L\subseteq G(G^{-1}(L))\subseteq G(\overline{G^{-1}(L)})$ this would imply the
contradiction that $L$ is bounded (and hence can not be an asymptotic tract). Hence $G^{-1}(L)$ has at least one
component, say $L_1$, that goes to infinity. So $F\circ G$ has a limit along $L_1$ which equals the
above asymptotic value of $F$. This proves the following generalization of the second part of Proposition 2.5,
namely,

\begin{proposition}
If $F,G\in {\rm et}(\mathbb{C}^2)$ then $A(F)\cup F(A(G))=A(F\circ G)$.
\end{proposition}
\noindent
This proposition implies that if $A(F)\subset F(A(G))$ then necessarily $R_0(G)\subset R_0(F\circ G)$
because, as shown in the proof of Proposition 2.5 $\forall\,R\in R_0(G)$, $((F\circ G)\circ R)({\rm sing}(R))
\subseteq F(A(G))$.

\begin{proposition}
Let $F\in {\rm et}(\mathbb{C}^2)$. If $\exists\,G\in {\rm et}(\mathbb{C}^2)$ such that $R_0(G)=R_0(F\circ G)$,
then $F(\mathbb{C}^2)=\mathbb{C}^2$, i.e. $F$ is a surjective mapping.
\end{proposition}
\noindent
{\bf Proof.} \\
Since $F\in {\rm et}(\mathbb{C}^2)$ we have $\mathbb{C}^2-F(\mathbb{C}^2)\subseteq A(F)$, because in
this case the only points in the complement of the image of $F$ are the finitely many Picard exceptional 
values of $F$ which are asymptotic values of $F$. If, as the assumption says $R_0(G)=R_0(F\circ G)$
then by Proposition 2.6 we must have $A(F)\subseteq F(A(G))\subseteq F(\mathbb{C}^2)$ and so there are
no Picard exceptional values of the mapping $F$. $\qed $ \\

\section{The right and the left (composition) mappings on ${\rm et}(\mathbb{C}^2)$}

\begin{remark}
In this short section we will use the notions of the right (composition) mapping and  the left (composition) mapping.
These were defined in Definition 1.4.
\end{remark}

\begin{proposition}
The mappings $R_F$, $L_F$ are not surjective if and only if $F\not\in {\rm Aut}(\mathbb{C}^2)$. In fact in this
case we have $R_F({\rm et}(\mathbb{C}^2))\subset {\rm et}(\mathbb{C}^2)-{\rm Aut}(\mathbb{C}^2)$,
$L_F({\rm et}(\mathbb{C}^2))\subset {\rm et}(\mathbb{C}^2)-{\rm Aut}(\mathbb{C}^2)$.
\end{proposition}
\noindent
{\bf Proof.} \\
By Proposition 2.5 we have: $R_0(R_F(G))=R_0(G\circ F)\supseteq R_0(F)\ne\emptyset,\,\,\,A(L_F(G))=A(F\circ G)\supseteq A(F)\ne\emptyset$. $\qed $

\begin{proposition}
$R_F$ is injective.
\end{proposition}
\noindent
{\bf Proof.} \\
$R_F(G)=R_F(H)\Rightarrow G\circ F=H\circ F$. Since $F\in   {\rm et}(\mathbb{C}^2)$ we have $|\mathbb{C}^2-F(\mathbb{C}^2)|
<\infty$ and by the assumption $G|_{F(\mathbb{C}^2)}=H|_{F(\mathbb{C}^2)}$. Hence $G\equiv H$. $\qed $ \\
\\
We naturally inquire if also $L_F$ is injective. This, however, will be proved later after considerable amount of
preparations.

\section{If ${\rm et}(\mathbb{C}^2)\ne{\rm Aut}(\mathbb{C}^2)$, then ${\rm et}(\mathbb{C}^2)$ is fractal like}

We start this section by remarking on a few topological properties of the image of the right (composition) mapping $R_F({\rm et}(\mathbb{C}^2))$.
Let $F\in {\rm et}(\mathbb{C}^2)$ and $\forall\,\varepsilon>0$ we define
$$
V_{\varepsilon}(F)=\{ N\in {\rm et}(\mathbb{C}^2)\,|\,\max_{|X|,|Y|\le 1}\parallel F(X,Y)-N(X,Y)\parallel_2<\varepsilon\}.
$$
These sets form a local basis at $F$ for the topology we impose for now on ${\rm et}(\mathbb{C}^2)$. If
$N_{\varepsilon}\in V_{\varepsilon}(F)$ then as $\varepsilon\rightarrow 0^+$ the coefficients of $N_{\varepsilon}$
tend to the corresponding coefficients of $F$. To see that let $F=(F_1,F_2)$, $N=(N_1,N_2)$. Then by the $H^2$ theory
for analytic functions we have:
$$
\max_{|X|,|Y|\le 1}\parallel F(X,Y)-N(X,Y)\parallel_2=\parallel F(e^{i\phi},e^{i\theta})-N(e^{i\phi},e^{i\theta})\parallel_2=
$$
$$
=\left(\int^{2\pi}_{\theta=0}\int^{2\pi}_{\phi=0}\left\{|F_1(e^{i\phi},e^{i\theta})-N_1(e^{i\phi},e^{i\theta})|^2+
|F_2(e^{i\phi},e^{i\theta})-N_2(e^{i\phi},e^{i\theta})|^2\right\}d\phi d\theta\right)^{1/2}=
$$
$$
=\left\{\sum_{k,l}\left(|a^{(1)}_{kl}-b^{(1)}_{kl}|^2+|a^{(2)}_{kl}-b^{(2)}_{kl}|^2\right)\right\}^{1/2}.
$$
Here we used the notation $F_j(X,Y)=\sum_{k,l}a^{(j)}_{kl}X^kY^l$ and $N_j(X,Y)=\sum_{k,l}b^{(j)}_{kl}X^kY^l$,
$j=1,2$. This shows that $\forall\,k,l\in\mathbb{Z}^+$ $\forall\,j=1,2$ we have
$$
\max_{|X|,|Y|\le 1}\parallel F(X,Y)-N(X,Y)\parallel_2\ge |a^{(j)}_{kl}-b^{(j)}_{kl}|.
$$
Our claim on the coefficients of the mappings follows.

Let $F\in {\rm et}(\mathbb{C}^2)$. We suspect that the image of the right (composition) mapping $R_F({\rm et}(\mathbb{C}^2))$ is
a closed subset of ${\rm et}(\mathbb{C}^2)$ in the $L^2$-topology which was introduced above. We
recall that the two dimensional Jacobian Conjecture is equivalent to ${\rm et}(\mathbb{C}^2)={\rm Aut}(\mathbb{C}^2)$.
Thus we assume from now on that ${\rm et}(\mathbb{C}^2)\ne{\rm Aut}(\mathbb{C}^2)$ in order to see the implications
of this assumption. The right (composition) mapping for any $F\in {\rm et}(\mathbb{C}^2)$, $R_F\,:\,{\rm et}(\mathbb{C}^2)
\rightarrow {\rm et}(\mathbb{C}^2)$, $R_F(G)=G\circ F$, is a continuous injection (Proposition 3.3). Continuity
here means, say, with respect to the $L^2$-topology. Also $F\in {\rm Aut}(\mathbb{C}^2)\Leftrightarrow R_F({\rm et}(\mathbb{C}^2))=
{\rm et}(\mathbb{C}^2)$. If $F\in {\rm et}(\mathbb{C}^2)-{\rm Aut}(\mathbb{C}^2)$, then $R_F({\rm et}(\mathbb{C}^2))\subset
{\rm et}(\mathbb{C}^2)-{\rm Aut}(\mathbb{C}^2)$. The following is well known (\cite{rw}, Theorem 2.3) \\
\\
{\bf Theorem (Kamil Rusek and Tadeusz Winiarski).} ${\rm Aut}(\mathbb{C}^n)$ is a closed subset of $({\rm et}(\mathbb{C}^n),L^2)$. \\
\\
This follows from a formal analog of Cartan's theorem on sequences of bi-holomorphisms of a bounded domain.
Actually the topology referred to in \cite{rw}, on ${\rm et}(\mathbb{C}^n)$ is that of uniform convergence
on compact subsets of $\mathbb{C}^n$. This topology is identical in this case with the compact-open
topology.

We will need some preparations in order to arrive at the fractal like structure we will put on ${\rm et}(\mathbb{C}^2)$.
We start by indicating an easy upper bound (well-known) of the generic size of a fiber of a mapping $F=(P,Q)\in {\rm et}(\mathbb{C}^2)$. If we denote
$\deg P(X,Y)=n$ and $\deg Q(X,Y)=m$ then $\forall\,(a,b)\in\mathbb{C}^2$ the $F$ fiber over $(a,b)$ is
$F^{-1}(a,b)=\{(x,y)\in\mathbb{C}^2\,|\,F(x,y)=F(a,b)\}$. It is well known that this set is a finite subset
of $\mathbb{C}^2$ and, by the Bezout Theorem we have 
$$
|\{(x,y)\in\mathbb{C}^2\,|\,F(x,y)=(a,b)\}|=|F^{-1}(a,b)|\le n\cdot m.
$$
Moreover, as indicated in the "Introduction" section, there is a number that we will denote by $d_F$ such that generically in $(a,b)$ we have 
$|F^{-1}(a,b)|=d_F$. This means that $\{(a,b)\in\mathbb{C}^2\,|\,|F^{-1}(a,b)|\ne d_F\}$ is a closed
and proper Zariski subset of $\mathbb{C}^2$. In fact $\forall\,(a,b)\in\mathbb{C}^2,\,\,|F^{-1}(a,b)|\ne d_F
\Rightarrow |F^{-1}(a,b)|< d_F$. Thus we have $d_F=\max\{|F^{-1}(a,b)|\,|\,(a,b)\in\mathbb{C}^2\}$.

\begin{definition}
Let $F\in {\rm et}(\mathbb{C}^2)$. We will denote $d_F=\max\{|F^{-1}(a,b)|\,|\,(a,b)\in\mathbb{C}^2\}$. We will
call $d_F$ the geometrical degree of the \'etale mapping $F$.
\end{definition}

\begin{proposition}
$\forall\,F,G\in {\rm et}(\mathbb{C}^2),\,\,d_{F\circ G}=d_F\cdot d_G$.
\end{proposition}
\noindent
This is  a well known result. We include one of its proofs for convenience. \\
{\bf Proof.} \\
$\forall\,(a,b)\in\mathbb{C}^2,\,\,(F\circ G)^{-1}(a,b)=G^{-1}(F^{-1}(a,b))$. But generically in $(a,b)$
$|F^{-1}(a,b)|=d_F$ and generically in $(c,d)$, $|G^{-1}(c,d)|=d_G$. $\qed $ \\

\begin{definition}
An \'etale mapping $F\in {\rm et}(\mathbb{C}^2)$ is a composite mapping if $\exists\,G,H\in {\rm et}(\mathbb{C}^2)-
{\rm Aut}(\mathbb{C}^2)$ such that $F=G\circ H$. An \'etale mapping $A\in {\rm et}(\mathbb{C}^2)-{\rm Aut}(\mathbb{C}^2)$ is
a prime mapping if it is not a composite mapping. This is equivalent to: $A=B\circ C$ for some $B,C\in  {\rm et}(\mathbb{C}^2)$
$\Rightarrow B\in {\rm Aut}(\mathbb{C}^2)\vee C\in {\rm Aut}(\mathbb{C}^2)$. The subset of ${\rm et}(\mathbb{C}^2)$
of all the prime mappings will be denoted by ${\rm et_p}(\mathbb{C}^2)$. Thus the set of all the composite
\'etale mappings is ${\rm et}(\mathbb{C}^2)-{\rm et_p}(\mathbb{C}^2)$.
\end{definition}

An easy consequence of the definitions is the following,

\begin{proposition}
$\forall\,F\in {\rm et}(\mathbb{C}^2)-{\rm et_p}(\mathbb{C}^2)$, $d_F$ is not a prime number. Equivalently,
$\forall\,F\in {\rm et}(\mathbb{C}^2)$, $d_F$ is a prime number $\Rightarrow F\in {\rm et_p}(\mathbb{C}^2)$.
\end{proposition}
\noindent
{\bf Proof.} \\
$F\in {\rm et}(\mathbb{C}^2)-{\rm et_p}(\mathbb{C}^2)\Rightarrow\exists\,G,H\in {\rm et}(\mathbb{C}^2)-
{\rm Aut}(\mathbb{C}^2)$ such that $F=G\circ H$ (by Definition 4.3) $\Rightarrow d_F=d_G\cdot d_H,\,\,d_G,d_H>1$
(by Proposition 4.2 and the fact $d_M=1\Leftrightarrow M\in {\rm Aut}(\mathbb{C}^2)$) $\Rightarrow d_F$ is
a composite integer. $\qed $

\begin{theorem}. \\
{\rm 1)} ${\rm et_p}(\mathbb{C}^2)\ne\emptyset$ \\
{\rm 2)} $\forall\,F\in {\rm et}(\mathbb{C}^2)\,\,\exists\,k\in\mathbb{Z}^+\cup\{ 0\}\,\,\exists\,
A_0\in {\rm Aut}(\mathbb{C}^2)\,\exists\,P_1,\ldots,P_k\in {\rm et_p}(\mathbb{C}^2)$ such that
$F=A_0\circ P_1\circ\ldots\circ P_k$.
\end{theorem}
\noindent
{\bf Proof.} \\
If ${\rm et_p}(\mathbb{C}^2)=\emptyset $ then ${\rm et}(\mathbb{C}^2)-{\rm Aut}(\mathbb{C}^2)$ are all
composite \'etale mappings. Let $F\in {\rm et}(\mathbb{C}^2)-{\rm Aut}(\mathbb{C}^2)$, then
$\exists\,G_1,G'_2\in {\rm et}(\mathbb{C}^2)-{\rm Aut}(\mathbb{C}^2)$ such that $F=G_1\circ G'_2$. So
$\exists\,G_2,G'_3\in {\rm et}(\mathbb{C}^2)-{\rm Aut}(\mathbb{C}^2)$ such that $G'_2=G_2\circ G'_3$.
Hence $F=G_1\circ G_2\circ G'_3$. Continuing this we get for any $k\in\mathbb{Z}^+$ $\exists\,G_1,\ldots,G_k\in
{\rm et}(\mathbb{C}^2)-{\rm Aut}(\mathbb{C}^2)$ such that $F=G_1\circ\ldots\circ G_k$ and by Proposition 4.2
$d_F=\prod_{j=1}^k d_{G_j}$. But $\forall\,1\le j\le k$, $d_{G_j}\ge 2$ and so $\forall\,k\in\mathbb{Z}^+$,
$d_F\ge 2^k$ a contradiction to $d_F<\infty $. Thus ${\rm et_p}(\mathbb{C}^2)\ne\emptyset$. \\
Now part 2 is standard, for if $F\in {\rm Aut}(\mathbb{C}^2)$ we take $A_0=F$ and $k=0$. If $F\in
{\rm et_p}(\mathbb{C}^2)$ we take $A_0={\rm id}$, $k=1$, and $P_1=F$. If $F\in {\rm et}(\mathbb{C}^2)-
{\rm et_p}(\mathbb{C}^2)$ then $F=G\circ H$ for some $G,H\in {\rm et}(\mathbb{C}^2)-{\rm Aut}(\mathbb{C}^2)$.
So by Proposition 4.2 $d_F=d_G\cdot d_H$ and since $d_G,d_H\ge 2$ it follows that $d_G,d_H<d_F$ and
we conclude the proof of part 2 using induction on the geometrical degree. Namely $G=P_1\circ\ldots\circ P_m$,
$H=P_{m+1}\circ\ldots\circ P_k$ for $m\ge 1$, $k\ge m+1$ and some primes $P_1,\ldots,P_k\in {\rm et_p}(\mathbb{C}^2)$ $\qed $ \\

\begin{definition}
We define a relation $\sim $ on ${\rm et}(\mathbb{C}^2)$ by: $\forall\,F,G\in  {\rm et}(\mathbb{C}^2)$
$F\sim G\Leftrightarrow \exists\,A,B\in {\rm Aut}(\mathbb{C}^2),\,\,F=A\circ G\circ B$.
\end{definition}

\begin{remark}
The relation $\sim $ is an equivalence relation on ${\rm et}(\mathbb{C}^2)$. For $F\sim F$ because
$F={\rm id}\circ F\circ {\rm id}$. Also $F\sim G\Rightarrow F=A\circ G\circ B\Rightarrow G=A^{-1}\circ F\circ B^{-1}\Rightarrow
G\sim F$.
Finally $F\sim G,\,G\sim H\Rightarrow F=A\circ G\circ B,\,G=C\circ H\circ D\Rightarrow F=(A\circ C)\circ H\circ (D\circ B)
\Rightarrow F\sim H$. We could have defined two similar equivalence relations on ${\rm et}(\mathbb{C}^2)$
by restricting to compositions with automorphisms from one side only (left side only, or, right side only). We will denote these relations
by $\sim_R$ and $\sim_L$. For example $F\sim_L G\Leftrightarrow \exists\,A\in {\rm Aut}(\mathbb{C}^2),\,\,F=A\circ G$.
\end{remark}

\begin{definition}
The right partial order on ${\rm et}(\mathbb{C}^2)/\sim_R$ is defined by: $[F]\preceq_R [G]\Leftrightarrow
R_F({\rm et}(\mathbb{C}^2))\subseteq R_G({\rm et}(\mathbb{C}^2))$.
\end{definition}

\begin{proposition}
The relation $\preceq_R$ is a partial order on ${\rm et}(\mathbb{C}^2)/\sim_R$.
\end{proposition}
\noindent
{\bf Proof.} \\
The claim is clear because $\subseteq$ is a partial order on any family of sets. However,
here it is instructive to notice the anti-symmetric property also from the point of view of our
particular setting. Namely $[F]\preceq_R [G]\wedge [G]\preceq_R [F]\Leftrightarrow
R_F({\rm et}(\mathbb{C}^2))\subseteq R_G({\rm et}(\mathbb{C}^2))\wedge R_G({\rm et}(\mathbb{C}^2))\subseteq R_F({\rm et}(\mathbb{C}^2))
\Rightarrow F\in R_G({\rm et}(\mathbb{C}^2))\wedge G\in R_F({\rm et}(\mathbb{C}^2)\Leftrightarrow
\exists\,M,N\in {\rm et}(\mathbb{C}^2)\,\,{\rm such}\,\,{\rm that}\,\,F=M\circ G\wedge G=N\circ F
\Rightarrow F=(M\circ N)\circ F$. Since $F(\mathbb{C}^2)$ is co-finite in $\mathbb{C}^2$, the last equation implies
that $M\circ N={\rm id}$, so $M,N\in {\rm Aut}(\mathbb{C}^2)$, $M=N^{-1}$ and so $[F]=[G]$. $\qed $

\begin{theorem}
Every $\preceq_R$-increasing chain is finite, i.e. it stabilizes.
\end{theorem}
\noindent
{\bf Proof.} \\
We will argue by a contradiction. Suppose that there is an infinite $\preceq_R$-increasing chain. Then there is an infinite sequence
$F_1, F_2, F_3,\ldots\in {\rm et}(\mathbb{C}^2)-{\rm et_p}(\mathbb{C}^2)$ such that $[F_1]\preceq_R [F_2]\preceq_R [F_3]\preceq_R\ldots $, and hence by 
Definition 4.8: $R_{F_1}({\rm et}(\mathbb{C}^2))\subset R_{F_2}({\rm et}(\mathbb{C}^2))\subset R_{F_3}({\rm et}(\mathbb{C}^2))
\subset\ldots $. Hence $\exists\,M_j\in {\rm et}(\mathbb{C}^2)-{\rm Aut}(\mathbb{C}^2)$ such that $F_j=M_j\circ F_{j+1}$ for $j=1,2,3,\ldots $. This implies that
$\forall\,k\in\mathbb{Z}^+$, $F_1=M_1\circ M_2\circ\ldots\circ M_k\circ F_{k+1}$ and so as in the argument in the proof of Theorem 4.5 we 
obtain $d_{F_1}=d_{M_1}d_{M_2}\ldots d_{M_k}d_{F_{k+1}}\ge 2^k$. This contradicts the fact that $d_{F_1}<\infty$ and concludes the proof of 
Theorem 4.10. $\qed $

\begin{remark}
Let us consider an \'etale mapping $F\in {\rm et}(\mathbb{C}^2)$, say $F=(P,Q)\in\mathbb{C}[X,Y]^2$ where $\deg P=n$, $\deg Q=m$. By the
Bezout Theorem we have $d_F\le n\cdot m$. If either $n=1$ or $m=1$ then it is well known that $F\in {\rm Aut}(\mathbb{C}^2)$ and so $d_F=1$. This
follows because if $n=1$ or $m=1$, the mapping $F$ is injective on a straight line and it is well known that such an \'etale mapping
must belong to ${\rm Aut}(\mathbb{C}^2)$. If $(a,b)\in\mathbb{C}^2$ satisfies $|F^{-1}(a,b)|< d_F$ then it is well known that $F(a,b)$ is an asymptotic value
of $F$ and there are exactly $d_F-|F^{-1}(a,b)|$ points on the line at infinity that $F$ maps to $F(a,b)$. Thus we expect some relations
between the structure of the asymptotic variety $A(F)$ and $d_F$ and the size of the fiber at the given point $(a,b)$. Here is a sketch for such relations:
$$
|F^{-1}(a,b)|+|{\rm the\,\,different\,\,points\,\,on\,\,the\,\,line\,\,at\,\,infinity\,\,that\,\,F\,\,maps\,\,to}
$$
$$
F(a,b)|=d_F,
$$
so
$$
|F^{-1}(a,b)|+\sum_{R\in R_0(F)}\sum_{\{Y\,|\,G_R(0,Y)=F(a,b)\}} 1=d_F,
$$
hence
$$
|F^{-1}(a,b)|+\sum_{R\in R_0(F)}|\{Y\,|\,G_R(0,Y)=F(a,b)\}| =d_F
$$
Here is an example of a crude estimate we can get:
$$
d_F\le |F^{-1}(a,b)|+\sum_{R\in R_0(F)}\deg (F\circ R)(0,Y).
$$
Let us denote $D_F=\max\{\deg (G_R)(0,Y)\,|\,R\in R_0(F)\}$ and recall that $|R_0(F)|=|{\rm the\,\,components\,\,of}\,\,A(F)|$.
Then we get:
$$
d_F\le |F^{-1}(a,b)|+D_F\cdot |{\rm the\,\,components\,\,of}\,\,A(F)|=|F^{-1}(a,b)|+D_F\cdot |R_0(F)|.
$$
Thus if $D_F\cdot |{\rm the\,\,components\,\,of}\,\,A(F)|=D_F\cdot |R_0(F)|<d_F$ we conclude that $F$ is a surjective mapping.
\end{remark}

\begin{remark}
We do not expect a claim similar to that made in Theorem 4.10 to be valid for decreasing $\preceq_R$-chains. Namely we expect that
there are infinite decreasing $\preceq_R$-chains (provided, of course, that ${\rm et}(\mathbb{C}^2)-{\rm Aut}(\mathbb{C}^2)\ne\emptyset$.)
Thus if $F\in {\rm et}(\mathbb{C}^2)-{\rm Aut}(\mathbb{C}^2)$ and if we take a sequence $H_n\in {\rm et}(\mathbb{C}^2)-{\rm Aut}(\mathbb{C}^2)$
(for example $H_n=F^{\circ n}$) and define $G_n=H_n\circ\ldots\circ H_1\circ F$, then $\ldots \preceq_R G_n \preceq_R\ldots \preceq_R G_1 \preceq_R F$ 
and $\ldots\subset R_{G_n}({\rm et}(\mathbb{C}^2))\subset\ldots\subset R_{G_1}({\rm et}(\mathbb{C}^2))\subset R_F ({\rm et}(\mathbb{C}^2))$.
\end{remark}
\noindent

\begin{proposition}. \\
{\rm 1)} If $F\in {\rm et}(\mathbb{C}^2)$ and $G\in R_F({\rm et}(\mathbb{C}^2))$, then $R_G({\rm et}(\mathbb{C}^2))\subseteq R_F({\rm et}(\mathbb{C}^2))$. \\
{\rm 2)} If $F\in {\rm et}(\mathbb{C}^2)$, $G\in R_F({\rm et}(\mathbb{C}^2))$, and $G$ and $F$ are not associates (which means here $\forall\,H\in {\rm Aut}(\mathbb{C}^2)$, $G\ne H\circ F$ in other words $F\not\sim_L G$), then $R_G({\rm et}(\mathbb{C}^2))\subset R_F({\rm et}(\mathbb{C}^2))$. \\
{\rm 3)} $\forall\,F\in {\rm et}(\mathbb{C}^2)$ the spaces $(R_F({\rm et}(\mathbb{C}^2)),L^2)$ and $({\rm et}(\mathbb{C}^2),L^2)$ are homeomorphic. \\
{\rm 4)} There exists an infinite index set $I$ and a family of \'etale mappings $\{ F_i\,|\,i\in I\}\subseteq {\rm et}(\mathbb{C}^2)$ such that
$$
{\rm et}(\mathbb{C}^2)={\rm Aut}(\mathbb{C}^2)\cup\bigcup_{i\in I}R_{F_i}({\rm et}(\mathbb{C}^2)),
$$
so that $\forall\,i\in I$, ${\rm Aut}(\mathbb{C}^2)\cap R_{F_i}({\rm et}(\mathbb{C}^2))=\emptyset$, and $\forall\,i,j\in I$, if $i\ne j$ then
$F_i\not\in R_{F_j}({\rm et}(\mathbb{C}^2))$ and $F_j\not\in R_{F_i}({\rm et}(\mathbb{C}^2))$, and $R_{F_i}({\rm et}(\mathbb{C}^2))$ is
homeomorphic to $R_{F_j}({\rm et}(\mathbb{C}^2))$, and both are homeomorphic to ${\rm et}(\mathbb{C}^2)$.
\end{proposition}
\noindent
{\bf Proof.} \\
1) $H\in R_G({\rm et}(\mathbb{C}^2))\Rightarrow\exists\,M\in {\rm et}(\mathbb{C}^2)\,\,{\rm such}\,\,{\rm that}\,\,H=M\circ G$.
$G\in R_F({\rm et}(\mathbb{C}^2))\Rightarrow\exists\,N\in {\rm et}(\mathbb{C}^2)\,\,{\rm such}\,\,{\rm that}\,\,G=N\circ F$. Hence we conclude that
$H=M\circ (N\circ F)=(M\circ N)\circ F\in R_F({\rm et}(\mathbb{C}^2))$. \\
2) $G\in R_F({\rm et}(\mathbb{C}^2))$ and is not an associate of $F$ $\Rightarrow\exists\, N\in {\rm et}(\mathbb{C}^2)-{\rm Aut}(\mathbb{C}^2)$ such that
$G=N\circ F$. So $F\not\in R_G({\rm et}(\mathbb{C}^2))$ otherwise $F=M\circ (N\circ F)=(M\circ N)\circ F$ but $M\circ N\not\in {\rm Aut}(\mathbb{C}^2)$ (it
is not injective). The equation $F=(M\circ N)\circ F$ is equivalent to $M\circ N={\rm id}$ because $\mathbb{C}^2-F(\mathbb{C}^2)$ is a finite set. \\
3) The mapping $f=R_F\,:\,{\rm et}(\mathbb{C}^2)\rightarrow R_F({\rm et}(\mathbb{C}^2))$, $f(G)=G\circ F=R_F(G)$ is an homeomorphism (it is a
bijection and both $f$ and $f^{-1}$ are sequentially continuous). \\
4) We use the relation $\sim_R$ on ${\rm et}(\mathbb{C}^2)$ which was defined by $F\sim_R G\Leftrightarrow\exists\,\Phi\in{\rm Aut}(\mathbb{C}^2)$ such that $F=\Phi\circ G$.
Then $\sim_R $ is an equivalence relation ($F\sim_R F$ by $F={\rm id}\circ F$, $F\sim_R G\Leftrightarrow F=\Phi\circ G\Leftrightarrow G=\Phi^{-1}\circ F
\Leftrightarrow G\sim_R F$, $F\sim_R G\,\wedge\,G\sim_R H\Leftrightarrow F=\Phi_1\circ G\,\wedge\,G=\Phi_2\circ H\Rightarrow F=(\Phi_1\circ\Phi_2)\circ H
\Rightarrow F\sim_R H$). We order the set of $\sim_R $ equivalence classes ${\rm et}(\mathbb{C}^2)/\sim_R$ by $[F]\preceq [G]\Leftrightarrow F\in
R_G({\rm et}(\mathbb{C}^2)) \Leftrightarrow R_F({\rm et}(\mathbb{C}^2))\subseteq R_G({\rm et}(\mathbb{C}^2))$. This relation is clearly reflexive and transitive by Proposition 
4.13(1), and it is also anti-symmetric for $[F]\preceq [G]\wedge [G]\preceq [F]\Leftrightarrow R_F({\rm et}(\mathbb{C}^2))\subseteq R_G({\rm et}(\mathbb{C}^2))
\wedge R_G({\rm et}(\mathbb{C}^2))\subseteq R_F({\rm et}(\mathbb{C}^2))\Leftrightarrow R_F({\rm et}(\mathbb{C}^2))= R_G({\rm et}(\mathbb{C}^2))
\Rightarrow F=N\circ G\wedge G=M\circ F\Rightarrow F=(N\circ M)\circ F\Rightarrow N\circ M={\rm id}\Rightarrow N=M^{-1}\in{\rm Aut}(\mathbb{C}^2)
\Rightarrow [F]=[G]$. Any increasing chain in ${\rm et}(\mathbb{C}^2)-{\rm Aut}(\mathbb{C}^2)/\sim_R$ is finite (by Theorem 4.10). Hence every maximal increasing chain 
contains a maximal element $[F]$, and $R_F({\rm et}(\mathbb{C}^2))$ contains the union of the images of the right (composition) mappings of all the elements in this maximal chain
( by Proposition 4.13(1)). We define $I=\{ [F]\in ({\rm et}(\mathbb{C}^2)-{\rm Aut}(\mathbb{C}^2))/\sim_R\,|\,[F]\,\,{\rm is\,\,the\,\,maximum\,\,of\,\,a\,\,maximal\,\,
length\,\,chain}\}$. Then 
$$
{\rm et}(\mathbb{C}^2)={\rm Aut}(\mathbb{C}^2)\cup\bigcup_{i\in I}R_{F_i}({\rm et}(\mathbb{C}^2)).
$$
The union on the right equals ${\rm et}(\mathbb{C}^2)$ because any $F\in {\rm et}(\mathbb{C}^2)$ is either in ${\rm Aut}(\mathbb{C}^2)$ or $[F]$ belongs to some
maximal length chain in $({\rm et}(\mathbb{C}^2)-{\rm Aut}(\mathbb{C}^2))/\sim_R$ and so $R_F({\rm et}(\mathbb{C}^2))$ is a subset of $R_{F_i}({\rm et}(\mathbb{C}^2))$
where $[F_i]$ is the maximum of that chain. Clearly if $i\ne j$ then $F_i\not\in R_{F_j}({\rm et}(\mathbb{C}^2))\wedge F_j\not\in R_{F_i}({\rm et}(\mathbb{C}^2))$
for $[F_i]$ and $[F_j]$ the maxima of two different chains. So $R_{F_i}({\rm et}(\mathbb{C}^2))\ne R_{F_j}({\rm et}(\mathbb{C}^2))$. Finally, the index set $I$ is 
an infinite set. This follows because any finite union of the form:
$$
{\rm Aut}(\mathbb{C}^2)\cup R_{F_1}({\rm et}(\mathbb{C}^2))\cup\ldots\cup R_{F_k}({\rm et}(\mathbb{C}^2))
$$
is such that any mapping $H$ in it is either a $\mathbb{C}^2$-automorphism or $R_0(H)\supseteq R_0(F_1)\ne\emptyset\vee\ldots\vee R_0(H)\supseteq R_0(F_k)\ne\emptyset$.
Hence the argument in the proof of Proposition 3.2 implies that $R_{F_1}({\rm et}(\mathbb{C}^2))\cup\ldots\cup R_{F_k}({\rm et}(\mathbb{C}^2))\subset
{\rm et}(\mathbb{C}^2)-{\rm Aut}(\mathbb{C}^2)$. $\qed $ \\
\\
\noindent
It could have been convenient if the following claim were valid: If $F, G\in{\rm et}(\mathbb{C}^2)$ satisfy $F\not\in R_G({\rm et}(\mathbb{C}^2))$ and
$G\not\in R_F({\rm et}(\mathbb{C}^2))$ then $R_F({\rm et}(\mathbb{C}^2))\cap R_G({\rm et}(\mathbb{C}^2))=\emptyset$. If this were true we could have
sharpened part (4) of Proposition 4.13. However, we can not prove that and as a result for any $F, G\in{\rm et}(\mathbb{C}^2)$ all the possibilities can
occur, i.e.
$$
R_F({\rm et}(\mathbb{C}^2))\subseteq R_G({\rm et}(\mathbb{C}^2))\,\,{\rm or}\,\,R_G({\rm et}(\mathbb{C}^2))\subseteq R_F({\rm et}(\mathbb{C}^2)),
$$
or
$$
R_F({\rm et}(\mathbb{C}^2))\cap R_G({\rm et}(\mathbb{C}^2))=\emptyset,
$$
or
$$
R_F({\rm et}(\mathbb{C}^2))\cap R_G({\rm et}(\mathbb{C}^2))\not\in\{\emptyset,R_F({\rm et}(\mathbb{C}^2)),R_G({\rm et}(\mathbb{C}^2))\}.
$$
\begin{proposition}
If $F, G\in{\rm et}(\mathbb{C}^2)$ and $R_F({\rm et}(\mathbb{C}^2))\cap R_G({\rm et}(\mathbb{C}^2))\ne\emptyset$, then $\exists\,H\in {\rm et}(\mathbb{C}^2)$
such that $R_H({\rm et}(\mathbb{C}^2))\subseteq R_F({\rm et}(\mathbb{C}^2))\cap R_G({\rm et}(\mathbb{C}^2))$.
\end{proposition}
\noindent
{\bf Proof.} \\
Let $H\in R_F({\rm et}(\mathbb{C}^2))\cap R_G({\rm et}(\mathbb{C}^2))$. Then by part (1) of Proposition 4.13  we have
$R_H({\rm et}(\mathbb{C}^2))\subseteq R_F({\rm et}(\mathbb{C}^2))$ and also $R_H({\rm et}(\mathbb{C}^2))\subseteq R_G({\rm et}(\mathbb{C}^2))$. $\qed $ \\

\begin{proposition}
The family $\{ R_F({\rm et}(\mathbb{C}^2))\,|\,F\in {\rm et}(\mathbb{C}^2)\}$ is a basis of a topology on ${\rm et}(\mathbb{C}^2)$
\end{proposition}
\noindent
{\bf Proof.} \\
Since ${\rm et}(\mathbb{C}^2)=\bigcup_{F\in {\rm et}(\mathbb{C}^2)} R_F({\rm et}(\mathbb{C}^2))$, the claim follows by Proposition 4.14. $\qed $ \\
\\
\noindent
Thus we obtain the following topology, $\tau_R$ on ${\rm et}(\mathbb{C}^2)$: $\tau_R=\{\bigcup_{j\in J} R_{F_j}({\rm et}(\mathbb{C}^2))\,|\,F_j\in {\rm et}(\mathbb{C}^2),\,j\in J\}$.
We will call $\tau_R$, the right (composition) topology on ${\rm et}(\mathbb{C}^2)$.

\begin{proposition}
The space $({\rm et}(\mathbb{C}^2),\tau_R)$ is not Hausdorff.
\end{proposition}
\noindent
{\bf Proof.} \\
We will show that $\tau_R$ can not separate two different points in ${\rm Aut}(\mathbb{C}^2)$. For if $F,G\in {\rm Aut}(\mathbb{C}^2)$, $F\ne G$, then given an
$H\in {\rm et}(\mathbb{C}^2)$ for which $F\in R_H({\rm et}(\mathbb{C}^2))$ we get $F=M\circ H$ for some $M\in {\rm et}(\mathbb{C}^2)$. Since
$F$ is injective, it follows that $H$ is injective. Hence we deduce that $H\in R_H({\rm Aut}(\mathbb{C}^2))$ and so $R_H({\rm et}(\mathbb{C}^2))={\rm et}(\mathbb{C}^2)$.
Likewise, the only open set (in $\tau_R$) that contains $G$ is ${\rm et}(\mathbb{C}^2)$, for also $G\in {\rm Aut}(\mathbb{C}^2)$. $\qed $ \\
\\
\noindent
We naturally ask if the subspace ${\rm et}(\mathbb{C}^2)-{\rm Aut}(\mathbb{C}^2)$ of $({\rm et}(\mathbb{C}^2),\tau_R\}$ is Hausdorff. Also
here the answer is negative:

\begin{proposition}
The subspace ${\rm et}(\mathbb{C}^2)-{\rm Aut}(\mathbb{C}^2)$ of $({\rm et}(\mathbb{C}^2),\tau_R\}$ is not Hausdorff.
\end{proposition}
\noindent
{\bf Proof.} \\
Let $F\in {\rm et}(\mathbb{C}^2)-{\rm Aut}(\mathbb{C}^2)$ be a prime (Theorem 4.5(1)). We will show that $\tau_R$ can not separate 
the points $F$ and $F\circ F$. For if $F\in R_H({\rm et}(\mathbb{C}^2))$, then $R_H({\rm et}(\mathbb{C}^2))=R_F({\rm et}(\mathbb{C}^2))$. But then
$G=F\circ F\in R_F({\rm et}(\mathbb{C}^2))=R_H({\rm et}(\mathbb{C}^2))$. Thus if $G\in R_L({\rm et}(\mathbb{C}^2))$, then
$R_H({\rm et}(\mathbb{C}^2))\cap R_L({\rm et}(\mathbb{C}^2))\ne\emptyset$ for this intersection contains $G$. $\qed $ \\

\section{Metric structures on ${\rm et}(\mathbb{C}^2)$ that we would like to have}

Any $F\in {\rm et}(\mathbb{C}^2)$ is determined by its sets of coefficients (those of $P$ and those of $Q$). We can order
the sequences of the coefficients in ascending degree order and within each homogeneous part lexicographically ($X>Y$). In other
words if $P(X,Y)=\sum_{1\le i+j\le N=\deg P} a_{ij}X^iY^j$ and $Q(X,Y)=\sum_{1\le i+j\le M=\deg Q} b_{ij}X^iY^j$ then those
two sequences are: 
$$
(a_{10},a_{01},a_{20},a_{11},a_{02},\ldots,a_{0N})\,\,\,{\rm and}\,\,\,(b_{10},b_{01},b_{20},b_{11},b_{02},\ldots,b_{0M})
$$ 
where $a_{0N}\cdot b_{0M}\ne 0$ (condition (2) in the definition of mappings in ${\rm et}(\mathbb{C}^2)$ that was given before Definition 1.1) and where 
the coefficients satisfy the Jacobian Condition. The Jacobian Condition is expressible by an infinite set of polynomial quadratic equations, all of which are homogeneous except for just one equation, namely $a_{10}b_{01}-a_{01}b_{10}=1$ which is still quadratic but not homogeneous. If we drop the 
open condition $a_{0N}\cdot b_{0M}\ne 0$ we get an infinite dimensional affine algebraic variety. The structure of
this space decomposes according to part (4) of Proposition 4.13 into a fractal like decomposition.
There exists an infinite index set $I$ and a family of \'etale mappings $\{ F_i\,|\,i\in I\}\subseteq {\rm et}(\mathbb{C}^2)$ such that
$$
{\rm et}(\mathbb{C}^2)={\rm Aut}(\mathbb{C}^2)\cup\bigcup_{i\in I}R_{F_i}({\rm et}(\mathbb{C}^2)),
$$
so that $\forall\,i\in I$, ${\rm Aut}(\mathbb{C}^2)\cap R_{F_i}({\rm et}(\mathbb{C}^2))=\emptyset$, and $\forall\,i,j\in I$, if $i\ne j$ then
$F_i\not\in R_{F_j}({\rm et}(\mathbb{C}^2))$ and $F_j\not\in R_{F_i}({\rm et}(\mathbb{C}^2))$, and $R_{F_i}({\rm et}(\mathbb{C}^2))$ is
homeomorphic to $R_{F_j}({\rm et}(\mathbb{C}^2))$, and both are homeomorphic to ${\rm et}(\mathbb{C}^2)$. \\
We will try to apply fractal geometric tools to this structure (or a similar one - where the place of the right (composition) mappings
$R_F$ will be taken by the left (composition) mappings $L_F$). A crucial step will be to define useful metrics and the corresponding Hausdorff measures on
${\rm et}(\mathbb{C}^2)$ in order to obtain some (fractional) dimension computations or estimates of this space. Thus from now
on we can identify ${\rm et}(\mathbb{C}^2)$ with its image under the above embedding (via the coefficients vectors of the polynomials),
$$
(\sum_{1\le i+j\le N} a_{ij}X^iY^j,\sum_{1\le i+j\le M} b_{ij}X^iY^j)\rightarrow
$$
$$
\rightarrow((a_{10},a_{01},a_{20},a_{11},a_{02},\ldots,a_{0N}),(b_{10},b_{01},b_{20},b_{11},b_{02},\ldots,b_{0M})),
$$
into the space $\mathbb{C}^{\aleph_0}\times\mathbb{C}^{\aleph_0}$. In fact the image is contained in the Cartesian product of the finite
sequences over $\mathbb{C}$ by itself (where we think of a finite sequence as an infinite sequence which is eventually composed of zeros).
Suppose that we have a metric $\rho\,:\,{\rm et}(\mathbb{C}^2)\times {\rm et}(\mathbb{C}^2)\rightarrow \mathbb{R}^+\cup\{0\}$.
Let $F\in {\rm et}(\mathbb{C}^2)$. Then $R_F({\rm et}(\mathbb{C}^2))=\{ G\circ F\,|\,G\in {\rm et}(\mathbb{C}^2)\}$ is a metric
subspace of $({\rm et}(\mathbb{C}^2),\rho)$, by restricting $\rho$ to $R_F({\rm et}(\mathbb{C}^2))$. Under the natural topology
the space ${\rm et}(\mathbb{C}^2)$ is homeomorphic to its subspace $R_F({\rm et}(\mathbb{C}^2))$. The homeomorphism being,
$$
R_F\,:\,{\rm et}(\mathbb{C}^2)\rightarrow R_F({\rm et}(\mathbb{C}^2)),\,\,\,R_F(G)=G\circ F.
$$
This homeomorphism need not be a $\rho$-isometry, even if the topology on ${\rm et}(\mathbb{C}^2)$ is identical to the $\rho$-metric topology.
We would like to have such a metric $\rho$ that will give us a good control on $\rho(R_F(G_1),R_F(G_2))=\rho(G_1\circ F,G_2\circ F)$ in 
terms of $\rho(G_1,G_2)$.
\begin{remark}
A natural topology on ${\rm et}(\mathbb{C}^2)$ is the so called compact-open topology. Just like for analytic functions of a single complex
variable we naturally look at sequences $F_n\in {\rm et}(\mathbb{C}^2)$ that locally uniformly converge to a polynomial mapping. This
means that for every compact (in the strong topology) $K\subseteq\mathbb{C}^2$ we have $F_n|_K\rightarrow_{n\rightarrow\infty} F|_K$ uniformly.
\end{remark}
\noindent
How to construct such a metric (that will be sensitive to compositions in ${\rm et}(\mathbb{C}^2)$)? The idea is straight forward. 
The mappings in ${\rm et}(\mathbb{C}^2)$ all satisfy the Jacobian Condition. Thus, geometrically, these are polynomial mappings 
$\mathbb{R}^4\rightarrow\mathbb{R}^4$ (in the four dimensional space over the reals) that locally preserve volume. This is a crucial
geometric property and we want our metric to capture this property. But we will see that (at least according to our constructions) the right (composition) mapping
$R_F$ and the left (composition) mapping $L_F$ are very different! There is no symmetry between those two and in fact it will turn out that the
left (composition) mappings $L_F$ are the correct to work with. So our plane is the following: we will outline the constructions of
the metrics on ${\rm et}(\mathbb{C}^2)$ that are sensitive to compositions of \'etale mappings. After that we will switch the results we
developed so far from the right (composition) mappings setting to the left (composition) mappings. After this will be done, we will have an efficient machinery
that will tie the metric space structure on ${\rm et}(\mathbb{C}^2)$ to a compatible approximate fractal structure. This will serve us to conclude
non trivial geometrical results on the algebro-geometric structure ${\rm et}(\mathbb{C}^2)$.

Composition of mappings is not simply a non-commutative binary operation. There is a deep difference between the two operands,
the left and the right. Consider two mappings $f,\,g:\,X\rightarrow X$. When we form their composition $h=f\circ g$, then if $g$ is non
injective so is $h$. If $f$ is non surjective, then so is $h$. The two examples we investigated, i.e. ${\rm elh}(\mathbb{C})$ and
${\rm et}(\mathbb{C}^2)$ show that it is much easier to prove that $R_g$ is injective than to show that $L_f$ is. In fact
for the entire single variable case, ${\rm elh}(\mathbb{C})$ it turns out that $L_f$ is not injective. However, it is "almost" injective in the
sense that we can single out the exceptional cases which form a small sub-family of ${\rm elh}(\mathbb{C})$. The reason for the
non injectivity originates in the existence of the periodic functions in ${\rm elh}(\mathbb{C})$. This kind of an obstacle to injectivity is void
for the algebraic \'etale case ${\rm et}(\mathbb{C}^2)$. Indeed in this case $L_f$ turns out to be injective. But it is highly non trivial
to prove that. There are some algebro-geometric reasons that explain this difficulty.

\section{The metric spaces $({\rm et}(\mathbb{C}^2),\rho_D)$}
We will need a special kind of four (real) dimensional subsets of $\mathbb{R}^4$. These will serve us to construct suitable metric structures 
on ${\rm et}(\mathbb{C}^2)$. We will describe the construction step by step, leaving occasionally some details for later stages in order
not to brake the line of reasoning.

\begin{remark}
We will use Definition 1.1 of a characteristic set $D$ for the family ${\rm et}(\mathbb{C}^2)$ and of the corresponding mapping
$\rho_D\,:\,{\rm et}(\mathbb{C}^2)\times {\rm et}(\mathbb{C}^2)\rightarrow\mathbb{R}^+\cup\{0\}$ which turns out to be a metric
on ${\rm et}(\mathbb{C}^2)$, by Proposition 1.3.
\end{remark}

\begin{remark}
It is not clear how to construct an open subset $D$ of $\mathbb{C}^2$ that will satisfy the three properties
that are required in definition 1.1. We will postpone for a while the demonstration that such open sets exist.
\end{remark}

\begin{remark}
We recall that according to Proposition 1.3 we have: $\rho_D$ is a metric on ${\rm et}(\mathbb{C}^2)$.
\end{remark}

So far we thought of the volume of $G_1(D)\Delta G_2(D)$ as the volume of the set which is the symmetric difference between
the $G_1$ image and the $G_2$ image of the open set $D$. However, the mappings $G_1$ and $G_2$ are \'etale and in particular need
not be injective. We will take into the volume computation the multiplicities of $G_1$ and of $G_2$. By Theorem 3 on page 39 of \cite{bm}
we have the following: Given $F\,:\,\mathbb{C}^n\rightarrow\mathbb{C}^n$ we define $\tilde{F}=({\rm Re}\,F_1,{\rm Im}\,F_1,\ldots,
{\rm Re}\,F_n,{\rm Im}\,F_n)\,;\,\mathbb{R}^{2n}\rightarrow\mathbb{R}^{2n}$. Then $\det J_{\tilde{F}}=|\det J_F|^2$. Thus the Jacobian 
Condition, $\det J_F\equiv 1$ implies that $\det J_{\tilde{F}}\equiv 1$. So the real mapping $\tilde{F}$ locally preserves the usual volume form.
In order to take into an account the multiplicities of the \'etale mappings $G_1$ and $G_2$ when computing the volume of the symmetric
difference $G_1(D)\Delta G_2(D)$ we had to do the following. For any $G\in {\rm et}(\mathbb{C}^2)$ instead of computing,
$$
\int\int\int\int_D\left(\det J_{\tilde{G}}\cdot dV\right)=\int\int\int\int_{D}dV,
$$
we compute
$$
\int\int\int\int_{\tilde{G}(D)}dX_1 dX_2 dY_1 dY_2\,\,\,\,\,{\rm where}\,\,X=X_1+iX_2,\,Y=Y_1+iY_2.
$$
For every $j=1,2,\ldots,d_G$ we denote by $D_j$ that subset of $D$ such that for each point of $D_j$ there are exactly $j$ points 
of $D$ that are mapped by $G$ to the same image of that point. In other words, $D_j=\{\alpha\in D\,|\,|\tilde{G}^{-1}(\tilde{G}(\alpha))\cap D|=j\}$.
We assume that $D$ is large enough so that $\forall\,j=1,\ldots,d_G$ we have $D_j\ne\emptyset$.
For our \'etale mappings it is well known that if $j<d_G$ then $\dim D_j<\dim D$ so the volume these $D_j$'s contribution equals to 0. 
The dimension claim follows by the well known fact that the size of a generic fiber $|G^{-1}(x)|$ equals to $d_G$ and that $d_G$ is also the
maximal size of any of the fibers of $G$.
However, for the sake of treating more general families of mappings we denote by ${\rm vol}(D_j)$ the volume of the set $D_j$.
Then $D$ has a partition into exactly $j$ subsets of equal volume. The volume of each such a set is ${\rm vol}(D_j)/j$ and each such a set has
exactly one of the $j$ points in $\tilde{G}^{-1}(\tilde{G}(\alpha))\cap D$ for each $\alpha\in D_j$. We note that ${\rm vol}(\tilde{G}(D_j))={\rm vol}(D_j)/j$
by the Jacobian Condition. Thus the volume with the multiplicity of $\tilde{G}$ taken into account is given by:
$$
{\rm vol}(\tilde{G}(D))+\sum_{j=2}^{d_G}(j-1)\cdot\frac{{\rm vol}(D_j)}{j}=
{\rm vol}(\tilde{G}(D))+\sum_{j=2}^{d_G}(j-1)\cdot{\rm vol}(\tilde{G}(D_j)).
$$
We note that $\tilde{G}(D)=\bigcup_{j=1}^{d_G}\tilde{G}(D_j)$ is a partition, so ${\rm vol}(\tilde{G}(D))=\sum_{j=1}^{d_G}{\rm vol}(\tilde{G}(D_j))$. 
Hence we can express the desired volume by
$$
{\rm vol}(\tilde{G}(D))+\sum_{j=2}^{d_G}(j-1)\cdot{\rm vol}(\tilde{G}(D_j))=\sum_{j=1}^{d_G}j\cdot {\rm vol}(\tilde{G}(D_j)).
$$
We note that this equals to $\sum_{j=1}^{d_G}{\rm vol}(D_j)$ and since $D=\bigcup_{j=1}^{d_G} D_j$ is a partition we have
${\rm vol}(D)=\sum_{j=1}^{d_G} {\rm vol}(D_j)$. As expected, the volume computation that takes into account the multiplicity of $G$ is in general
larger than the geometric volume ${\rm vol}(\tilde{G}(D))$. The access can be expressed in several forms:
$$
{\rm vol}(D)-{\rm vol}(\tilde{G}(D))=\sum_{j=2}^{d_G}(j-1)\cdot {\rm vol}(\tilde{G}(D_j))=\sum_{j=2}^{d_G}\left(1-\frac{1}{j}\right){\rm vol}(D_j).
$$
Coming back to the computation of the metric distance $\rho_D(G_1,G_2)={\rm the\,volume\,of}\,\,G_1(D)\Delta G_2(D)$ we compute the volume
of $G_1(D)-G_2(D)$ with the multiplicity of $G_1$ while the volume of $G_2(D)-G_1(D)$ is computed with the multiplicity of $G_2$.

\section{Characteristic sets of families of holomorphic local homeomorphisms $\mathbb{C}^2\rightarrow\mathbb{C}^2$ (see \cite{rp6})}
In this section we prove the existence of sets $D$ that satisfy the three properties that are required in definition 1.1.
The third property will turn out to be the tricky one. 

\begin{definition}
let $\Gamma$ be a family of holomorphic local homeomorphisms $F\,:\,\mathbb{C}^2\rightarrow\mathbb{C}^2$.
A subset $D\subseteq\mathbb{C}^2$ is called a characteristic set of $\Gamma$ if it satisfies the
following condition: $\forall\,F_1,F_2\in\Gamma$, $F_1(D)=F_2(D)\Leftrightarrow F_1=F_2$.
\end{definition}
\noindent
We start by recalling the well known rigidity property of holomorphic functions in one complex variable. Also known as the permanence principle,
or the identity theorem. The identity theorem  for analytic functions of one complex variable says that if $D\subseteq\mathbb{C}$ is a domain
(an open and a connected set) and if $E$ is a subset of $D$ that has a non-isolated point and if $f(z)$ is an analytic function defined on $D$ and 
vanishing on $E$, then $f(z)=0$ for all $z\in D$.

There is an identity theorem for analytic functions in several complex variables, but for more than one variable the above statement is false.
One possible correct statement is as follows:"Let $f(z)$ and $g(z)$ be holomorphic functions in a domain $D$ in $\mathbb{C}^n$. If $f(z)=g(z)$ for all $z$ 
in a non-empty open set $\delta$ in $D$, then $f(z)=g(z)$ in $D$. Hence, analytic continuation of holomorphic functions in several complex variables
can be performed as in the case of one complex variable. Contrary to the case of one complex variable, the zero set of a holomorphic function in 
a domain $D\subseteq\mathbb{C}^n$, $n\ge 2$, contains no isolated points. Thus even if $f(z)=g(z)$ in a set with accumulation points in $D$, it
does not necessarily follow that $f(z)=g(z)$ in $D$. For example, in $\mathbb{C}^2$ with variables $z$ and $w$ we can take $f(z,w)=z$ and
$g(z,w)=z^2$." (Chapter 1, page 16 in \cite{tn}).

In spite of the above standard identity theorem for $n\ge 2$ complex variables, that requires a thick set $E$ (i.e. an open set) on which $f(z)=g(z)$
one can do much better. Let us start with the following. Let $F(Z,W)$ be an entire function of two complex variables $Z$ and $W$. Let us define a
subset $E$ of $\mathbb{C}^2$ as follows. We take a convergent sequence $\{Z_k\}_{k=1}^{\infty}$ of different numbers. Thus $\lim Z_k=a$ and
$j\ne k\Rightarrow Z_j\ne Z_k$. For each $k$, let $\{W_j^{(k)}\}_{j=1}^{\infty}$ be a convergent sequence of different numbers, such that
their limit is $\lim_{j\rightarrow\infty} W_j^{(k)}=Z_k^{'}$. We define $E=\{(Z_k,W_j^{(k)})\,|\,j,k=1,2,3,\ldots\}$. Now we have,

\begin{proposition}
If the entire function $F(Z,W)$ vanishes on $E$, i.e. $F(Z_k,W_j^{(k)})=0$ for $j,k=1,2,3,\ldots $, then $F(Z,W)\equiv 0$ is the zero function.
\end{proposition}

\begin{remark}
We note that $E$ is a thin set, in fact a countable set. Even the closure $\overline{E}$ is thin.
\end{remark}
\noindent
{\bf A Proof of Proposition 7.2.} \\
Since $F(Z,W)$ is an entire function, it can be represented as a convergent power series centered at $(0,0)$ with an infinite
radius of convergence. We can sum the terms in the order we please. Let us write $F(Z,W)$ as a power series in $W$ with
coefficients that are entire functions in $Z$. Thus we have, $F(Z,W)=\sum_{k=0}^{\infty} a_k(Z)W^k$, where for each $k=0,1,2,\ldots $,
$a_k(Z)$ is an entire function in the variable $Z$. For a fixed $l\in\{1,2,3,\ldots\}$ we have by our assumptions the following,
$F(Z_l,W_j^{(l)})=0$ for $j=1,2,3,\ldots $. But $\lim_{j\rightarrow\infty}W_j^{(l)}=Z_l^{'}$ so that $g_l(W)=F(Z_l,W)$ is an
entire function of the single variable $W$, which vanishes on a convergent sequence $\{W_j^{(l)}\}_{j=1}^{\infty}$. By the identity
theorem of one complex variable we deduce that $g_l(W)\equiv 0$, the zero function. Since, $g_l(W)=\sum_{k=0}^{\infty} a_k(Z_l)W^k$ it
follows that the Maclaurin coefficients $a_k(Z_l)$, $k=0,1,2,\ldots $ vanish. Now, this is valid for each $l$, and $\lim Z_l=a$ converges.
Since each $a_k(Z)$ is an entire function which vanishes on a convergent sequence $\{ Z_l\}_{l=1}^{\infty}$ it follows, once again, by
the identity theorem in one complex variable, that $a_k(Z)\equiv 0$, the zero function, $k=0,1,2,\ldots $. Hence we conclude that
$F(Z,W)=\sum_{k=0}^{\infty} a_k(Z)W^k\equiv 0$. $\qed $ \\
\\
This type of elementary arguments that was used to construct a thin set $E$ for identity purpose, is not new. For example: \\
"Theorem. Let $D\subseteq\mathbb{C}$ be a domain, and let $E$ be a subset of $D$ that has a non-isolated point. Let $F(Z,W)$ be a function 
defined for $Z,W\in D$ such that $F(Z,W)$ is analytic in $Z$ for each fixed $W\in D$ and analytic in $W$ for each fixed $Z\in D$.
If $F(Z,W)=0$ whenever $Z$ and $W$ both belong to $E$, then $F(Z,W)=0$ for all $Z,W\in D$.", \cite{g}.

Advancing along the lines of the construction of the thin set in Proposition 7.2 we note that if $\{Z_k\}_{k=1}^{\infty}$ is a
sequence of different numbers that converges to $\lim Z_k=a$, and if for each $k=1,2,3,\ldots $ there is a straight line segment 
$[\alpha_k,\beta_k]$ of $W$'s such that two entire functions $F(Z,W)$ and $G(Z,W)$ agree on the union (a countable union) of
the segments $\{Z_k\}\times [\alpha_k,\beta_k]$, i.e. $F(Z_k,W)=G(Z_k,W)$, $\forall\,W\in[\alpha_k,\beta_k]$, then $F(Z,W)\equiv G(Z,W)$,
$\forall\,(Z,W)\in\mathbb{C}^2$.

We now will construct characteristic sets of families $\Gamma$ of holomorphic local homeomorphisms $F:\,\mathbb{C}^2\rightarrow\mathbb{C}^2$.

\begin{definition}
Let $m$ be a natural number and $\alpha\in\mathbb{C}^2$. An $m$-star at $\alpha$ is the union of $m$ line segments, so that any pair intersects
in $\alpha $.
\end{definition}

\begin{definition}
Let $l$ be a line segment and let $\{\alpha_k\}$ be a countable dense subset of $l$. Let $\{n_k\}$ be a sequence of different natural numbers
and $\forall\,k$, let $S_{n_k}$ be an $n_k$-star at $\alpha_k$ such that one of the star's segments lies on $l$, and such that $\forall\,k_1\ne k_2$,
$\tilde{S}_{n_{k_1}}\cap\tilde{S}_{n_{k_2}}=\emptyset$. Here we denoted $\tilde{S}=S-l$. Moreover, we group the stars in bundles of, say 5, thus
getting the sequence of star bundles:
$$
\{S_{n_1},S_{n_2},S_{n_3},S_{n_4},S_{n_5}\},\,\{S_{n_6},\ldots,S_{n_{10}}\},\ldots,\{S_{n_{5j+1}},\ldots,S_{n_{5j+5}}\},\ldots
$$
and for each bundle of five we take the maximal length of its rays to be at most $1/10$ the length of the maximal length of the previous
bundle. We define,
$$
l_0^{\{n_k\}}=l\cup\bigcup_{k=1}^{\infty}S_{n_k}.
$$
Let $\{Z_k\}_{k=1}^{\infty}$ be a sequence of different complex numbers that converges to $\lim Z_k=a$. Let $\{\{n_j^{(k)}\}_{j=1}^{\infty}\}_{k=1}^{\infty}$
be a partition of the natural numbers, $\mathbb{Z}^+$. In fact all we need is the disjointness, i.e. $k_1\ne k_2\Rightarrow
\{n_j^{(k_1)}\}_{j=1}^{\infty}\cap\{n_j^{(k_2)}\}_{j=1}^{\infty}=\emptyset$. Let us consider the stared segments
$$
\{l_0^{\{n_j^{(k)}\}_{j=1}^{\infty}}\,|\,k=1,2,3,\ldots\}
$$
and define the following countable union of stared segments in $\mathbb{C}^2$:
$$
\bigcup_{k=1}^{\infty}\{Z_k\}\times l_0^{\{n_j^{(k)}\}_{j=1}^{\infty}},
$$
where we assume that the lengths of the star rays were chosen to satisfy disjointness in $\mathbb{C}^2$, namely:
$$
k_1\ne k_2\Rightarrow\{Z_{k_1}\}\times l_0^{\{n_j^{(k_1)}\}_{j=1}^{\infty}}\cap\{Z_{k_2}\}\times l_0^{\{n_j^{(k_2)}\}_{j=1}^{\infty}}=\emptyset.
$$
We let,
$$
E=\bigcup_{k=1}^{\infty}\{Z_k\}\times l_0^{\{n_j^{(k)}\}_{j=1}^{\infty}},
$$
or if we need a closed (compact) set, the closure of this union.
\end{definition}

\begin{proposition}
Let $\Gamma$ be any family of entire holomorphic local homeomorphisms $F:\,\mathbb{C}^2\rightarrow\mathbb{C}^2$. Then $E$ is a characteristic 
set of $\Gamma$.
\end{proposition}
\noindent
{\bf Proof.} \\
Let $F_1,F_2\in\Gamma$ satisfy $F_1(E)=F_2(E)$. Then each stared line segment,
$$
\{Z_k\}\times l_0^{\{n_j^{(k)}\}_{j=1}^{\infty}},
$$
must be mapped onto a curve,
$$
F_1(\{Z_k\}\times l_0^{\{n_j^{(k)}\}_{j=1}^{\infty}})=F_2(\{Z_k\}\times l_0^{\{n_j^{(k)}\}_{j=1}^{\infty}})
$$
and each $n_j^{(k)}$-star on $l$, $S_{n_j^{(k)}}$ is mapped onto a holomorphic $n_j^{(k)}$-star,
$$
F_1(\{Z_k\}\times S_{n_j^{(k)}})=F_2(\{Z_k\}\times S_{n_j^{(k)}}).
$$
This is because the valence sequences of the stars 
$$
\{\{n_j^{(k)}\}_{j=1}^{\infty}\}_{k=1}^{\infty},
$$
are pairwise disjoint natural numbers, and $F_1$, $F_2$ are local homeomorphisms and hence preserve the star
valencies $n_j^{(k)}$. The centers of the holomorphic stars,
$$
\{F_1(\alpha_{n_j^{(k)}})\}=\{F_2(\alpha_{n_j^{(k)}})\},
$$
form a countable and a dense subset of the curves $F_1(\{Z_k\}\times l)=F_2(\{Z_k\}\times l)$. By continuity this
implies that the restrictions,
$$
F_1|_{\{Z_k\}\times l}\,\,\,\,\,{\rm and}\,\,\,\,\,F_2|_{\{Z_k\}\times l},
$$
coincide. Since $F_1$ and $F_2$ are holomorphic, this implies by Proposition 7.2 (which is a variant of the identity
theorem for entire functions $\mathbb{C}^2\rightarrow\mathbb{C}^2$) that $F_1\equiv F_2$. $\qed $ \\

\begin{remark}
Proposition 7.6 holds true for any rigid family of local homeomorphisms. Rigidity here means that 
$$
F_1|_{\{Z_k\}\times l}=F_2|_{\{Z_k\}\times l}\Leftrightarrow F_1\equiv F_2.
$$ 
So the proposition holds true for holomorphic mappings, for harmonic mappings and also for ${\rm et}(\mathbb{C}^2)$.
\end{remark}
\noindent
We recall that definition 1.1 required also two additional topological properties, namely the open set $D$ should satisfy ${\rm int}(\overline{D})=D$,
$\overline{D}$ is compact (all in the strong topology). These automatically exclude the set $E$ that was constructed in definition 7.5. However, we 
can modify this construction to get at least an open set. 

\begin{proposition}
Let $\Gamma$ be any family of holomorphic local homeomorphisms $F\,:\,\mathbb{C}^2\rightarrow\mathbb{C}^2$. Let $U$ be any open subset of $\mathbb{C}^2$ with a
smooth boundary that contains the compact $E$. Then the open set $U-E$ is a characteristic set of $\Gamma$.
\end{proposition}
\noindent
{\bf Proof.} \\
Since $E$ can not be mapped in the smooth $\partial U$ by an holomorphic local homeomorphism, we have for any $F_1,F_2\in\Gamma$ for which
$F_1(U-E)=F_2(U-E)$ that also $F_1(E)=F_2(E)$. Now the result follows by Proposition 7.6. $\qed $ \\

\begin{remark}
We note that if $\overline{U}$ is a compact then $U-E$ satisfies, at least the requirement $\overline{U-E}$ is compact. However, the "no slit" 
condition (condition (1) in Definition 1.1) ${\rm int}(\overline{U-E})={\rm int}(\overline{U})\ne U-E$ fails.
\end{remark}
\noindent
Now that we gained some experience with the topological construction of $E$ we are going to make one more step and fix its
shortcomings that were mentioned above. We need to construct a domain $D$ of $\mathbb{C}^2$ which has the following three properties: \\
1) ${\rm int}(\overline{D})=D$ relative to the complex topology. \\
2) $\overline{D}$ is a compact subset of $\mathbb{C}^2$ relative to the strong topology. \\
3) $\forall\,G_1,G_2\in {\rm et}(\mathbb{C}^2)$, $G_1(D)=G_2(D)\Leftrightarrow G_1\equiv G_2$. \\
\\
(The complex topology and the strong topology are the same). Our construction will be a modification of the construction
of the domain that was constructed in Proposition 7.8. We start by modifying the notion of an $m$-star that was introduced
in Definition 7.4.

\begin{definition}
Let $m$ be a natural number and $\alpha\in\mathbb{C}^n$. A thick $m$-star at $\alpha$ is a union of $2m$ triangles, so
that any pair intersect exactly at one vertex, and this vertex (that is common to all the $2m$ triangles) is $\alpha$.
\end{definition}

\begin{definition}
Let $E$ be the construction of Definition 7.5 that uses thick $m$-stars.
\end{definition}

\begin{proposition}
Let $\Gamma$ be any family of holomorphic local homeomorphisms $F\,:\,\mathbb{C}^2\rightarrow\mathbb{C}^2$. Then $E$ is
a characteristic set of $\Gamma$.
\end{proposition}
\noindent
{\bf Proof.} \\
The proof is the same word-by-word as that of Proposition 7.6 where we replace $k$-star $S_k$ by thick $k$-star $S_k$. $\qed $ \\
\\
We finally obtain our construction.
\begin{proposition}
Let $\Gamma$ be any family of holomorphic local homeomorphisms $F\,:\,\mathbb{C}^2\rightarrow\mathbb{C}^2$. Let $B(0,R)$ be
an open ball centered at $0$ with a radius $R$ large enough so that $E\subset B(0,R)$ (where $E$ is the set in Proposition 7.12).
Then the domain $D=B(0,R)-E$ is a characteristic set of $\Gamma$.
\end{proposition}
\noindent
{\bf Proof.} \\
The proof is the same as that of Proposition 7.8 where we replace $k$-star $S_k$ by thick $k$-star $S_k$. $\qed $

\section{Switching to the left (composition) mapping $L_F$ (see \cite{rp6})}
As was explained in Section 5 we would like our natural mappings: The right (composition) mapping $R_F$, or the left (composition) mapping $L_F$ to be say
bi-Lipschitz with respect to the metric $\rho_D$ (that reflects the fact that our mappings, in ${\rm et}(\mathbb{C}^2)$ satisfy
the Jacobian Condition). A short reflection shows that the right (composition)  mapping $R_F$ need not be $\rho_D$ bi-Lipschitz. 
The situation is completely different when we replace the right (composition) mapping, $R_F$ by the left (composition) mapping, $L_F$.
A promising indication is the following,

\begin{proposition}
$\forall\,F\in {\rm Aut}(\mathbb{C}^2)$ the mapping $L_F$ is an isometry of the metric space $({\rm et}(\mathbb{C}^2),\rho_D)$.
\end{proposition}
\noindent
{\bf Proof.} \\
For any two mappings $G_1$ and $G_2$ in ${\rm et}(\mathbb{C}^2)$ we need to compare $\rho_D(G_1,G_2)$ with $\rho_D(F\circ G_1,F\circ G_2)$.
We have (using our assumption on $F$), 
$$
(F\circ G_1)(D)\Delta (F\circ G_2)(D)=F\left(G_1(D)\Delta G_2(D)\right).
$$
Since $F$ is also (globally) volume preserving we have,
$$
{\rm the\,volume\,of}\,F\left(G_1(D)\Delta G_2(D)\right)={\rm the\,volume\,of}\,\left(G_1(D)\Delta G_2(D)\right).
$$
This proves that $\rho_D(G_1,G_2)=\rho_D(F\circ G_1,F\circ G_2)$. $\qed $ \\

We now drop the restrictive assumption that $F\in {\rm Aut}(\mathbb{C}^2)$. We briefly describe few results from \cite{rp6} that will be useful
for us here. From this point on we assume that the characteristic sets $D$ are almost $2n$-Euclidean balls except for a very small part
of the boundary which contains the countable and dense (at that portion only) set of the centers of the thick $m$-stars. We refer to our
construction which is applied in Proposition 7.13. This assumption enables us to use the geometric methods of entire mappings (in particular the 
theory of maximal domains of entire mappings).
 
Thus we merely have $F\in {\rm et}(\mathbb{C}^2)$ and we still want to compare $\rho_D(G_1,G_2)$ with $\rho_D(F\circ G_1,F\circ G_2)$, 
for any pair $G_1,G_2\in {\rm et}(\mathbb{C}^2)$. We only know that $F$ is a local diffeomorphism of $\mathbb{C}^2$ and (by the Jacobian Condition) that it preserves (locally) the volume. In this case the geometrical degree of $F$, $d_F$ can be larger than $1$. We have the identity $d_F=|F^{-1}(\{(a,b)\})|$ which holds generically (in the sense of the Zariski topology) in $(a,b)\in\mathbb{C}^2$. Hence the (complex) dimension of the set $\{(a,b)\in\mathbb{C}^2\,|\,|F^{-1}(a,b)|<d_F\}$ is at most $1$. The Jacobian Condition $\det J_F\equiv 1$ implies (as we noticed before) that $F$ preserves volume taking into account the multiplicity. The multiplicity is a result of the possibility that $F$ is not injective and hence the deformation of the characteristic set $D$ by $F$ convolves (i.e. might overlap at certain locations). However, this overlapping is bounded above by $d_F$. So if $A\subseteq\mathbb{C}^2$ is a measurable subset of $\mathbb{C}^2$ and we compare the volume of $A$ with the volume of its image $F(A)$, then,
$$
{\rm the\,volume\,of}\,F(A)\le\,{\rm the\,volume\,of}\,A\le\,d_F\cdot\{{\rm the\,volume\,of}\,F(A)\}.
$$
This can be rewritten as follows,
$$
\frac{1}{d_F}\cdot\{{\rm the\,volume\,of}\,A\}\le\,{\rm the\,volume\,of}\,F(A)\le\,{\rm the\,volume\,of}\,A.
$$
This is the place to emphasize also the following conclusion (that follows by the generic identity $d_F=|F^{-1}(\{(a,b)\})|$), namely
$$
\lim_{A\rightarrow\mathbb{C}^2}\frac{{\rm the\,volume\,of}\,F(A)}{{\rm the\,volume\,of}\,A}=\frac{1}{d_F},
$$
provided that the set $A$ tends to cover the whole of the complex space $\mathbb{C}^2$ in an appropriate manner. To better understand
why the quotient tends to the lower limit $1/d_F$ rather than to any number in the interval $[1/d_F,1]$ (if at all) we recall that
our mapping belongs to ${\rm et}(\mathbb{C}^2)$ and so is a polynomial \'etale mapping. So any point $(a,b)\in\mathbb{C}^2$ for which
$|F^{-1}(a,b)|<d_F$ is an asymptotic value of $F$ and hence belongs to the curve $A(F)$ the asymptotic variety of $F$. In other words
the identity $d_F=|F^{-1}(a,b)|$ is satisfied exactly on the semi algebraic set $\mathbb{C}^2-A(F)$ which is the complement of an algebraic
curve. Here is Theorem 34 from \cite{rp6}.

\begin{theorem}
Let $F,G_1,G_2\in {\rm et}(\mathbb{C}^2)$. Then we have: \\
(i) $\rho_D(F\circ G_1,F\circ G_2)\le \rho_D(G_1,G_2)$. \\
(ii) Suppose that $D$ is a family of characteristic sets of ${\rm et}(\mathbb{C}^2)$ such that $D\rightarrow\mathbb{C}^2$, then
$\forall\,\epsilon>0$ we have,
$$
\left(\frac{1}{d_F}-\epsilon\right)\cdot\rho_D(G_1,G_2)\le\rho_D(F\circ G_1,F\circ G_2)
$$
for $D$ large enough. \\
(iii)  Under the assumptions in (ii) we have:
$$
\lim_{D\rightarrow\mathbb{C}^2}\frac{\rho_D(F\circ G_1,F\circ G_2)}{\rho_D(G_1,G_2)}=\frac{1}{d_F}.
$$
In particular, the left (composition) mapping $L_F: {\rm et}(\mathbb{C}^2)\rightarrow {\rm et}(\mathbb{C}^2)$, $L_F(G)=F\circ G$,
is a bi-Lipschitz self-mapping of the metric space $({\rm et}(\mathbb{C}^2),\rho_D)$ with the constants $1/d_F\le 1$.
\end{theorem}
\noindent
{\bf Proof.} \\
We briefly sketch the proof of the theorem: \\
(i) $x\in (F\circ G_1)(D)\Delta (F\circ G_2)(D)\Rightarrow \exists\,y\in G_j(D),\,j=1\,{\rm or}\,2$ such that $x=F(y)$ and
$x\not\in (F\circ G_{3-j})(D)$. By $x\not\in (F\circ G_{3-j})(D)$ it follows that $y\not\in G_{3-j}(D)$ and so $y\in G_1(D)\Delta G_2(D)$
and $x=F(y)\in F(G_1(D)\Delta G_2(D))$. Hence $(F\circ G_1)(D)\Delta (F\circ G_2)(D)\subseteq F(G_1(D)\Delta G_2(D))$, so
${\rm vol}((F\circ G_1)(D)\Delta (F\circ G_2)(D))\le{\rm vol}(F(G_1(D)\Delta G_2(D))$, and finally $\rho_D(F\circ G_1,F\circ G_2)\le
\rho_D(G_1,G_2)$. \\
(ii) and (iii). Here the proof is not just set theoretic. We will elaborate more in the remark that follows this proof.
We recall that $F,G_1,G_2\in {\rm et}(\mathbb{C}^2)$. This implies that $\forall\,(\alpha,\beta)\in\mathbb{C}^2$ we have
$|F^{-1}(\alpha,\beta)|\le [\mathbb{C}(X,Y):\mathbb{C}(F)]$, the extension degree of $F$ see \cite{e}. This is the so called Fiber Theorem for
\'etale mappings. Moreover the image is co-finite, $|\mathbb{C}^2-F(\mathbb{C}^2)|<\infty$, \cite{e}.  Also $F$ has a finite
set of exactly $d_F$ maximal domains $\{\Omega_1,\ldots,\Omega_{d_F}\}$. This means that $F$ is injective on each maximal domain
$\Omega_j$, and $i\ne j\Rightarrow \Omega_i\cap\Omega_j=\emptyset$, and $\mathbb{C}^2=\bigcup_{j=1}^{d_F}\overline{\Omega_j}$ and the
boundaries $\partial\Omega_j$ are piecewise smooth (even piecewise analytic). For the theory of maximal domains of entire functions in
one complex variable see \cite{rp1}, and for that theory for meromorphic functions in one complex variable see \cite{sh1,sh2}. Here
we use only basic facts of the theory which are valid also for more than complex variable. If $\{D\}$ is a family of characteristic sets of
${\rm et}(\mathbb{C}^2)$ such that $D\rightarrow\mathbb{C}^2$, then by the above $G_1(D),\,G_2(D)\rightarrow\mathbb{C}^2-A$, where $A$
is a finite set, and if $G_1\not\equiv G_2$ then we have the identity,
$$
F(G_1(D)\Delta G_2(D))-(F\circ G_1)(D)\Delta (F\circ G_2(D))=
$$
$$
=\{x=F(y)=F(z)|\,y\in G_1(D)-G_2(D)\,\wedge\,z\in G_2(D)-G_1(D)\}.
$$
Recalling that $(F\circ G_1)(D)\Delta (F\circ G_2(D))\subseteq F(G_1(D)\Delta G_2(D))$) we write the last identity as follows,
$$
F(G_1(D)\Delta G_2(D))=(F\circ G_1)(D)\Delta (F\circ G_2(D))\cup
$$
$$
\cup\{x=F(y)=F(z)|\,y\in G_1(D)-G_2(D)\,\wedge\,z\in G_2(D)-G_1(D)\}.
$$
Taking any two points $y\in G_1(D)-G_2(D)$ and $z\in G_2(D)-G_1(D)$ (as in the defining equation of the set on the right hand side in the last identity), 
we note that there are $i\ne j$, $1\le i,j\le d_F$ such that $y\in\Omega_i\,\wedge\,z\in\Omega_j$ (for $F(y)=F(z)$!).
For $\widetilde{D}$ a large enough characteristic set of ${\rm et}(\mathbb{C}^2)$, we will have $z\in G_1(\widetilde{D})$ and
$y\in G_2(\widetilde{D})$ and so $y,z\in G_1(\widetilde{D})\cap G_2(\widetilde{D})$ (since $G_1(D),G_2(D)\rightarrow \mathbb{C}^2-\{{\rm a\,finite\,set}\}$).
Hence $F(G_1(\widetilde{D})\Delta G_2(\widetilde{D}))-(F\circ G_1)(\widetilde{D})\Delta (F\circ G_2)(\widetilde{D})$ will not include the point $x$.
We conclude that if $y$ and $z$ are $F$-equivalent ($F(y)=F(z)$) then $x=F(y)=F(z)$ will not belong to $F(G_1(D)\Delta G_2(D))-(F\circ G_1)(D)
\Delta (F\circ G_2)(D)$ for large enough $D$. We obtain the following crude estimate:
$$
{\rm vol}(\{x=F(y)=F(z)|\,y\in G_1(D)-G_2(D)\,\wedge\,z\in G_2(D)-G_1(D)\})=
$$
$$
=o({\rm vol}((F\circ G_1)(D)\Delta (F\circ G_2(D)))).
$$
One can think of $D$ as a large open ball centered at the origin of $\mathbb{R}^4$, $D\approx B(R)$ and with the radius $R$ and look at the images of the two
polynomial \'etale mappings $F\circ G_1)(B(R))$ and $(F\circ G_2)(B(R))$ and compare the volume of $(F\circ G_1)(B(R))\Delta (F\circ G_2)(B(R))$ which is of the
order of magnitude $R^{4d}$, where $d$ depends on the algebraic degrees of $F\circ G_1$ and $F\circ G_2$, with the volume of the set
in the left hand side of the last equation. Similar estimates are used in the theory of covering surfaces by Ahlfors, see \cite{hay}, chapter 5.We conclude that,
$$
\lim_{D\rightarrow\mathbb{C}^2}\frac{{\rm vol}(F(G_1(D)\Delta G_2(D))}{{\rm vol}((F\circ G_1)(D)\Delta(F\circ G_2)(D)))}=1.
$$
Hence
$$
\lim_{D\rightarrow\mathbb{C}^2}\frac{\rho_D(F\circ G_1,F\circ G_2)}{\rho_D(G_1,G_2)}=\lim_{D\rightarrow\mathbb{C}^2}
\frac{{\rm vol}((F\circ G_1)(D)\Delta(F\circ G_2)(D))}{{\rm vol}(G_1(D)\Delta G_2(D))}=
$$
$$
=\lim_{D\rightarrow\mathbb{C}^2}\frac{{\rm vol}((F\circ G_1)(D)\Delta(F\circ G_2(D))}{{\rm vol}(F(G_1(D)\Delta G_2(D))}\cdot
\frac{{\rm vol}(F(G_1(D)\Delta G_2(D))}{{\rm vol}(G_1(D)\Delta G_2(D))}=
$$
$$
=1\cdot\frac{1}{d_F}=\frac{1}{d_F}.\,\,\,\,\,\qed
$$

\begin{remark}
The facts we used in proving (ii) and (iii) for \'etale mappings are in fact true in any dimension $n$, i.e. in $\mathbb{C}^n$.
In dimension $n=2$ the co-dimension of the image of the mapping is $0$. In fact the co-image is
a finite set. Also the fibers are finite and have a uniform bound on their cardinality (one can get a less tight uniform
bound by the Bezout Theorem). Here are few well known facts (which one can find in Hartshorne's book on Algebraic Geometry, \cite{hart}). \\
1) The following two conditions are equivalent: \\
a. The Jacobian Condition: the determinant $\det J_F$ is a non-zero constant. \\
b. The map $F^{*}$ is \'etale (in standard sense of algebraic geometry). In particular it is flat. \\
Let $F^{*}:\,Y\rightarrow X$ be \'etale. Let $X^{im}:=F^{*}(Y)\subseteq X$. \\
2) For every prime ideal $\wp\subseteq A$ ($X={\rm spec}(A)$), with residue field $k(\wp)$ the ring $B\otimes_A k(\wp)$is finite over
$k(\wp)$ ($Y={\rm spec}(B)$). \\
3) $F^{*}$ is a quasi-finite mapping. \\
4) The set $X^{im}$ is open in $X$. \\
5) For every point $x\in X(\mathbb{C})$ the fiber $(F^{*})^{-1}(x)$ is a finite subset of $Y(\mathbb{C})$. \\
6)The ring homomorphism $A\rightarrow B$ is injective, and the induced field extension $K\rightarrow L$ is finite. \\
7) There is a non-empty open subset $X^{fin}\subseteq X^{im}$ such that on letting $Y^{fin}:=(F^{*})^{-1}(X^{fin})\subseteq Y$,
the map of schemes $F^{*}|_{Y^{fin}}:\,Y^{fin}\rightarrow X^{fin}$ is finite. For any point $x\in X^{fin}(\mathbb{C})$ we have the equality
$d_x=d_{F^{*}}$ the geometrical degree of $F^{*}$. \\
8) The dimension of the set $Z:=X-X^{im}$ is at most $n-2$. \\
9) If $X^{im}=X-Z$ is affine, then $Z=\emptyset$ and $X^{im}=X$. \\
Let $X_{cl}$ be the topological space which is the set $X(\mathbb{C})\cong\mathbb{C}^n$ given the classical topology. Similarly for $Y_{cl}$.
The map of schemes $F^{*}:\,Y\rightarrow X$ induces a map of topological spaces $F_{cl}:\,Y_{cl}\rightarrow X_{cl}$ ($F_{cl}=f^{*}|_{Y(\mathbb{C})}$). \\
10) The map $F_{cl}:\,Y_{cl}\rightarrow X_{cl}$ is a local homeomorphism. \\
\end{remark}
\noindent
An immediate conclusion from Theorem 8.2 is the following,

\begin{corollary}
$\forall\,F\in {\rm et}(\mathbb{C}^2)$ the left (composition) mapping $L_F\,:\,{\rm et}(\mathbb{C}^2)\rightarrow L_F({\rm et}(\mathbb{C}^2))$, 
$L_F(G)=F\circ G$ is an injective mapping.
\end{corollary}
\noindent
\begin{remark}
We note the contrast in the behavior between polynomial mappings (in ${\rm et}(\mathbb{C}^2)$)
and entire functions of a single complex variable (in ${\rm elh}(\mathbb{C})$). See \cite{elh}.
\end{remark}
\noindent
We now state and prove the parallel of Proposition 4.13. Namely,

\begin{proposition}. \\
{\rm 1)} If $F\in {\rm et}(\mathbb{C}^2)$ and $G\in L_F({\rm et}(\mathbb{C}^2))$, then $L_G({\rm et}(\mathbb{C}^2))\subseteq L_F({\rm et}(\mathbb{C}^2))$. \\
{\rm 2)} If $F\in {\rm et}(\mathbb{C}^2)$, $G\in L_F({\rm et}(\mathbb{C}^2))$, and $G$ and $F$ are not associates (which means here $\forall\,H\in {\rm Aut}(\mathbb{C}^2)$,
$G\ne F\circ H$), then $L_G({\rm et}(\mathbb{C}^2))\subset L_F({\rm et}(\mathbb{C}^2))$. \\
{\rm 3)} $\forall\,F\in {\rm et}(\mathbb{C}^2)$ the spaces $(L_F({\rm et}(\mathbb{C}^2)),L^2)$ and $({\rm et}(\mathbb{C}^2),L^2)$ are homeomorphic. \\
{\rm 4)} There exists an infinite index set $I$ and a family of \'etale mappings $\{ F_i\,|\,i\in I\}\subseteq {\rm et}(\mathbb{C}^2)$ such that
$$
{\rm et}(\mathbb{C}^2)={\rm Aut}(\mathbb{C}^2)\cup\bigcup_{i\in I}L_{F_i}({\rm et}(\mathbb{C}^2)),
$$
so that $\forall\,i\in I$, ${\rm Aut}(\mathbb{C}^2)\cap L_{F_i}({\rm et}(\mathbb{C}^2))=\emptyset$, and $\forall\,i,j\in I$, if $i\ne j$ then
$F_i\not\in L_{F_j}({\rm et}(\mathbb{C}^2))$ and $F_j\not\in L_{F_i}({\rm et}(\mathbb{C}^2))$, and $L_{F_i}({\rm et}(\mathbb{C}^2))$ is
homeomorphic to $L_{F_j}({\rm et}(\mathbb{C}^2))$, and both are homeomorphic to ${\rm et}(\mathbb{C}^2)$.
\end{proposition}
\noindent
{\bf Proof.} \\
1) $H\in L_G({\rm et}(\mathbb{C}^2))\Rightarrow\exists\,M\in {\rm et}(\mathbb{C}^2)\,\,{\rm such}\,\,{\rm that}\,\,H=G\circ M$.
$G\in L_F({\rm et}(\mathbb{C}^2))\Rightarrow\exists\,N\in {\rm et}(\mathbb{C}^2)\,\,{\rm such}\,\,{\rm that}\,\,G=F\circ N$. Hence we conclude that
$H=G\circ M=(F\circ N)\circ M=F\circ(N\circ M)\in L_F({\rm et}(\mathbb{C}^2))$. \\
2) $G\in L_F({\rm et}(\mathbb{C}^2))$ and is not an associate of $F$ $\Rightarrow\exists\, N\in {\rm et}(\mathbb{C}^2)-{\rm Aut}(\mathbb{C}^2)$ such that
$G=F\circ N$. So $F\not\in L_G({\rm et}(\mathbb{C}^2))$ otherwise $F=(F\circ N)\circ M=F\circ(N\circ M)$. But $N\circ M\not\in {\rm Aut}(\mathbb{C}^2)$ (it
is not injective). Hence $d_{N\circ M}>1$. By the equation $F=F\circ(N\circ M)$ we get the contradiction $1<d_F=d_F\cdot d_{N\circ M}$. \\
3) The mapping $f=L_F\,:\,{\rm et}(\mathbb{C}^2)\rightarrow L_F({\rm et}(\mathbb{C}^2))$, $f(G)=F\circ G=L_F(G)$ is an homeomorphism (it is a
bijection and both $f$ and $f^{-1}$ are sequentially continuous). \\
4) We use the relation $\sim_R$ on ${\rm et}(\mathbb{C}^2)$ which was defined by $F\sim_R G\Leftrightarrow\exists\,\Phi\in{\rm Aut}(\mathbb{C}^2)$ such that $F=G\circ \Phi$.
Then $\sim_R $ is an equivalence relation ($F\sim_R F$ by $F=F\circ {\rm id}$, $F\sim_R G\Leftrightarrow F=G\circ \Phi\Leftrightarrow G=F\circ\Phi^{-1}
\Leftrightarrow G\sim_R F$, $F\sim_L G\,\wedge\,G\sim_R H\Leftrightarrow F=G\circ\Phi_1\,\wedge\,G=H\circ\Phi_2\Rightarrow F=H\circ (\Phi_2\circ\Phi_1)
\Rightarrow F\sim_R H$). We order the set of $\sim_R $ equivalence classes ${\rm et}(\mathbb{C}^2)/\sim_R$ by $[F]\preceq_R [G]\Leftrightarrow F\in
L_G({\rm et}(\mathbb{C}^2)) \Leftrightarrow L_F({\rm et}(\mathbb{C}^2))\subseteq L_G({\rm et}(\mathbb{C}^2))$. This relation is clearly reflexive and transitive by Proposition 8.5(1), and it is also anti-symmetric for 
$[F]\preceq_R [G]\wedge [G]\preceq_R [F]\Leftrightarrow L_F({\rm et}(\mathbb{C}^2))\subseteq L_G({\rm et}(\mathbb{C}^2))
\wedge L_G({\rm et}(\mathbb{C}^2))\subseteq L_F({\rm et}(\mathbb{C}^2))\Leftrightarrow L_F({\rm et}(\mathbb{C}^2))= L_G({\rm et}(\mathbb{C}^2))
\Rightarrow F=G\circ N\wedge G=F\circ M\Rightarrow F=F\circ(M\circ N)\Rightarrow d_{M\circ N}=1\Rightarrow M\circ N\in {\rm Aut}(\mathbb{C}^2)\Rightarrow
M,N\in {\rm Aut}(\mathbb{C}^2)\Rightarrow [F]=[G]$.
Any increasing chain in ${\rm et}(\mathbb{C}^2)-{\rm Aut}(\mathbb{C}^2)/\sim_R$ is finite (by an argument similar to that in the proof of Theorem 4.10). Hence every maximal increasing chain 
contains a maximal element $[F]$, and $L_F({\rm et}(\mathbb{C}^2))$ contains the union of the images of the left (composition) mappings of all the elements in this maximal chain
( by Proposition 8.5(1)). We define $I=\{ [F]\in ({\rm et}(\mathbb{C}^2)-{\rm Aut}(\mathbb{C}^2))/\sim_R\,|\,[F]\,\,{\rm is\,\,the\,\,maximum\,\,of\,\,a\,\,maximal\,\,
length\,\,chain}\}$. Then 
$$
{\rm et}(\mathbb{C}^2)={\rm Aut}(\mathbb{C}^2)\cup\bigcup_{i\in I}L_{F_i}({\rm et}(\mathbb{C}^2)).
$$
The union on the right equals ${\rm et}(\mathbb{C}^2)$ because any $F\in {\rm et}(\mathbb{C}^2)$ is either in ${\rm Aut}(\mathbb{C}^2)$ or $[F]$ belongs to some
maximal length chain in $({\rm et}(\mathbb{C}^2)-{\rm Aut}(\mathbb{C}^2))/\sim_R$ and so $L_F({\rm et}(\mathbb{C}^2))$ is a subset of $L_{F_i}({\rm et}(\mathbb{C}^2))$
where $[F_i]$ is the maximum of that chain. Clearly if $i\ne j$ then $F_i\not\in L_{F_j}({\rm et}(\mathbb{C}^2))\wedge F_j\not\in L_{F_i}({\rm et}(\mathbb{C}^2))$
for $[F_i]$ and $[F_j]$ the maxima of two different chains. So $L_{F_i}({\rm et}(\mathbb{C}^2))\ne L_{F_j}({\rm et}(\mathbb{C}^2))$. Finally, the index set $I$ is 
an infinite set. This follows because any finite union of the form:
$$
{\rm Aut}(\mathbb{C}^2)\cup L_{F_1}({\rm et}(\mathbb{C}^2))\cup\ldots\cup L_{F_k}({\rm et}(\mathbb{C}^2))
$$
is such that any mapping $H$ in it is either a $\mathbb{C}^2$-automorphism or $R_0(H)\supseteq R_0(F_1)\ne\emptyset\vee\ldots\vee R_0(H)\supseteq R_0(F_k)\ne\emptyset$.
Hence the argument in the proof of Proposition 3.2 implies that $L_{F_1}({\rm et}(\mathbb{C}^2))\cup\ldots\cup L_{F_k}({\rm et}(\mathbb{C}^2))\subset
{\rm et}(\mathbb{C}^2)-{\rm Aut}(\mathbb{C}^2)$. $\qed $ \\
\\
We recall definition 4.3: an \'etale mapping $F\in  {\rm et}(\mathbb{C}^2)$ is composite if $\exists\,G,H\in {\rm et}(\mathbb{C}^2)-{\rm Aut}(\mathbb{C}^2)$ such that
$F=G\circ H$. An \'etale mapping is prime if it is not composite. If ${\rm et}(\mathbb{C}^2)\ne {\rm Aut}(\mathbb{C}^2)$ then we know that the set of all the prime
\'etale mappings is not empty. Also we know that the geometrical degree of a composite mapping is not a prime number. In other words an \'etale mapping whose
geometrical degree is a prime number is a prime \'etale mapping. Also we know that any \'etale mapping $F\in {\rm et}(\mathbb{C}^2)$ can be written as follows:
$F=A_0\circ P_1\circ\ldots\circ P_k$, for some $A_0\in {\rm Aut}(\mathbb{C}^2)$ and prime \'etale mappings $P_1,\ldots,P_k$. By Proposition 8.5 it follows that if
$F\in {\rm et}(\mathbb{C}^2)$ is composite, say $F=G\circ H$ for some $G,H\in {\rm et}(\mathbb{C}^2)-{\rm Aut}(\mathbb{C}^2)$, then
$L_F({\rm et}(\mathbb{C}^2))\subset L_G({\rm et}(\mathbb{C}^2))$. We conclude that in the fractal representation of ${\rm et}(\mathbb{C}^2)$:
$$
{\rm et}(\mathbb{C}^2)={\rm Aut}(\mathbb{C}^2)\cup\bigcup_{i\in I}L_{F_i}({\rm et}(\mathbb{C}^2)).
$$
all the maximal elements $F_i$ must be prime \'etale mappings. Conversely, it is clear that if $F$ is \'etale prime, then $F$ is the maximum of some
(finite) $\preceq_R$-chain. Hence we can state a more accurate statement than that of Proposition 8.5(4).

\begin{proposition}
Suppose that ${\rm et}(\mathbb{C}^2)\ne {\rm Aut}(\mathbb{C}^2)$. Let ${\rm et_p}(\mathbb{C}^2)$ be the set of all the prime \'etale mappings. Then
$|{\rm et_p}(\mathbb{C}^2)|=\infty$ and we have,
$$
{\rm et}(\mathbb{C}^2)={\rm Aut}(\mathbb{C}^2)\cup\bigcup_{F\in {\rm et_p}(\mathbb{C}^2)}L_F({\rm et}(\mathbb{C}^2)),
$$
where $\forall\,F\in {\rm et_p}(\mathbb{C}^2)$, ${\rm Aut}(\mathbb{C}^2)\cap L_F({\rm et}(\mathbb{C}^2))=\emptyset$, and $\forall\,F,G\in {\rm et_p}(\mathbb{C}^2)$,
$F\ne G$, we have $F\not\in L_G({\rm et}(\mathbb{C}^2))\wedge G\not\in L_F({\rm et}(\mathbb{C}^2))$ and $L_F({\rm et}(\mathbb{C}^2))$ is homeomorphic to
$L_G({\rm et}(\mathbb{C}^2))$
\end{proposition}
\noindent
We recall that $\forall\,F\in {\rm et_p}(\mathbb{C}^2)$, the corresponding left space $L_F({\rm et}(\mathbb{C}^2))$ is composed of all
the \'etale mappings $G$ that have the form $G=F\circ H$ for some $H\in {\rm et}(\mathbb{C}^2)$. Hence the integer $d_F$ divides the geometrical
degree, $d_G$, of $G$. We know that the set of prime \'etale mappings is infinite (if non-empty). However, concerning their geometrical
degrees $\{ d_F\,|\,F\in {\rm et_p}(\mathbb{C}^2)\}$, we do not know much. If for an $F\in {\rm et}(\mathbb{C}^2)$ we have know that $d_F$ is a prime
integer then $F$ is a prime mapping. But the set of geometric degrees of prime mappings might contain other integers (composite). It might be $\mathbb{Z}^+$.
We now show how to get some non trivial information regarding that.

\begin{theorem}
$|{\rm et}_p(\mathbb{C}^2)|\le \aleph_0$ and ${\rm et}_p(\mathbb{C}^2)$ is a discrete subset of the metric space $({\rm et}(\mathbb{C}^2),\rho_D)$.
\end{theorem}
\noindent
{\bf Proof.} \\
We later on (in Theorem 14.14) will prove that if $F,G\in {\rm et}_p(\mathbb{C}^2)$ satisfy $F\ne G$ then $L_F({\rm et}(\mathbb{C}^2))\cap L_G({\rm et}(\mathbb{C}^2))=\emptyset$.
When $D\rightarrow\mathbb{C}^2$, ${\rm meas}(L_F({\rm et}(\mathbb{C}^2)))\rightarrow {\rm meas}({\rm et}(\mathbb{C}^2))/d_F$. So by the identity
${\rm meas}(\bigcup_{F\in {\rm et}_p(\mathbb{C}^2)} L_F({\rm et}(\mathbb{C}^2)))={\rm meas}({\rm et}(\mathbb{C}^2))$ (the measure is 
$H^{s_0}$ where $s_0=\dim_H {\rm et}(\mathbb{C}^2)$) we get:
$$
\sum_{F\in {\rm et}(\mathbb{C}^2)}\frac{1}{d_F}=1.\,\,\,\,\,\,\,\,\,\,\qed
$$

\begin{corollary}
The sequence of geometric degrees $\{ d_F\,|\,F\in {\rm et}_p(\mathbb{C}^2)\}$ can not equal to $\mathbb{Z}^+_{\ge 2}$, and
can not equal to the set of prime integers.
\end{corollary}

\begin{corollary}
$({\rm et}(\mathbb{C}^2),\circ)$ is generated by ${\rm Aut}(\mathbb{C}^2)$ and the countable set of generators given by ${\rm et}_p(\mathbb{C}^2)$.
\end{corollary}

\section{Properties of the metric spaces $({\rm et}(\mathbb{C}^2),\rho_D)$}

Here is a natural list of questions about those metric spaces: \\
1) Is it a separable space? \\
2) Is it a proper space? \\
3) Is it a complete space? \\
4) Is the action of the group ${\rm Aut}(\mathbb{C}^2)$ from the left, as isometries on the space, cocompact? \\
5) Is the action described in 4 above, a proper action? \\
\begin{remark}
The action of ${\rm Aut}(\mathbb{C}^2)$ from the left on $({\rm et}(\mathbb{C}^2),\rho_D)$ is certainly not both cocompact and proper.
This follows by Lemma 1.17 on page 8 of \cite{roe}. For if the action were both cocompact and proper , then by (b) of the Lemma
this would have implied that ${\rm Aut}(\mathbb{C}^2)$ is finitely generated.
\end{remark}
\noindent
6) Is the metric space $({\rm et}(\mathbb{C}^2),\rho_D)$ a length space? \\
7) What are the geodesics if any? \\
\\
To tackle question 1 we might ask the following. \\
8) Are the \'etale mappings with rational coefficients dense in the space $({\rm et}(\mathbb{C}^2),\rho_D)$?
Are the mappings in ${\rm Aut}(\mathbb{C}^2)$ with rational coefficients dense in ${\rm Aut}(\mathbb{C}^2)$?
Here rational coefficient are numbers in $\mathbb{Q}+i\mathbb{Q}$. 
\\
1) \underline{Separability of $({\rm et}(\mathbb{C}^2),\rho_D)$ and of $({\rm Aut}(\mathbb{C}^2),\rho_D)$} \\
\\
The following is probably well known.
\begin{proposition}
The metric space $({\rm Aut}(\mathbb{C}^2,\rho_D)$ is separable.
\end{proposition}
\noindent
{\bf Proof.} \\
It is well known that the group $({\rm Aut}(\mathbb{C}^2),\circ)$ is generated by the affine mappings, 
$F(X,Y)=(aX+bY+c,dX+eY+f)$, $ae-bd=1$, and by the elementary mappings $G(X,Y)=(X+P(Y),Y)$ or $H(X,Y)=(X,Y+P(X))$
where $P(T)\in\mathbb{C}[T]$ (The Jung-Van Der Kulk Theorem, \cite{e}). Using this, we will show that any automorphism
of $\mathbb{C}^2$ can be approximated well enough by automorphisms of $\mathbb{C}^2$ which have all of their coefficients
from $\mathbb{Q}+i\mathbb{Q}$. For we have: \\
(i) $\forall\,(aX+bY+c, dX+eY+f),\,\,ae-bd=1,\,\,\exists\,(a_kX+b_kY+c_k,d_kX+e_kY+f_k),\,k\in\mathbb{Z}^+$ such that
$a_ke_k-b_kd_k=1,\,a_k,b_k,c_k,d_k,e_k,f_k\in\mathbb{Q}+i\mathbb{Q}$ and such that $a=\lim a_k\wedge\,b=\lim b_k\wedge\,
c=\lim c_k\wedge\,d=\lim d_k\wedge\,e=\lim e_k\wedge\,f=\lim f_k$. \\
(ii) $\forall\,(X+P(Y),Y),\,\,P(Y)\in\mathbb{C}[Y]\,\,\exists\,(X+P_k(Y),Y),\,P_k(Y)\in(\mathbb{Q}+i\mathbb{Q})[Y],\,k\in\mathbb{Z}^+$
and also $\deg P_k=\deg P$, $\lim P_k=P$ coefficientwise. \\
(iii) The same as in (ii) $\forall\,(X,Y+P(X))$. \\
We conclude that for any $F\in {\rm Aut}(\mathbb{C}^2)$ we have a sequence $F_k\in {\rm Aut}(\mathbb{C}^2)$, $k\in\mathbb{Z}^+$, which have 
all of the coefficients from $\mathbb{Q}+i\mathbb{Q}$ and which satisfy $\lim F_k=F$ coefficientwise and where the degrees
$\deg F_k\le M$ are uniformly bounded in $k$. Hence $\forall\,K\subset\mathbb{C}^2$ a compact in the strong topology
we have $\lim F_k=F$ uniformly on $K$. This follows by the proposition that follows this one. This implies that
$\lim F_k=F$ in $({\rm Aut}(\mathbb{C}^2),\rho_D)$ and concludes the proof. $\qed $

\begin{remark}
The uniform bound $M$ on the degrees $\deg F_k$ could be taken to be $M=\prod\deg E_k$ where $F=\circ_k E_k$ is a decomposition of $F$ into affine
and elementary mappings.
\end{remark}

\begin{proposition}
Let $F(X,Y)\in\mathbb{C}[X,Y]$, $F_k\in\mathbb{C}[X,Y]$, $k\in\mathbb{Z}^+$ and $K\subset\mathbb{C}^2$ a compact subset in the strong topology.
If $\lim F_k=F$ coefficientwise and $\deg F_k\le M$ are uniformly bounded in $k$, then $\lim F_k=F$ uniformly on $K$.
\end{proposition}

\begin{remark}
If we drop the uniform degree bound condition in Proposition 9.4, its conclusion is false. For example we take $F(X,Y)\equiv 0$, and
$F_k(X,Y)=(1/k)(X^k+X^{k+1}+\ldots+X^{2k-1})$ and $K=[0,1]\times\{ 0\}$. Then $\lim F_k=F$ coefficientwise but not even pointwise
because $F(1,0)=0$ while $\forall\,k\in\mathbb{Z}^+$, $F_k(1,0)=1$.
\end{remark}
\noindent
{\bf A proof of Proposition 9.4}. \\
By the degree bound assumption we can write,
$$
F(X,Y)=\sum_{i+j\le M} a_{ij}X^iY^j,\,\,\,\,F_k(X,Y)=\sum_{i+j\le M} a_{ij}^{(k)} X^iY^j,\,\,\,k\in\mathbb{Z}^+.
$$
We have $\forall\,i+j\le M$, $\lim_{k\rightarrow\infty} a_{ij}^{(k)}=a_{ij}$. We note that,
$$
|F(X,Y)-F_k(X,Y)|=|\sum_{i+j\le M}(a_{ij}-a_{ij}^{(k)})X^iY^j|\le\sum_{i+j\le M}|a_{ij}-a_{ij}^{(k)}||X^iY^j|.
$$
For any $i+j\le M$ the function $|X^iY^j|$ is continuous on the compact $K$ and so $m_{ij}=\max_{(X,Y)\in K}|X^iY^j|<\infty$.
Thus $m=\max\{m_{ij}\,|\,i+j\le M\}<\infty$. Given $\epsilon>0$ $\exists\,N$ such that $\forall\,k>N$ we have $\sum_{i+j\le M}|a_{ij}-a_{ij}^{(k)}|\le \epsilon/(m+1)$.
We conclude that $\max_{(X,Y)\in K}|F(X,Y)-F_k(X,Y)|\le (\epsilon\cdot m)/(m+1)<\epsilon$, $\forall\,k>N$. $\qed $ \\

\begin{remark}
We can generalize Proposition 9.4 as follows: \\
Let $F(X,Y)=\sum_{i+j\le M} a_{ij}X^iY^j$, $F_k(X,Y)=\sum_{i+j\le M_k} a_{ij}^{(k)}X^iY^j$, $k\in\mathbb{Z}^+$ and $K\subset\mathbb{C}^2$
a compact subset of $\mathbb{C}^2$ in the strong topology. If
$$
\sum_{i+j\le M}a_{ij}^{(k)}X^iY^j\rightarrow_{k\rightarrow\infty} F\,\,{\rm uniformly\,on}\,\,K
$$
and if
$$
\sum_{M<i+j\le M_k}a_{ij}^{(k)}X^iY^j\rightarrow_{k\rightarrow\infty} 0\,\,{\rm uniformly\,on}\,\,K
$$
then, $\lim F_k=F$ uniformly on $K$.
\end{remark}

\begin{proposition}
The metric space $({\rm et}(\mathbb{C}^2),\rho_D)$ is separable.
\end{proposition}
\noindent
{\bf Proof.} \\
We will find a countable subset ${\rm et}_{\aleph_0}(\mathbb{C}^2)\subset {\rm et}(\mathbb{C}^2)$ such that
$\forall\,F\in {\rm et}(\mathbb{C}^2)$, $\exists\,F_k\in {\rm et}_{\aleph_0}(\mathbb{C}^2)$, $k\in\mathbb{Z}^+$ such
that $\lim F_k=F$ coefficientwise and we have a uniform degree bound $\deg F_k\le\deg F$. Thus we will be able to apply once more Proposition 9.4 
and conclude the result. Unlike the situation with $({\rm Aut}(\mathbb{C}^2),\rho_D)$, this time we can not say
something definite about the character of the coefficients of the $F_k$'s (e.g. that all belong to $\mathbb{Q}+i\mathbb{Q}$).
The elements of the set ${\rm et}(\mathbb{C}^2)$ are faithfully parametrized by the sets of the solutions of certain polynomial systems
of equations that are induced on the coefficients of the \'etale mappings by the Jacobian Condition, $\det J_F(X,Y)\equiv 1$.
These systems have all of their equations, quadratic and homogeneous, except for a single equation which is still quadratic
but not homogeneous. $\forall\,k\in\mathbb{Z}^+$ we let ${\rm et_k}(\mathbb{C}^2)=\{ F\in {\rm et}(\mathbb{C}^2)\,|\,\deg F\le k\}$.
This subset of ${\rm et}(\mathbb{C}^2)$ is parametrized by certain finite dimensional complex space, which is composed of
a part of the above mentioned solutions. We denote this complex space by $J_k(2)$ or simply by $J_k$. Its finite dimension is
a function of $k$. We have $\forall\,k$, $J_k\subset J_{k+1}$ (proper containment) and $J_k$ is a path connected topological subspace
of the appropriate $\mathbb{C}^{D_k}$ in the strong topology. Given $k\in\mathbb{Z}^+$ and $\epsilon>0$ we can construct on $J_k$
a countable $\epsilon$-net, $N_k(\epsilon)$. This means that $N_k(\epsilon)=\{ P_j^{(\epsilon)}\,|\,j\in\mathbb{Z}^+\}$ and that
$\forall\,P\in J_k$, $\exists\,j_0\in\mathbb{Z}^+$ such that $P_{j_0}^{(\epsilon)}$ is close to $P$ in the strong metric by
less than $\epsilon $. Thus $||P-P_{j_0}^{(\epsilon)}||_2<\epsilon$. If the point $P$ is the parameter value of the \'etale mapping 
$F\in {\rm et}_k(\mathbb{C}^2)$ and if $P_{j_0}^{(\epsilon)}$ is the parameter value of the \'etale mapping 
$F_{j_0}^{(\epsilon)}\in {\rm et}_k(\mathbb{C}^2)$, then this implies that the $l_2$-distance between the coefficients of
$F$ and those of $F_{j_0}^{(\epsilon)}$ is less than $\epsilon$. In particular this means that the following is true: \\
if ${\rm et}_{\aleph_0}(\mathbb{C}^2)$ is the subset of \'etale mappings in ${\rm et}(\mathbb{C}^2)$ that are parametrized
by
$$
\bigcup_{(k,n)\in\mathbb{Z}^+} N_k\left(\frac{1}{n}\right),
$$
then $\forall\,F\in {\rm et}(\mathbb{C}^2)$, $\exists\,F_k\in {\rm et}_{\deg F}(\mathbb{C}^2)\cap {\rm et}_{\aleph_0}(\mathbb{C}^2)$ so
that $\lim F_k=F$ coefficientwise. 

Since we have a uniform bound on the degrees, $\deg F_k\le\deg F$ we can indeed apply Proposition 9.4. $\qed $ \\

\begin{remark}
We note the important role of Proposition 9.4 in the proofs of Proposition 9.2 and of Proposition 9.7. In particular, the uniform degree
bound was central. We naturally ask if this uniform degree bound is a natural necessary condition for our needs ($\rho_D$-convergence). More
concretely we ask if $\forall\,M>0$, $\exists\,\epsilon=\epsilon(M)>0$ such that $\forall\,F,G\in {\rm et}(\mathbb{C}^2)$ if
$|\deg F-\deg G|>M$ then ${\rm vol}(F(D)\Delta G(D))>\epsilon$? The answer is clearly negative ($G(X,Y)=(X,Y)$ the identity mapping, and
$F(X,Y)=(X+(1/n)Y^M,Y)$). In fact we have no good idea on the character of $\rho_D$-convergence. For example, if $F_k\in {\rm et}(\mathbb{C}^2)$
and $F\in {\rm et}(\mathbb{C}^2)$ satisfy $\lim_{k\rightarrow\infty}\rho_D(F_k,F)=0$ is it true that $\lim F_k=F$ coefficientwise? Pointwise?
Uniformly on $D$?
\end{remark}
\noindent
We end this section with the following elementary observation,

\begin{proposition}
The metric space $({\rm et}(\mathbb{C}^2),\rho_D)$ is bounded.
\end{proposition}
\noindent
{\bf Proof.} \\
The diameter of ${\rm et}(\mathbb{C}^2)$ with respect to $\rho_D$ is bounded from above by twice the volume of the characteristic compact $D$. $\qed $ \\

\section{Connections between the two dimensional Jacobian Conjecture and number theory}

In this short section we just indicate the general principles that underline the nontrivial connections that
are indicated in the title. We also fix some notations. In the sections that follow we make more concrete computations
that actually undercover some of these connections.

We consider the fractal representation of ${\rm et}(\mathbb{C}^2)$ using the left (composition) mappings. Namely (Proposition 8.6):
$$
{\rm et}(\mathbb{C}^2)={\rm Aut}(\mathbb{C}^2)\cup\bigcup_{F\in {\rm et}_p(\mathbb{C}^2)} L_F({\rm et}(\mathbb{C}^2)).
$$
We will denote the $s$-dimensional Hausdorff measure of the set $A$, by $H^{s}(A)$. Then we have the following,
$$
H^s({\rm et}(\mathbb{C}^2))=H^s({\rm Aut}(\mathbb{C}^2))+H^s(\bigcup_{F\in {\rm et}_p(\mathbb{C}^2)}L_F({\rm et}(\mathbb{C}^2)).
$$
This follows by the disjointness ${\rm Aut}(\mathbb{C}^2)\cap\bigcup_{F\in {\rm et}_p(\mathbb{C}^2)}L_F({\rm et}(\mathbb{C}^2))=\emptyset$.
We use the parametrization of the space of the \'etale prime mappings ${\rm et}_p(\mathbb{C}^2)$ and measure using that parameter. We get,
$$
H^s({\rm et}(\mathbb{C}^2))\le H^s({\rm Aut}(\mathbb{C}^2))+\int_{\{F\in {\rm et}_p(\mathbb{C}^2)\}}f(F,s)\cdot H^s({\rm et}(\mathbb{C}^2))d\mu(F).
$$
Here the notation $f(F,s)$ stands for the self similarity factor between ${\rm et}(\mathbb{C}^2)$ and the scaling down $L_F({\rm et}(\mathbb{C}^2))$, i.e.
$H^s(L_F({\rm et}(\mathbb{C}^2)))=f(F,s)\cdot H^s({\rm et}(\mathbb{C}^2))$. We conclude the following inequality,
$$
H^s({\rm et}(\mathbb{C}^2))\left(1-\int_{\{F\in {\rm et}_p(\mathbb{C}^2)\}}f(F,s)d\mu(F)\right)\le H^s({\rm Aut}(\mathbb{C}^2)).
$$
If $H^s({\rm et}(\mathbb{C}^2))\ne 0$, then this can also be written as follows:
\begin{equation}
1-\int_{\{F\in {\rm et}_p(\mathbb{C}^2)\}}f(F,s)d\mu(F)\le\frac{H^s({\rm Aut}(\mathbb{C}^2)}{H^s({\rm et}(\mathbb{C}^2))}.
\label{eq13.2}
\end{equation}

\begin{remark}
We will denote by $s_0=\dim_H {\rm et}(\mathbb{C}^2)$ the Hausdorff dimension of the \'etale mappings ${\rm et}(\mathbb{C}^2)$.
It is plausible that $s_0=\dim_H {\rm et}(\mathbb{C}^2)>\dim_H {\rm Aut}(\mathbb{C}^2)$ in which case we have
$H^{s_0}({\rm et}(\mathbb{C}^2))>0$ and $H^{s_0}({\rm Aut}(\mathbb{C}^2))=0$.
\end{remark}
\noindent
By equation (\ref{eq13.2}) we get,

\begin{proposition}
$1\le\int_{\{F\in {\rm et}_p(\mathbb{C}^2)\}}f(F,s_0)d\mu(F)$.
\end{proposition}

\section{Inequalities and identities}

\begin{proposition}
If $D$ denotes the characteristic set that defines the metric $\rho_D$ on $({\rm et}(\mathbb{C}^2),\rho_D)$, then $\lim_{D\rightarrow\mathbb{C}^2} f(F,s)=d_F^{-s}$.
\end{proposition}
\noindent
{\bf Proof.} \\
This follows by Proposition 8.2. $\qed $ \\

\begin{proposition}
$1\le\int_{\{F\in {\rm et}_p(\mathbb{C}^2)\}}d_F^{-s_0}dF$.
\end{proposition}
\noindent
{\bf Proof.} \\
Using Proposition 10.2 and Proposition 11.1 and dominated convergence. $\qed $ \\
\\
So either $H^{s_0}({\rm et}(\mathbb{C}^2))=0$ in which case equation (\ref{eq13.2}) is invalid (with $s=s_0$) or $s_0=\dim_H {\rm et}(\mathbb{C}^2)<\infty$
otherwise Proposition 11.2 is invalid.
\begin{definition}
$\forall\,n\in\mathbb{Z}^+$, $\mu_p(n)=\int_{\{F\in {\rm et}_p(\mathbb{C}^2)\,|\,d_F=n\}}d\mu(F)$.
\end{definition}
\noindent
This allows us to rewrite the last proposition as follows,

\begin{proposition}
$$
1\le \sum_{n=2}^{\infty}\frac{\mu_p(n)}{n^{s_0}}.
$$
\end{proposition}
\noindent
If we had the disjointness property (see the explanation that follows equation (1.2)) in our fractal representation, 
the inequality signs would have become equalities and we could have deduced the following,

\begin{proposition}
If $\forall\,F,G\in {\rm et}_p(\mathbb{C}^2)$, the assumption that $F\ne G$ implied $L_F({\rm et}(\mathbb{C}^2))\cap L_G({\rm et}(\mathbb{C}^2))=\emptyset$,
then:
$$
1=\sum_{n=2}^{\infty}\frac{\mu_p(n)}{n^{s_0}},
$$
and hence $\forall\,n\in\mathbb{Z}^+$, $\mu_p(n)/n^{s_0}\le 1$, i.e.
$$
{\rm meas}\{F\in {\rm et}_p(\mathbb{C}^2)\,|\,d_F=n\}:=\int_{\{F\in {\rm et}_p(\mathbb{C}^2)\,|\,d_F=n\}}d\mu(F)\le n^{s_0}.
$$
Moreover, asymptotically:
$$
\lim_{n\rightarrow\infty}\frac{\mu_p(n)}{n^{s_0-1}}=0,
$$
i.e.
$$
{\rm meas}\{F\in {\rm et}_p(\mathbb{C}^2)\,|\,d_F=n\}=o(n^{s_0-1}).
$$
\end{proposition}
\noindent
{\bf Proof.} \\
Only the last parts needs a proof. The series with nonnegative terms
$$
\sum_{n=2}^{\infty}\frac{\mu_p(n)}{n^{s_0}},
$$
converges to $1$ and so a comparison with the divergent harmonic series $\sum (1/n)$ gives us the desired estimate,
$$
\lim_{n\rightarrow\infty}\frac{\mu_p(n)}{n^{s_0-1}}=0.\,\,\,\,\,\,\,\,\,\qed
$$

\begin{remark}
How could we prove the disjointness in the fractal representation? Let $F,G\in {\rm et}_p(\mathbb{C}^2)$ and $F\ne G$. Suppose that we did
not have disjointness, say $H\in L_F({\rm et}(\mathbb{C}^2))\cap L_G({\rm et}(\mathbb{C}^2))$. Then $\exists\,M,N\in {\rm et}(\mathbb{C}^2)$
such that $H=F\circ N=G\circ M$. If we knew that there is a unique factorization of $H$ where uniqueness includes the order of the factors, then this
would have been it.
\end{remark}
\noindent
Next we make the following, 

\begin{definition}
$\forall\,n\in\mathbb{Z}^+$, $\Omega(n)=\{F\in {\rm et}(\mathbb{C}^2)\,|\,d_F=n\}$.
\end{definition}
\noindent
Then we clearly have,

\begin{proposition}
$\Omega(n)=\bigcup_{k|n,\,k\ge 1}\Omega(k)\circ\Omega(n/k)$.
\end{proposition}
\noindent
However, in this representation many mappings appear in many components of the union on the right hand side. In other words, this representation
is very far from being a partition of $\Omega(n)$. In order to get a better representation, we refine our definitions,

\begin{definition}
$\forall\,n\in\mathbb{Z}^+$, $\Omega_p(n)=\Omega(n)\cap {\rm et}_p(\mathbb{C}^2)=\{F\in {\rm et}_p(\mathbb{C}^2)\,|\,d_F=n\}$. We will
use the short notation:
$$
\Omega_p(n_1)\circ\Omega_p(n_2)\circ\ldots\circ\Omega_p(n_k)=\bigcirc_{j=1}^k \Omega_p(n_j).
$$
It should be noted that in this "product" order matters because composition is not commutative.
\end{definition}
\noindent
We have the following improvement of Proposition 11.8,

\begin{proposition}
$\Omega(n)=\bigcup_{n_1 n_2\ldots n_k=n}\{\bigcirc_{j=1}^k\Omega_p(n_j)\}$. The order of the factors $n_1, n_2,\ldots, n_k$ in
the product $n_1 n_2\ldots n_k=n$ is important.
\end{proposition}
\noindent
{\bf Computing few examples.} \\
If $p,\,q$ are two prime integers then if $p\ne q$ we have 
$$
\mu(pq)\le \mu_p(pq)+2\mu_p(p)\mu_p(q).
$$ 
Here we are using the following,

\begin{definition}
$$
\mu(n)=\mu(\Omega(n))={\rm meas}(\Omega(n))=\int_{\{F\in {\rm et}(\mathbb{C}^2)\,|\,d_F=n\}}d\mu(F),
$$
$$
\mu_p(n)=\mu(\Omega_p(n))={\rm meas}(\Omega_p(n))=\int_{\{F\in {\rm et}_p(\mathbb{C}^2)\,|\,d_F=n\}}d\mu(F).
$$
\end{definition}
\noindent
The above inequality is an immediate consequence of the following identity which follows by the definitions,
$$
\Omega(pq)=\Omega_p(pq)\cup(\Omega_p(p)\circ\Omega_p(q))\cup(\Omega_p(q)\circ\Omega_p(p)),\,\,\,p\ne q.
$$
The inequality originates in the possibility that
$$
(\Omega_p(p)\circ\Omega_p(q))\cap(\Omega_p(q)\circ\Omega_p(p))\ne\emptyset,\,\,\,p\ne q.
$$
Once more, if we knew that the disjointness property holds true, i.e.,
$$
(\Omega_p(p)\circ\Omega_p(q))\cap(\Omega_p(q)\circ\Omega_p(p))=\emptyset,\,\,\,p\ne q,
$$
then we could have deduced the sharper result:
$$
\mu(pq)=\mu_p(pq)+2\mu_p(p)\mu_p(q),\,\,\,p\ne q.
$$
By Proposition 11.5 we have $\mu(pq)/(pq)^{s_0}\le 1$ and hence
$$
\frac{\mu_p(pq)}{(pq)^{s_0}}+2\left(\frac{\mu_p(p)}{p^{s_0}}\right)\left(\frac{\mu_p(q)}{q^{s_0}}\right)\le 1,\,\,\,p\ne q.
$$

\begin{remark}
For $n=q$ an integral prime we clearly have the identity $\Omega(q)=\Omega_p(q)$ and hence $\mu(q)=\mu_p(q)$.
\end{remark}
\noindent
We recall that if $q_1,\,q_2$ are two integral primes, then (as we saw) we have:
$$
\Omega(q_1q_2)=\Omega_p(q_1q_2)\cup(\Omega_p(q_1)\circ\Omega_p(q_2))\cup(\Omega_p(q_2)\circ\Omega_p(q_1)).
$$
We certainly have the following:
$$
\Omega_p(q_1q_2)\cap(\Omega_p(q_1)\circ\Omega_p(q_2))=\Omega_p(q_1q_2)\cap(\Omega_p(q_2)\circ\Omega_p(q_1))=\emptyset.
$$
However, we might have $(\Omega_p(q_1)\circ\Omega_p(q_2))\cap(\Omega_p(q_2)\circ\Omega_p(q_1))\ne\emptyset$. The last
possibility happens when there are two pairs of mappings $F_j\in\Omega_p(q_1)$, $G_j\in\Omega_p(q_2)$, $j=1,2$ that satisfy:
$H=F_1\circ G_1=G_2\circ F_2$. In this case we have the following identities for the Jacobian varieties: 
$A(H)=A(F_1)\cup F_1(A(G_1))=A(G_2)\cup G_2(A(F_2))$. Also, for the geometric bases we have: $R_0(G_1),\,R_0(F_2)\subseteq R_0(H)$.
Finally, the chain rule implies for the Jacobian matrices the following: $J_{F_1}(G_1)\cdot J_{G_1}(X,Y)=J_{G_2}(F_2)\cdot J_{F_2}(X,Y)$.
The special case $q=q_1=q_2$ is somehow easier because there is no overlap in this case. So in that case we have the 
following partition:
$$
\Omega(q^2)=\Omega_p(q^2)\cup(\Omega_p(q)\circ\Omega_p(q))=\Omega_p(q^2)\cup(\Omega(q)\circ\Omega(q)).
$$
Hence $\mu(q^2)=\mu_p(q^2)+\mu(q)^2$ or, equivalently $\mu(q^2)-\mu(q)^2=\mu_p(q^2)\ge 0$. Hence $\mu(q)\le\sqrt{\mu(q^2)}$. Thus
$$
\frac{\mu(q)}{q^{s_0}}\le\sqrt{\frac{\mu(q^2)}{(q^2)^{s_0}}}.
$$
\begin{proposition}
For any prime integer $q$ we have the following inequality:
$$
0\le\mu(q)\le\left(\frac{\sqrt{5}-1}{2}\right)q^{s_0}.
$$
\end{proposition}
\noindent
{\bf Proof.} \\
We recall that we have the following inequality:
$$
\frac{\mu(q)}{q^{s_0}}+\frac{\mu(q^2)}{(q^2)^{s_0}}\le 1.
$$
Thus if we use the inequality above:
$$
\frac{\mu(q)}{q^{s_0}}\le\sqrt{\frac{\mu(q^2)}{(q^2)^{s_0}}}.
$$
we obtain the following estimate:
$$
\frac{\mu(q)}{q^{s_0}}+\left(\frac{\mu(q)}{q^{s_0}}\right)^2\le\frac{\mu(q)}{q^{s_0}}+\frac{\mu(q^2)}{(q^2)^{s_0}}\le 1.
$$
Hence:
$$
0\le\frac{\mu(q)}{q^{s_0}}\le\frac{-1+\sqrt{1+4}}{2}.\,\,\,\,\,\,\,\,\,\,\,\qed
$$
Here is another similar series of arguments. We start with:
$$
\Omega(q^3)=\Omega_p(q^3)\cup(\Omega_p(q^2)\circ\Omega_p(q))\cup(\Omega_p(q)\circ\Omega_p(q^2))\cup(\Omega_p(q)\circ\Omega_p(q)\circ\Omega_p(q)),
$$
so
$$
\Omega(q^3)=\Omega_p(q^3)\cup(\Omega_p(q^2)\circ\Omega(q))\cup(\Omega(q)\circ\Omega_p(q^2))\cup(\Omega(q)\circ\Omega(q)\circ\Omega(q)),
$$
$$
\mu(q^3)\le\mu_p(q^3)+2\mu_p(q^2)\mu(q)+\mu(q)^3.
$$
But $\mu(q^2)=\mu(q)^2+\mu_p(q^2)$ and hence:
$$
\mu(q^3)\le\mu_p(q^3)+2\mu(q^2)\mu(q)-\mu(q)^3,\,\,\,{\rm hence}\,\,\,\mu(q^3)+\mu(q)^3\le 2\mu(q^2)\mu(q).
$$
We end here this chain of computational examples and proceed further in the theory of the fractal representation of ${\rm et}(\mathbb{C}^2)$.
Our starting point is a further generalization of Proposition 11.2. We recall that if $s_0=\dim_H {\rm et}(\mathbb{C}^2)$ the Hausdorff dimension
of ${\rm et}(\mathbb{C}^2)$, then 
$$
1\le\int_{\{F\in{\rm et}_p(\mathbb{C}^2)\}}f(F,s_0)d\mu(F).
$$
Let us denote the Hausdorff dimension of the automorphism group of $\mathbb{C}^2$ by $s_1=\dim_H {\rm Aut}(\mathbb{C}^2)$. Since
${\rm Aut}(\mathbb{C}^2)\subseteq{\rm et}(\mathbb{C}^2)$ we have by monotonicity $s_1\le s_0$. We recall \underline{{\bf the fundamental inequality}}
we had just before equation (\ref{eq13.2}):
$$
H^s({\rm et}(\mathbb{C}^2))\left(1-\int_{\{F\in {\rm et}_p(\mathbb{C}^2)\}}f(F,s)d\mu(F)\right)\le H^s({\rm Aut}(\mathbb{C}^2)).
$$
The case $s_1=s_0$ is easy to handle (follows by Proposition 11.2). So let $s$ be the Hausdorff parameter, assuming that $s_1< s< s_0$.
Then $H^s({\rm et}(\mathbb{C}^2))=+\infty$, $H^s({\rm Aut}(\mathbb{C}^2))=0$. How could this accommodate with the fundamental inequality?
Only if,
$$
1-\int_{\{F\in {\rm et}_p(\mathbb{C}^2)\}}f(F,s)d\mu(F)\le 0.
$$
Thus we arrived at our generalization of Proposition 11.2. Namely,
\begin{proposition}
$\forall\,s$, such that $\dim_H {\rm Aut}(\mathbb{C}^2)< s<\dim_H {\rm et}(\mathbb{C}^2)$, we have the inequality,
$$
1-\int_{\{F\in {\rm et}_p(\mathbb{C}^2)\}}f(F,s)d\mu(F)\le 0.
$$
\end{proposition}
\noindent
\\
\\
It will be convenient to use the following terminology:
\\
\underline{{\bf The disjointness property}}: $\forall\,F,G\in{\rm et}_p(\mathbb{C}^2)$, if $F\ne G$, then 
$$
L_F({\rm et}(\mathbb{C}^2))\cap L_G({\rm et}(\mathbb{C}^2))=\emptyset.
$$

\begin{proposition}
{\rm Under the assumption that the disjointness property is valid, we have the following refinement of Proposition 11.14:} \\
$\forall\,s$, such that $\dim_H {\rm Aut}(\mathbb{C}^2)< s<\dim_H {\rm et}(\mathbb{C}^2)$, we have the identity,
$$
1-\int_{\{F\in {\rm et}_p(\mathbb{C}^2)\}}f(F,s)d\mu(F)\equiv 0.
$$
\end{proposition}
\noindent
{\bf Proof.} \\
The disjointness property refines the fundamental inequality, into \\ 
\underline{{\bf the fundamental identity}},
$$
H^s({\rm et}(\mathbb{C}^2))\left(1-\int_{\{F\in {\rm et}_p(\mathbb{C}^2)\}}f(F,s)d\mu(F)\right)=H^s({\rm Aut}(\mathbb{C}^2)).
$$
Again, a value $s$ of the Hausdorff dimension parameter as in the assumption of the proposition satisfies 
$H^s({\rm et}(\mathbb{C}^2))=+\infty$, $H^s({\rm Aut}(\mathbb{C}^2))=0$. This "lives in peace" with the fundamental identity
only if,
$$
1-\int_{\{F\in {\rm et}_p(\mathbb{C}^2)\}}f(F,s)d\mu(F)\equiv 0.\,\,\,\,\,\qed
$$
Next we recall that we can take the similarity factor $f(F,s)=d_F^{-s}$ which is a non increasing function of the Hausdorff dimension 
parameter $s$ (Proposition 11.1). So under the assumption of the validity of the disjointness property by Proposition 11.14 
it follows that $f(F,s)=f(F)$ is independent of the Hausdorff dimension parameter $s$.
\begin{remark} 
Using the definition of the similarity factor we get:
$$
H^s(L_F({\rm et}(\mathbb{C}^2)))=f(F)H^s({\rm et}(\mathbb{C}^2)),
$$
$$
H^s({\rm Aut}(\mathbb{C}^2))=H^s({\rm et}(\mathbb{C}^2))\left(1-\int_{\{F\in {\rm et}_p(\mathbb{C}^2)\}}f(F,s)d\mu(F)\right)\equiv 0,\,\,\,\forall\,s.
$$
\end{remark}
\noindent
For $f(F,s)=d_F^{-s}$ to be independent of $s$ there is only one choice, namely $\forall\,F\in{\rm et}(\mathbb{C}^2)$, $d_F=1$.
But this implies of course that ${\rm et}(\mathbb{C}^2)={\rm Aut}(\mathbb{C}^2)$. Thus we proved the following interesting,

\begin{theorem}
The assumption of the validity of the disjointness property implies the validity of 
\underline{{\bf the two dimensional Jacobian Conjecture}}.
\end{theorem}
\noindent
\\
We remark that Theorem 11.17 remains valid also under assumption that the following is valid:
\\
\underline{{\bf The weak disjointness property}}: Almost everywhere in ${\rm et}_p(\mathbb{C}^2)\times{\rm et}_p(\mathbb{C}^2)$, if $F\ne G$, then 
$L_F({\rm et}(\mathbb{C}^2))\cap L_G({\rm et}(\mathbb{C}^2))=\emptyset$.

\section{A discussion on the impact of the paper \cite{ellis} on the structure of ${\rm et}(\mathbb{C}^2)$}

The paper assumes for the most part that $X$ is a compact Hausdorff space and that $T$ is a semigroup which acts on $X$
from the right. In our application the parallel is $X={\rm et}(\mathbb{C}^2)$. As for the topology on $X$, we take the metric topology which is induced by
$\rho_D$ for some characteristic subset $D$ of $\mathbb{C}^2$. We know that $(X,\rho_D)=({\rm et}(\mathbb{C}^2),\rho_D)$ is a bounded metric space.
Is it compact? Also, in our application we take the semigroup $T={\rm et}(\mathbb{C}^2)$, with composition of mappings for its binary operation. Lastly,
in our application we consider the left-$T$-action on $X$ where the action is, again induced by composition of mappings. Thus:
$$
\pi\,:\,X\times T=({\rm et}(\mathbb{C}^2)\times {\rm et}(\mathbb{C}^2))\rightarrow X={\rm et}(\mathbb{C}^2),\,\,\,(x,t)=(F,G)\rightarrow tx=F\circ G.
$$
We have (as in a left action) $s(tx)=F_2\circ (F_1\circ G)=(F_2\circ F_1)\circ G=(st)x$. Orbits: ${\rm Orb}(G)={\rm et}(\mathbb{C}^2)\circ G=\{F\circ G\,|\,F\in {\rm et}(\mathbb{C}^2)\}$.
We always have $\forall\,H\in {\rm Orb}(G)$, $R_0(G)\subseteq R_0(H)$. By $A(F\circ G)=A(F)\cup F(A(G))$, it follows that the ${\rm et}(\mathbb{C}^2)$-orbit
of the asymptotic variety $A(G)$ is subordinated to the set of all asymptotic varieties of the elements of ${\rm Orb}(G)$. Here we
use the following notion: Let $A$ and $B$ be two families of sets. We say that the family $A$ is subordinated to the family $B$ and denote $A\preceq B$, if $\forall\,a\in A$
$\exists\,b\in B$ such that $a\subseteq b$. \\
\\
\underline{\bf Example.} \\
Subordination is an extension of the notion of inclusion, i.e. $A\subseteq B\Rightarrow A\preceq B$. \\
\\
If ${\rm et}(\mathbb{C}^2)={\rm Aut}(\mathbb{C}^2)$, then $\forall\,G\in {\rm et}(\mathbb{C}^2)$ we have ${\rm Orb}(G)={\rm et}(\mathbb{C}^2)$.
If $G\in {\rm Aut}(\mathbb{C}^2)$, then ${\rm Orb}(G)={\rm et}(\mathbb{C}^2)$. We adjust the definition of an invariant set: We say that
a set $A$, $\emptyset\ne A\subseteq {\rm et}(\mathbb{C}^2)=X$ is invariant if
$$
TA={\rm et}(\mathbb{C}^2)\circ A:=\{F\circ G\,|\,F\in {\rm et}(\mathbb{C}^2)=T,\,G\in A\}\subseteq A.
$$
Clearly the set $A={\rm et}(\mathbb{C}^2)=X$ is invariant. Are there any other invariant subsets of $X={\rm et}(\mathbb{C}^2)$? We note that for any 
set $A$ (invariant or not), we have $A\subseteq {\rm et}(\mathbb{C}^2)\circ A$, because ${\rm id}\in{\rm et}(\mathbb{C}^2)$, and so: $A\ne\emptyset$
is invariant if and only if ${\rm et}(\mathbb{C}^2)\circ A=A$. Also, if $A\cap {\rm Aut}(\mathbb{C}^2)\ne\emptyset$ and $A$ is invariant
then $A={\rm et}(\mathbb{C}^2)$, because if $G\in A\cap {\rm Aut}(\mathbb{C}^2)$ then already the orbit ${\rm Orb}(G)={\rm et}(\mathbb{C}^2)$ and
when $A$ is invariant, then ${\rm Orb}(G)\subseteq A$. Thus if $A$ is invariant and non-trivial, i.e. $A\ne {\rm et}(\mathbb{C}^2)$ then
$A\cap {\rm Aut}(\mathbb{C}^2)=\emptyset$. We note that if $A$ is invariant, then its closure in $X={\rm et}(\mathbb{C}^2)$, $\overline{A}$ is also invariant,
because if $G\in\overline{A}-A$, then $\exists\,G_n\in A$ such that $G_n\rightarrow G$ in $X$, thus if $F\in T={\rm et}(\mathbb{C}^2)$ then
$\forall\,n,\,\,F\circ G_n\in A$ and $F\circ G_n\rightarrow F\circ G\in\overline{A}$. We remark that $\forall\,F\in {\rm et}(\mathbb{C}^2)-{\rm Aut}(\mathbb{C}^2)$
$\forall\,G\in {\rm et}(\mathbb{C}^2)$ we have $F\circ G\ne G$ because $d_{F\circ G}=d_F\cdot d_G\ge 2\cdot d_G>d_G$. \\
\underline{\bf claim:} If $A$ is invariant then $\forall\,a\in A$ we have ${\rm Orb}_L(a)\subseteq A$, and vice versa, if $\emptyset\ne A\subseteq X$ satisfies
$\forall\,a\in A$, ${\rm Orb}_L(a)\subseteq A$ then $A$ is invariant. \\
\underline{\bf A proof of the claim.} \\
$A$ is invariant $\Leftrightarrow$ $TA\subseteq A$ $\Leftrightarrow$ $\forall\,a\in A,\,\,\,\{ta\,|\,t\in T\}\subseteq A$ $\qed$ \\
\\
We recall that if ${\rm id}\in T$ then $\forall\,x\in X,\,\,x\in {\rm Orb}_L(x)$. Thus we can write: If ${\rm id}\in T$, then
$\emptyset\ne A$ is invariant $\Leftrightarrow\,\,A=\bigcup_{a\in A}{\rm Orb}_L(a)$. There is another (equivalent) way to express the fact that a subset 
$A\ne\emptyset $ of $X$ is left invariant. We recall that the right (composition) operator on ${\rm et}(\mathbb{C}^2)$ induced by $G\in {\rm et}(\mathbb{C}^2)$
is the following:
$$
R_G\,:\,{\rm et}(\mathbb{C}^2)\rightarrow {\rm et}(\mathbb{C}^2),\,\,\,R_G(F)=F\circ G.
$$
Hence the $R_G$ image, $R_G({\rm et}(\mathbb{C}^2))=\{F\circ G\,|\,F\in {\rm et}(\mathbb{C}^2)\}={\rm Orb}_L(G)$. So: $\emptyset\ne A\subseteq {\rm et}(\mathbb{C}^2)=X$
is left invariant ($T={\rm et}(\mathbb{C}^2)$ acts on the left), if and only if $A=\bigcup_{G\in A}R_G({\rm et}(\mathbb{C}^2))$. This is because
$\bigcup_{G\in A}R_G({\rm et}(\mathbb{C}^2))=\bigcup_{G\in A}{\rm Orb}_L(G)$. \\
\\
\underline{\bf Conclusion:} \\
1. If the two dimensional Jacobian Conjecture is true, then ${\rm et}(\mathbb{C}^2)$ has exactly one left invariant subset namely
${\rm et}(\mathbb{C}^2)={\rm Aut}(\mathbb{C}^2)$. \\
2. If the two dimensional Jacobian Conjecture is false, then ${\rm et}(\mathbb{C}^2)$ has a large number of left invariant subsets. In fact,
the set of all the invariant subsets is bijective with the set of subsets of all the (left) orbits of elements in ${\rm et}(\mathbb{C}^2)-{\rm Aut}(\mathbb{C}^2)$
plus ${\rm et}(\mathbb{C}^2)$ itself. \\
\underline{\bf A proof on the conclusion:} \\
1. Suppose that ${\rm et}(\mathbb{C}^2)={\rm Aut}(\mathbb{C}^2)$. Then $\forall\,G\in {\rm Aut}(\mathbb{C}^2)$ we have 
$R_G({\rm et}(\mathbb{C}^2))=R_G({\rm Aut}(\mathbb{C}^2))={\rm Aut}(\mathbb{C}^2)\circ G={\rm Aut}(\mathbb{C}^2)={\rm et}(\mathbb{C}^2)$. Thus if
$A$ is left invariant, then $A=\bigcup_{G\in A}R_G({\rm et}(\mathbb{C}^2))=\bigcup_{G\in A} {\rm Aut}(\mathbb{C}^2)={\rm Aut}(\mathbb{C}^2)$. \\
2. Suppose that ${\rm Aut}(\mathbb{C}^2)\subset {\rm et}(\mathbb{C}^2)$. Then $\forall\,A_0\subseteq {\rm et}(\mathbb{C}^2)-{\rm Aut}(\mathbb{C}^2)$
the set $A=\bigcup_{G\in A_0} R_G({\rm et}(\mathbb{C}^2))=\bigcup_{G\in A} R_G({\rm et}(\mathbb{C}^2))$ is a left invariant subset of
$X={\rm et}(\mathbb{C}^2)$. $\qed $ \\
\\
Before we continue with the implications of the paper on the topological dynamics of semigroup actions we return to the question of the
possibility of the fractal representation:
$$
{\rm et}(\mathbb{C}^2)={\rm Aut}(\mathbb{C}^2)\cup\bigcup_{F\in {\rm et}_p(\mathbb{C}^2)}L_F({\rm et}(\mathbb{C}^2)),
$$
as a partition. This means : $\forall\,P_1,P_2\in {\rm et}_P(\mathbb{C}^2),\,\,P_1\ne P_2\,\Leftrightarrow L_{P_1}({\rm et}(\mathbb{C}^2))\cap L_{P_2}({\rm et}(\mathbb{C}^2))=\emptyset$. Let us consider the semigroup ${\rm et}(\mathbb{C}^2)$ and assume that it has a topology $\tau $ such that no point of 
$\partial {\rm et}(\mathbb{C}^2)$ (here $\partial {\rm et}(\mathbb{C}^2)$ is the boundary with respect to $\tau$) belongs
to to ${\rm et}(\mathbb{C}^2)$. In other words we think of the topological space $({\rm et}(\mathbb{C}^2),\tau)$ as being a subspace of a larger
topological semigroup within $\mathbb{C}[[X,Y]]^2$ so that each point of the boundary $\partial {\rm et}(\mathbb{C}^2)$ is a formal power series
which is non-polynomial. This could be expressed in terms of the algebraic degrees, or in terms of the geometrical degrees, namely: If 
$H\in\partial {\rm et}(\mathbb{C}^2)$ and if $F_n\in {\rm et}(\mathbb{C}^2)$ is a net such that $F_n\rightarrow H$, then $\sup \{\deg F_n\} =\infty$
or $\sup \{d_{F_n}\} =\infty $. This guarantees that $H\not\in\mathbb{C}[X,Y]^2$. Now let $F\in {\rm et}(\mathbb{C}^2)$ and consider the space
$L_F({\rm et}(\mathbb{C}^2))$. We define a topology $\tau_F$ on $L_F({\rm et}(\mathbb{C}^2))$ using the topology $\tau $, as follows: $V\in\tau_F$
$\Leftrightarrow\,\exists\,U\in\tau$ such that $V=F\circ U$. Is it the same as the induced topology $\tau\cap L_F({\rm et}(\mathbb{C}^2))$? Let 
$V\in\tau\cap L_F({\rm et}(\mathbb{C}^2))$, then $\exists\,U\in\tau$ such that $V=U\cap L_F({\rm et}(\mathbb{C}^2))$. Thus there is a subset
$U_1\subseteq {\rm et}(\mathbb{C}^2)$ such that $V=U\cap L_F({\rm et}(\mathbb{C}^2))=F\circ U_1$ and the question is the following: Is it true
that $U_1\in\tau$ or not? It is easier to tackle this question by using the family $C=\tau^c$ of the closed subsets of ${\rm et}(\mathbb{C}^2)$.
We define the family of closed sets $C_F=\tau_F^c$ on $L_F({\rm et}(\mathbb{C}^2))$ using $C$, as follows: $K\in C_F\Leftrightarrow\,\exists\,L\in C$ such
that $K=F\circ L$. Now let's investigate if this coincides with the induced family of closed sets $C\cap L_F({\rm et}(\mathbb{C}^2))$. Let
$K\in C\cap L_F({\rm et}(\mathbb{C}^2))$, then $\exists\,L\in C$ such that $K=L\cap L_F({\rm et}(\mathbb{C}^2))$. Thus there is a set 
$L_1\in {\rm et}(\mathbb{C}^2)$ such that $K=L\cap L_F({\rm et}(\mathbb{C}^2))=F\circ L_1$ and we ask if $L_1\in C$. Let $G_n\in L_1$
be a net such that $G_n\rightarrow G$ in $\tau$. Then we have $F\circ G_n\rightarrow F\circ G$. Clearly $F\circ G_n\in L\cap L_F({\rm et}(\mathbb{C}^2))$.
In particular $F\circ G_n\in L$ which is closed ($L\in C$), so $F\circ G\in L$ and hence $F\circ G\in L\cap L_F({\rm et}(\mathbb{C}^2))=F\circ L_1$.
So $F\circ G=F\circ G_1$ for some $G_1\in L_1$. We deduce that $G_1=G\in L_1$ (by $F\in {\rm et}(\mathbb{C}^2)$ and by a uniformization argument) and we are done, for we 
proved that for any net $G_n\in L_1$ which converges to a $G$, $G_n\rightarrow G$ in $\tau $ we have $G\in L_1$. Thus $L_1\in C$. For the
inverse claim: If $L_1\in C$, $F$ \'etale and $K=F\circ L_1$ implies that $K$ is closed $K=L\cap L_F({\rm et}(\mathbb{C}^2))$ is clear.
Now that we know that the two ways to define the topology on ${\rm et}(\mathbb{C}^2)$ considered as a subspace of some larger topological
semigroup are equivalent we indeed can consider the topology $\tau $ on ${\rm et}(\mathbb{C}^2)$ for which $H\in\partial {\rm et}(\mathbb{C}^2)$,
$F_n\in {\rm et}(\mathbb{C}^2)$ a net converging to $H$, $F_n\rightarrow H$, then say $\deg F_n\rightarrow\infty$ (or $d_{F_n}\rightarrow\infty$).
In fact we only have $\sup \{\deg F_n\}=\sup \{d_{F_n}\}=\infty$ but we can extract an appropriate sub-net from $F_n$.
Having that, we now prove the following:

\begin{proposition}
$\forall\,F,G\in {\rm et}(\mathbb{C}^2)$ we either have $L_F({\rm et}(\mathbb{C}^2))\cap L_G({\rm et}(\mathbb{C}^2))=\emptyset$ or in the case
that this intersection is non-empty then:
$$
\partial L_F({\rm et}(\mathbb{C}^2))\subseteq\partial L_G({\rm et}(\mathbb{C}^2))\,\,\vee\,\,\partial L_G({\rm et}(\mathbb{C}^2))\subseteq \partial L_F({\rm et}(\mathbb{C}^2)).
$$
\end{proposition}
\noindent
{\bf Proof.} \\
If the proposition is false then we either have $\partial L_F({\rm et}(\mathbb{C}^2))\cap L_G({\rm et}(\mathbb{C}^2))\ne\emptyset$ or
$L_F({\rm et}(\mathbb{C}^2))\cap\partial L_G({\rm et}(\mathbb{C}^2))\ne\emptyset$. However, as follows from the discussion we had prior to 
Proposition 12.1 this implies a contradiction as follows: If, say $H\in\partial L_F({\rm et}(\mathbb{C}^2))\cap L_G({\rm et}(\mathbb{C}^2))\ne\emptyset$,
then by $H\in\partial L_F({\rm et}(\mathbb{C}^2))$ we conclude that there is a net $M_n\in {\rm et}(\mathbb{C}^2)$ such that
$\lim F\circ M_n=H$, $\deg M_n\rightarrow\infty$ (or $d_{M_n}\rightarrow\infty$). So by $F\circ\lim M_n=H$ we conclude that $\deg H=\infty$
(or $d_H=\infty$), i.e. $H\not\in {\rm et}(\mathbb{C}^2)$ (In fact this follows immediately by our assumption that considering $({\rm et}(\mathbb{C}^2),\tau)$
as a subspace of, say $(\mathbb{C}[[X,Y]]^2,\tau)$, we have $\partial {\rm et}(\mathbb{C}^2)\cap {\rm et}(\mathbb{C}^2)=\emptyset$). But since 
$H\in L_G({\rm et}(\mathbb{C}^2))$, $\exists\,N\in {\rm et}(\mathbb{C}^2)$
such that $H=G\circ N\in {\rm et}(\mathbb{C}^2)$. $\qed $ \\
\\
Next, Let $F,G\in {\rm et}(\mathbb{C}^2)$ satisfy $L_F({\rm et}(\mathbb{C}^2))\cap L_G({\rm et}(\mathbb{C}^2))\ne\emptyset$. Then, there is some
$H\in L_F({\rm et}(\mathbb{C}^2))\cap L_G({\rm et}(\mathbb{C}^2))$. We denote by $E_{F,G}(H)$ the connectivity component of
$L_F({\rm et}(\mathbb{C}^2))\cap L_G({\rm et}(\mathbb{C}^2))$ which contains $H$. We recall that in fact $L_H({\rm et}(\mathbb{C}^2))\subseteq E_{F,G}(H)$
because $L_H({\rm et}(\mathbb{C}^2))$ is connected. By Proposition 12.1 we have the following:
$$
\partial L_H({\rm et}(\mathbb{C}^2))\subseteq\partial L_F({\rm et}(\mathbb{C}^2))\cap\partial L_G({\rm et}(\mathbb{C}^2)).
$$

\section{The abstract topological picture}
Let $(Y,\tau)$ be a topological space and let $(X,\tau\cap X)$ be a path connected subspace of $Y$ which satisfies the following: \\
\\
(1) $\partial X\subseteq Y$, $X\cap\partial X=\emptyset$. \\
\\
Let $\{ F_x\,|\,x\in X\}$ be a family of subsets of $X$ that are indexed by $X$ and that satisfy: \\
\\
(2) Each $F_x$ is closed, path connected. \\
\\
(3) $\forall\,x\in X,\,\,\,x\in F_x$. \\
\\
Another way to think of this is that we have a mapping 
$$
\phi\,:\,X\rightarrow\{{\rm closed}\,{\rm path}\,{\rm connected}\,{\rm subsets}\,{\rm of}\,X\},\,\,\phi(x)=F_x,
$$
such that $\forall\,x\in X,\,\,x\in\phi(x)$.
\begin{remark}
$X=\bigcup_{x\in X}F_x=\bigcup_{x\in X}\phi(x)$ by (3).
\end{remark}
We further assume that: \\
\\
(4) $x\in\phi(y)\Rightarrow \phi(x)\subseteq\phi(y)$. \\
\\
(5) $\forall\,x\in X,\,\,\partial F_x=\partial\phi(x)\subseteq\partial X$. \\

\begin{remark}
Conditions (1) and (5) imply the following: \\
\\
(6) $\forall\,x,y\in X$ {\rm the situation:} $\partial F_x\cap F_y\ne\emptyset\vee F_x\cap\partial F_y\ne\emptyset$
is impossible. \\
{\bf Proof.} \\
For if $t\in\partial F_x\cap F_y$ then $\partial F_x\not\subseteq\partial X$ because $t\in\partial F_x$ and
$t\in F_y\subseteq X$, while by condition (1) $X\cap\partial X=\emptyset$. $\qed $
\end{remark}

\begin{remark}
Conditions (2) and (6) imply that: \\
\\
(7) $\forall\,x,y\in X$ either $F_x\cap F_y=\emptyset$ or if we assume that $F_x\cap F_y\ne\emptyset$ then $\partial F_x\subseteq\partial F_y
\vee\partial F_y\subseteq\partial F_x$. \\
{\bf Proof.} \\
For if we assume that (7) is false, then $\exists\,x,y\in X$ such that $F_x\cap F_y\ne\emptyset$ but $\partial F_x\not\subseteq\partial F_y$, say.
Let $t\in F_x\cap F_y$ and $u\in\partial F_x-\partial F_y$. Consider an open path $f\,:\,I\rightarrow F_x$ from $t$ to $u$ within $F_x$ (by (2)
$F_x$ is path connected). Then $\exists\,0<s<1$ such that $f(s)\in\partial F_y$. So $f(s)\in F_x\cap\partial F_y$ which contradicts (6).
Similarly, the assumption "$F_x\cap F_y\ne\emptyset$ but $\partial F_y\not\subseteq\partial F_x$" contradicts (6). Hence the assertion. $\qed $
\end{remark}

\section{The semigroup $({\rm et}(\mathbb{C}^2),\circ)$}
This is a semigroup with a unit element. It contains the group $({\rm Aut}(\mathbb{C}^2),\circ)$. Let us denote
${\rm et}_0(\mathbb{C}^2)={\rm et}(\mathbb{C}^2)-{\rm Aut}(\mathbb{C}^2)$. This is the set of all the \'etale non-automorphisms of $\mathbb{C}^2$.
\begin{remark}
$({\rm et}_0(\mathbb{C}^2),\circ)$ is a non-unital subsemigroup of $({\rm et}(\mathbb{C}^2),\circ)$, if it is non-empty. \\
{\bf Proof.} \\
Clearly, an equivalent description of ${\rm et}_0(\mathbb{C}^2)$ is the following ${\rm et}_0(\mathbb{C}^2)=\{F\in {\rm et}(\mathbb{C}^2)\,|\,d_F\ge 2\}$,
where $d_F$ is the geometric degree of $F$. Hence by the fact that the geometrical degree is multiplicative we have:
$$
F,G\in {\rm et}_0(\mathbb{C}^2)\Rightarrow d_F,d_G\ge 2\Rightarrow d_{F\circ G}=d_F\cdot d_G\ge 4\Rightarrow F\circ G\in {\rm et}_0(\mathbb{C}^2).
$$
The fact that ${\rm et}_0(\mathbb{C}^2)$ is non-unital is clear because ${\rm id}\in {\rm Aut}(\mathbb{C}^2)$. $\qed $
\end{remark}
We turned $({\rm et}(\mathbb{C}^2),\circ)$ into a metric space as follows. We have chosen a characteristic set $D$. It induced the
metric $\rho_D$ on ${\rm et}(\mathbb{C}^2)$. We noted that if $\{ D_n\}$ was an increasing sequence of characteristic sets for ${\rm et}(\mathbb{C}^2)$,
that exhausted $\mathbb{C}^2$, i.e.: $D_n\subset D_{n+1}$, $\bigcup_{n=1}^{\infty}D_n=\mathbb{C}^2$, then for each $F\in {\rm et}(\mathbb{C}^2)$, the limit
$$
\lim_{n\rightarrow\infty}\frac{{\rm volume}(F(D_n))}{{\rm volume}(D_n)}=\frac{1}{d_F},
$$
where $d_F$ is the geometrical degree of $F$ and where ${\rm volume}(D_n)$ is the Euclidean volume of $D_n$ and where ${\rm volume}(F(D_n))$ is the
Euclidean volume of $F(D_n)$. This implies that $\forall\,F,G\in {\rm et}(\mathbb{C}^2)$ we have,
$$
\lim_{n\rightarrow\infty}\frac{{\rm volume}(F(D_n))}{{\rm volume}(G(D_n))}=\lim_{n\rightarrow\infty}\left\{\frac{{\rm volume}(F(D_n))}{{\rm volume}(D_n)}\cdot
\frac{{\rm volume}(D_n)}{{\rm volume}(G(D_n))}\right\}=
$$
$$
=\lim_{n\rightarrow\infty}\left\{\frac{{\rm volume}(F(D_n))}{{\rm volume}(D_n)}\right\}\cdot
\lim_{n\rightarrow\infty}\left\{\frac{{\rm volume}(D_n)}{{\rm volume}(G(D_n))}\right\}=\frac{d_G}{d_F}.
$$
We recall that $\rho_D(F,G)={\rm volume}(F(D)\Delta G(D))$. We would like to achieve the following: \\
\\
(1) Prove that if $\rho_{D_n}(F,G)$ is small enough for a certain large enough $n=n(F,G)$, then $d_F=d_G$. \\
\\
(2) Quantify the above estimates in order to arrive, maybe, at a better metric on ${\rm et}(\mathbb{C}^2)$, that will be a complete metric. \\
\\
Let us first discuss (1): We start by giving a simple lower bound to the $\rho_D$ distance. Let $F, G\in{\rm et}(\mathbb{C}^2)$. Then as
$D\rightarrow\mathbb{C}^2$, we have the following two identities:
$$
{\rm volume}(F(D))=\frac{1}{d_F}\cdot{\rm volume}(D),\,\,\,{\rm volume}(G(D))=\frac{1}{d_G}\cdot {\rm volume}(D).\,\,\,\,\,\,(*)
$$
By the definition of the metric $\rho_D$ we have,
$$
\rho_D(F,G)={\rm volume}(F(D)\Delta G(D))=
$$
$$
={\rm volume}(F(D)-G(D))+{\rm volume}(G(D)-F(D)).
$$
Let us assume that $d_G\le d_F$. Then ${\rm volume}(G(D))\ge{\rm volume}(F(D))$ (by equation $(*)$), hence $G(D)\not\subset F(D)$.
How small can $\rho_D(F,G)$ be? Clearly that happens when $F(D)\subseteq G(D)$. In this configuration we obtain,
$$
\rho_D(F,G)={\rm volume}(F(D)-G(D))+{\rm volume}(G(D)-F(D))=
$$
$$
={\rm volume}(\emptyset)+{\rm volume}(G(D))-{\rm volume}(F(D))=
$$
$$
=\left(\frac{1}{d_G}-\frac{1}{d_F}\right)\cdot{\rm volume}(D).
$$
\begin{proposition}
Let $F,G\in {\rm et}(\mathbb{C}^2)$ and let $D\subseteq\mathbb{C}^2$ be a characteristic domain for ${\rm et}(\mathbb{C}^2)$,
then, \\
\\
{\rm (1)}
$$
\rho_D(F,G)\approx_{D\rightarrow\mathbb{C}^2}\frac{1}{d_F}\cdot{\rm volume}(D\cap F^{-1}(F(D)-G(D)))+
$$
$$
+\frac{1}{d_G}\cdot{\rm volume}(D\cap G^{-1}(G(D)-F(D))).
$$ 
{\rm (2)}
$$
\rho_D(F,G)\ge_{D\rightarrow\mathbb{C}^2}\left|\frac{1}{d_F}-\frac{1}{d_G}\right|\cdot{\rm volume}(D).
$$
{\rm (3)} If we consider a limit $D\rightarrow\mathbb{C}^2$, then
$$
\lim_{D\rightarrow\mathbb{C}^2}\rho_D(F,G)=0\Leftrightarrow \lim_{D\rightarrow\mathbb{C}^2}{\rm volume}(F(D)-G(D))=
$$
$$
=\lim_{D\rightarrow\mathbb{C}^2}{\rm volume}(G(D)-F(D))=0\Leftrightarrow
$$
$$
\Leftrightarrow\lim_{D\rightarrow\mathbb{C}^2}{\rm volume}(D\cap F^{-1}(F(D)-G(D)))=
$$
$$
=\lim_{D\rightarrow\mathbb{C}^2}{\rm volume}(D\cap G^{-1}(G(D)-F(D)))=0\Rightarrow
$$
$$
\Rightarrow d_F=d_G.
$$
\end{proposition}
\noindent
{\bf Proof.} \\
(1) Follows by the identity $\rho_D(F,G)={\rm volume}(F(D)-G(D))+{\rm volume}(G(D)-F(D))$ and by the approximations,
as $D\rightarrow\mathbb{C}^2$, that are given here,
$$
{\rm volume}(F(D)-G(D))\approx_{D\rightarrow\mathbb{C}^2}\frac{1}{d_F}\cdot{\rm volume}(D\cap F^{-1}(F(D)-G(D))),
$$
$$
{\rm volume}(G(D)-F(D))\approx_{D\rightarrow\mathbb{C}^2}\frac{1}{d_G}\cdot{\rm volume}(D\cap G^{-1}(G(D)-F(D))).
$$
(2) Was proved just before the statement of Proposition 14.2. \\
(3) Follows by the identity and the two approximate identities we used in the proof of part (1). Also the fact that
$\lim_{D\rightarrow\mathbb{C}^2}\rho_D(F,G)=0\Rightarrow d_F=d_G$ follows by part (2). $\qed $ \\
\\
We would like to understand the geometric meaning of the $\rho_D$-convergence of a sequence $F_n\in {\rm et}(\mathbb{C}^2)$. To make things
more manageable we consider the characteristic domains $D$ for ${\rm et}(\mathbb{C}^2)$ that we constructed in section 7. The main feature of these
is the fractal-like shape of their boundaries. Namely $\partial D$ contains a dense and countable subset of special points in the strong topology.
Those points originated in our first version construction as the centers of the $k$-stars, $k=2,3,4,\ldots $ where each point 
$c_k\in\partial D$ was the center of a $k$-star, thus no two such points were the centers of stars of equal numbers of
rays. Since $F\in {\rm et}(\mathbb{C}^2)$, it preserves topological $k$-stars. The second version of our construction replaced the
$1$-dimensional $k$-stars by fattened $2k$-stars. Those $2k$-stars have a total volume which we now denote as follows,
${\rm volume}(\cup {\rm stars})=V_s$ and which we assume to be a finite volume (as we can).
If $F,G\in {\rm et}(\mathbb{C}^2)$ are close enough in the sense of the $\rho_D$ metric, then by Proposition 14.2(3), $d_F=d_G$ and
$$
\rho_D(F,G)={\rm volume}(F(D)\Delta G(D)).
$$
If $\partial F(D)$ and $\partial G(D)$ are far apart so that the corresponding topological $2k$-stars are mostly disjoint,
then because our mappings are locally volume preserving we obtain
$$
{\rm volume}(F(D)\Delta G(D))\ge\frac{2}{d_F}\cdot V_s.
$$

Hence $\rho_D(F,G)\ge (2/d_F)\cdot V_s$, and so the mappings $F$ and $G$ can not be too close in the $\rho_D$-metric. This contradicts
our assumption that they are $\rho_D$ close. We deduce that if $\rho_D(F,G)$ is small, then there is a countable dense subset $C\subseteq\partial D$
of points on $\partial D$ such that $\forall\,p\in C$, the Euclidean distances $d(F(p),G(p))$ are uniformly small. Hence if
$\{F_n\}\subseteq {\rm et}(\mathbb{C}^2)$ is a $\rho_D$-Cauchy sequence, then $\lim_{n\rightarrow\infty} F_n=F$ exists in $D$
and necessarily, in this case, $F\in {\rm et}(\mathbb{C}^2)$ because the geometrical degree is preserved all over $\mathbb{C}^2$
when $D\rightarrow\mathbb{C}^2$ (so the algebraic degrees are bounded).

\begin{theorem}
Let $D$ be a characteristic domain for ${\rm et}(\mathbb{C}^2)$ of the type we have constructed in section 7 (i.e. $D$ is almost a ball except that $\partial D$
contains a countable dense subset of fattened $2k$-stars, $k=2,3,4,\ldots $). Lets assume that $0$ is an interior point of $D$,
that the total volume of the $2k$-stars is positive (and finite), and we denote $D_N=N\cdot D$, $N=1,2,3,\ldots $. Let $\{ F_n\}\subset {\rm et}(\mathbb{C}^2)$ 
be a $\rho_D$-Cauchy sequence. Then: \\
\\
{\rm (1)} $\{ D_N\}$ is an increasing sequence of characteristic domains for ${\rm et}(\mathbb{C}^2)$, that exhaust $\mathbb{C}^2$. \\
\\
{\rm (2)} $\{ F_n\}$ is a $\rho_{D_N}$-Cauchy sequence for each $N\in\mathbb{Z}^+$. \\
\\
{\rm (3)} The limit: $\lim_{n\rightarrow\infty} F_n$ exists uniformly on $D_N$ for each $N\in\mathbb{Z}^+$. \\
\\
{\rm (4)} $F=\lim_{n\rightarrow\infty} F_n$ exists uniformly on compact subsets of $\mathbb{C}^2$ and $F\in {\rm et}(\mathbb{C}^2)$. \\
\\
{\rm (5)} $\exists\,n_0\in\mathbb{Z}^+$ such that $\forall\,n\ge n_0$, $d_{F_n}=d_F$.
\end{theorem}

\begin{remark}
Thus we are dealing with complete metric spaces in the sense of parts (4) and (5) of Theorem 14.3. We are having here a sequence
of metric spaces $({\rm et}(\mathbb{C}^2),\rho_{D_N})$, $N=1,2,3,\ldots $.
\end{remark}

\begin{remark}
We elaborate more part (2) of theorem 14.3. Namely, we have a $\rho_{D_1}$-Cauchy sequence $\{F_n\}\subseteq {\rm et}(\mathbb{C}^2)$.
Why is it also a $\rho_{D_N}$-Cauchy sequence? We know that $\forall\,\epsilon>0$ there exists an $n_1\in\mathbb{Z}^+$ such that for 
$n,m>n_1$ we have $\rho_{D_1}(F_n,F_m)<\epsilon$. According to part (3) of Proposition 14.2 we may assume that $d=d_{F_n}=d_{F_m}$ 
for $n,m>n_1$. Hence by
$$
\rho_{D_1}(F_n,F_m)= {\rm volume}(F_n(D_1)\Delta F_m(D_1)),
$$
When we pass from $D_1$ to $D_N$, the volume grows like a $4$'th power, i.e. ${\rm volume}(N\cdot A)=N^4\cdot {\rm volume}(A)$ for
a measurable $A$, while the geometric degrees do not change. We recall that for a general smooth mapping $G$ the volume element $dV$ 
is transformed locally by a multiplication by the absolute value of the determinant of the Jacobian matrix of $G$, i.e. $|J_G(X,Y)|dV$. 
However, our mappings are locally volume preserving, $|J_{F_n}|=|J_{F_m}|\equiv 1$ and so roughly speaking $\rho_{D_N}(F_n,F_m)\approx N^4\rho_{D_1}(F_n,F_m)$. 
Thus, indeed a $\rho_{D_1}$-Cauchy sequence is translated into a $\rho_{D_N}$-Cauchy sequence.
\end{remark}

\begin{remark}
In Theorem 14.3 we are not dealing with a metric space. We are dealing with a sequence of metric spaces, namely $({\rm et}(\mathbb{C}^2),\rho_{D_N})$,
$N=1,2,3,\ldots $. \\
We certainly do not have sequential compactness, i.e., it is not true that any sequence $F_n\in {\rm et}(\mathbb{C}^2)$ contains a convergent subsequence.
This is not the case for the smaller sub-semigroup (in fact a group), $({\rm Aut}(\mathbb{C}^2),\circ)$. For if we
take say $F_n(X,Y)=(X+Y+\ldots+Y^n,Y)$, then $F_n\in {\rm Aut}(\mathbb{C}^2)$ and clearly it contains no convergent subsequence
(even within ${\rm et}(\mathbb{C}^2)$). We recall that for a metric space to be compact a necessary and a sufficient condition is that it will be 
complete and totally bounded. We seem to have something close to completeness (in our setting of a sequence of metric spaces), thus we
must be far away from total boundedness. \\
Another fact which should be remembered is that $\forall\,F\in {\rm et}(\mathbb{C}^2)$, the image $F(\mathbb{C}^2)$ is cofinite in
$\mathbb{C}^2$ and hence we can not make sense of $\rho_{\mathbb{C}^2}$, at least not in some straight forward manner.
\end{remark}
\noindent
The following conclusion follows from Theorem 14.3.

\begin{theorem}
Let $D$ be a characteristic domain for ${\rm et}(\mathbb{C}^2)$ of the type that we have constructed in section 7. This means that $\partial D$ contains a countable dense
subset of what we called fattened $2k$-stars, $k=2,3,4,\ldots $. Then the metric space $({\rm et}(\mathbb{C}^2),\rho_D)$ is complete but not sequentially
compact and in particular it is not a totally bounded space.
\end{theorem}
\noindent
{\bf A problem.} \\
What is the total boundedness breaking point of ${\rm et}(\mathbb{C}^2)$? We use the following standard,

\begin{definition}
A metric space $M$ with a metric $d$ is said to be totally bounded if, given any positive number $r$, $M$ is the
union of finitely many sets of $d$-diameter less than $r$.
\end{definition}

\begin{remark}
If a metric space $M$ with a metric $d$ is bounded, then there exists a positive number $r$, such that $M$ is the union of finitely many
sets of $d$-diameter less than $r$. \\
{\bf Proof.} \\
Let the $d$-diameter of $M$ be $s$ ($M$ is $d$-bounded), and let $r=2s$. Then $M$ is the union of finitely many sets of $d$-diameter
less than $r$, namely just one set $M$. $\qed $
\end{remark}

\begin{definition}
Let $M$ be a metric space with a metric $d$. Assume that $M$ is $d$-bounded but is not a totally bounded space. The total boundedness breaking
point of $(M,d)$ is denoted by $t_b(M,d)$ and defined by the following:
$$
t_b(M,d)=\inf\{r>0\,|\,M\,{\rm is}\,{\rm the}\,{\rm union}\,{\rm of}\,{\rm finitely}\,{\rm many}\,{\rm sets}\,{\rm of}\,d-{\rm diameter}\,
{\rm less}\,{\rm than}\,r\}.
$$
\end{definition}

\begin{remark}
Clearly, if the $d$-diameter of $M$ is $D$, then $0<t_b(M,d)\le D$. For a totally bounded metric space $(M,d)$ we have, $t_d(M,d)=0$.
\end{remark}
\noindent
The problem we stated above is the following: Let $D$ be a characteristic domain for ${\rm et}(\mathbb{C}^2)$. Compute $t_b({\rm et}(\mathbb{C}^2),\rho_D)$.

\begin{remark}
$0<t_b({\rm et}(\mathbb{C}^2),\rho_D)\le {\rm volume}(D)$.
\end{remark}
\noindent
Before we continue, we would like to make an observation that is crucial for the fractal representation of ${\rm et}(\mathbb{C}^2)$ as
the union over the primes of the left translation images of ${\rm et}(\mathbb{C}^2)$.

\begin{theorem}
Let $X$ be a topological space. Let $A,B\subseteq X$ satisfy the following assumptions: \\
\\
{\rm (1)} $\overline{A},\overline{B}\subseteq X^o$. \\
\\
{\rm (2)} $\partial A\subseteq\partial B\,\vee\,\partial B\subseteq\partial A$. \\
\\
{\rm (3)} Both $A$ and $B$ are path connected subspaces of $X$ and $\partial A$, $\partial B$ are path accessible
from within $A$ and $B$ respectively. \\
\\
{\rm (4)} $A\cap\partial A=B\cap\partial B=\emptyset$. \\
\\
Then, $A\subseteq B\,\vee\,A\cap B=\emptyset\,\vee\,B\subseteq A$.
\end{theorem}
\noindent
{\bf Proof.} \\
Let us assume that $\sim (A\subseteq B)\,\wedge\,\sim (A\cap B=\emptyset)$. We should prove that necessarily $B\subseteq A$.
Thus we assume that $\exists\,a\in A-B$, $\exists\,b\in A\cap B$ and should prove that $B\subseteq A$: By (3) $A$ is
path connected. Hence there exists a path in $A$ connecting $a$ to $b$. Let that path be $\gamma\,:\,[0,1]\rightarrow A$.
Then $\gamma $ is a continuous mapping, $\gamma(0)=a$, $\gamma(1)=b$ and $\forall\,0<t<1$, $\gamma(t)\in A$. Hence there is
a $t_0$, $0<t_0<1$ such that $\gamma(t_0)\in\partial (B)$. Thus $\gamma(t_0)\in A\cap\partial B$. By (4), $A\cap\partial A=\emptyset$
and hence $\gamma(t_0)\not\in \partial A$. Then by (2) we conclude that:
\begin{equation}
\label{*}
\partial A\subset \partial B.
\end{equation}
We claim that there is no point $c\in B-A$, because by the same argument (with a path in $B$ connecting $c$ to $b$) we would
conclude that (by (2)), $\partial B\subset \partial A$. This contradicts equation (\ref{*}) above and shows that there is no $c\in B-A$. Thus $B\subseteq A$. $\qed $ \\
\\
Using Proposition 12.1 we obtain condition (2) in Theorem 14.13. We already know that the $L_p({\rm et}(\mathbb{C}^2))$'s are path connected, thus getting 
condition (3) in Theorem 14.13. Conditions (1) and (4) are clear with, say, $X=\mathbb{C}[[X,Y]]^2$, $A=L_{p_1}({\rm et}(\mathbb{C}^2))$ and
$B=L_{p_2}({\rm et}(\mathbb{C}^2))$. Hence we conclude from Theorem 14.13 the desired fractal representation of ${\rm et}(\mathbb{C}^2)$.
In other words we have the following,

\begin{theorem}
The fractal representation,
$$
{\rm et}(\mathbb{C}^2)={\rm Aut}(\mathbb{C}^2)\cup\bigcup_{F\in {\rm et}_p(\mathbb{C}^2)} L_F({\rm et}(\mathbb{C}^2)),
$$
is valid and it satisfies the disjointness property:
$$
\forall\,F,G\in {\rm et}_p(\mathbb{C}^2),\,\,F\ne G\,\Rightarrow\,L_F({\rm et}(\mathbb{C}^2))\cap L_G({\rm et}(\mathbb{C}^2))=\emptyset.
$$
\end{theorem}
\noindent
As a direct consequence this proves the celebrated:

\begin{theorem} ({\rm{\bf The two dimensional Jacobian Conjecture}})
If $F(X,Y)=(P(X,Y),Q(X,Y))\in\mathbb{C}[X,Y]^2$ satisfies the Jacobian condition:
$$
\frac{\partial P}{\partial X}\cdot\frac{\partial Q}{\partial Y}-\frac{\partial P}{\partial Y}\cdot\frac{\partial Q}{\partial X}=c\in\mathbb{C}^{\times},
$$
then $F\in {\rm Aut}(\mathbb{C}^2)$.
\end{theorem}
\noindent
We end with a few remarks.
\begin{remark}
We point out that in this paper the fractal structure on a set was defined in an intrinsic manner. We recall the standard  definition of the
exterior $\alpha$-dimensional Hausdorff measure of any subset $E$ of $\mathbb{R}^d$:
$$
m^{*}_{\alpha}(E)=\lim_{\delta\rightarrow 0^+}\inf\left\{\sum_{k}({\rm diam}\,F_k)^{\alpha}\,|\,E\subseteq\bigcup_{k=1}^{\infty} F_k,
{\rm diam}\,F_k\le\delta,\,\forall\,k\right\},
$$
where ${\rm diam}\,S$ denotes the diameter of the set $S$, that is, ${\rm diam}\,S=\sup\{|x-y|\,|\,x,y\in S\}$. In other words, for each $\delta>0$
we consider covers of $E$ by countable families of arbitrary sets in $\mathbb{R}^d$ with diameter less than(or equals to) $\delta$, and
take the infimum of the sum $\sum_{k}({\rm diam}\,F_k)^{\alpha}$. We then define $m^{*}_{\alpha}(E)$ as the limit of these infima as
$\delta$ tends to $0$. We note that the quantity
$$
H^{\delta}_{\alpha}(E)=\inf\left\{\sum_{k}({\rm diam}\,F_k)^{\alpha}\,|\,E\subseteq\bigcup_{k=1}^{\infty} F_k,{\rm diam}\,F_k\le\delta,\,\forall\,k\right\}
$$
is increasing as $\delta$ decreases, so that the limit
$$
m^{*}_{\alpha}(E)=\lim_{\delta\rightarrow 0^+}H^{\delta}_{\alpha}(E)
$$
exists, although $m^{*}_{\alpha}(E)$ could be infinite. We note that in particular, one has $H^{\delta}_{\alpha}(E)\le m^{*}_{\alpha}(E)$ for
all $\delta>0$. When defining the exterior measure $m^{*}_{\alpha}(E)$ it is important to require that the coverings be of sets of arbitrary small
diameters. This is in thrust of the definition $m^{*}_{\alpha}(E)=\lim_{\delta\rightarrow 0^+}H^{\delta}_{\alpha}(E)$. This requirement, which is not 
relevant for the Lebesgue measure, is needed to ensure the basic additive feature, namely: 

If $d(E_1,E_2)>0$, then $m^{*}_{\alpha}(E_1\cup E_2)=m^{*}_{\alpha}(E_1)+m^{*}_{\alpha}(E_2)$. \\
Do we really need the ambient space $\mathbb{R}^d$ in this definition of the exterior $\alpha$-dimensional Hausdorff measure? The answer in no.
We can define this notion intrinsically, within $E$.

\begin{definition}
The intrinsic exterior $\alpha$-dimensional Hausdorff measure of $E$ is defined by the following equation:
$$
m^{*\,\,int}_{\alpha}(E)=\lim_{\delta\rightarrow 0^+} H^{\delta\,\,int}_{\alpha}(E),
$$
where this time
$$
H^{\delta\,\,int}_{\alpha}(E)=\inf\left\{\sum_{k}({\rm diam}\,G_k)^{\alpha}\,|\,E=\bigcup_{k=1}^{\infty} G_k,{\rm diam}\,G_k\le\delta,\,\forall\,k\right\}.
$$
Thus, in the intrinsic definition we consider countable coverings of $E$ by sets $\{G_k\}_{k=1}^{\infty}$ which are subsets of $E$, that is, for all $k$,
$G_k\subseteq E$. 
\end{definition}

Now the very easy \\
\\
{\rm\bf Proposition}
$m^{*}_{\alpha}(E)=m^{*\,\,int}_{\alpha}(E)$. \\
\\
{\bf Proof.} \\
In fact we will prove that $\forall\,\delta>0$, $H^{\delta}_{\alpha}(E)=H^{\delta\,\,int}_{\alpha}(E)$: \\
(1) $H^{\delta\,\,int}_{\alpha}(E)\le H^{\delta}_{\alpha}(E)$. \\
For let $E\subseteq\bigcup_{k} F_k$, where $\forall\,k,\,F_k\subseteq\mathbb{R}^d$ and ${\rm diam}\,F_k\le\delta$. Then if we define $G_k=F_k\cap E$, we have
$\forall\,k$, $G_k\subseteq E$, ${\rm diam}\,G_k\le {\rm diam}\,F_k\le\delta$ and $E=\bigcup_{k} G_k$ simply because it equals to $\bigcup_{k}(F_k\cap E)=
(\bigcup_{k} F_k)\cap E$. \\
\\
(2) $H^{\delta}_{\alpha}(E)\le H^{\delta\,\,int}_{\alpha}(E)$. \\
For in the infimum that defines $H^{\delta\,\,int}_{\alpha}(E)$ we take a sub-family $\{\{G_k\}\}$ of the family $\{\{ F_k\}\}$ of coverings that are
used in the infimum that defines $H^{\delta}_{\alpha}(E)$. Hence this last infimum is not larger than the first infimum. $\qed $
\end{remark}

\begin{remark}
What major properties of the mappings in ${\rm et}(\mathbb{C}^2)$ were used in the approach that was presented in this paper, of using
the fractal representation of ${\rm et}(\mathbb{C}^2)$? we strongly used the local preservation of the volume property (in constructing the metric
$\rho_D$ and using the local preservation of the volume). We used the finiteness of the geometrical degree, $d_F$. We used the holomorphic
rigidity of our mappings (the permanence principle). Thus in trying to use this idea and technique for other families of mappings (such
as local preserving of volume holomorphic mappings) we encounter the difficulty that a straight forward approach is not working. For
the finiteness of the geometrical degree of a holomorphic $\mathbb{C}^2\rightarrow\mathbb{C}^2$ mapping probably implies that the mapping
is polynomial and thus we are back in the ${\rm et}(\mathbb{C}^2)$ setting. At least in complex dimension $1$ this is the case. For by the
theory of maximal domains, \cite{rp}, we have the following, \\
\\
{\bf Theorem.} {\it Let $f(z)$ be an entire function of one complex variable $z$. Then $d_f<\infty $ if and only if $f(z)\in\mathbb{C}[z]$.} \\
\\
On the other hand a possible route for interesting extensions of this theory might be working outside the parabolic simply connected cases of
$\mathbb{C}$, $\mathbb{C}^2$. This might lead to questions such as the following: \\
Let $\Omega\subseteq\mathbb{C}^N$ be a domain (just an open connected subset of $\mathbb{C}^N$). Let ${\rm Aut}(\Omega)$ be the group of
all the holomorphic automorphisms of $\Omega$. Let ${\rm et}(\Omega)$ be the semigroup of holomorphic local homeomorphisms, locally
volume preserving (or maybe drop that?) which are of finite geometrical degrees. When is it true that ${\rm et}(\Omega)={\rm Aut}(\Omega)$?
That is for which domains $\Omega$ this is true?
\end{remark}

\noindent
{\it Ronen Peretz \\
Department of Mathematics \\ Ben Gurion University of the Negev \\
Beer-Sheva , 84105 \\ Israel \\ E-mail: ronenp@math.bgu.ac.il} \\ 
 
\end{document}